\documentclass[nonblindrev]{informs3noheader}

\OneAndAHalfSpacedXI

\usepackage{mathtools}
\usepackage[inline]{enumitem}
\usepackage{subfig}
\usepackage{booktabs}
\usepackage{array}
\usepackage{multirow}
\usepackage{booktabs,caption,fixltx2e}
\usepackage[flushleft]{threeparttable}
\usepackage{lscape}
\usepackage{makecell}
\usepackage{amssymb}
\usepackage{amsmath}
\usepackage[mathscr]{euscript}
\usepackage{algorithm}
\usepackage{algpseudocode}
\usepackage{pifont}
\usepackage{caption}
\usepackage{algorithmicx}
\usepackage{hyperref}
\usepackage{mathrsfs}
\usepackage{arydshln}
\usepackage{tabularx}
\usepackage{adjustbox}
\usepackage{bm, bbm}
\usepackage{soul}
\usepackage{blindtext}
\usepackage{longtable}
\hypersetup{
colorlinks=true,
linkcolor=black,
citecolor=black,
filecolor=black,
urlcolor=black,
}

\usepackage{comment}

\usepackage{tikz}
\usetikzlibrary{shapes.geometric,shapes.multipart}
\usetikzlibrary{positioning}
\usetikzlibrary{decorations.pathreplacing}
\definecolor{mygreen}    {RGB}{0,90,0}
\definecolor{myblue}     {RGB}{0,51,140}
\definecolor{myorange}   {RGB}{238,118,0}
\definecolor{myred}      {RGB}{126,0,0}
\definecolor{mygray}     {RGB}{100,100,105}
\definecolor{mygrayblue} {RGB}{0,128,128}
\definecolor{mygraygreen}{RGB}{128,128,0}
\definecolor{DarkPurple}     {RGB}{142, 36, 170}
\definecolor{LightPurple}    {RGB}{57, 130, 7}

\usepackage{natbib}
\bibpunct[, ]{(}{)}{,}{a}{}{,}%
\def\bibfont{\small}%

\usepackage[utf8]{inputenc}
\usepackage[english]{babel}
\usepackage{graphicx}
\newcommand{\md}[1]{\mathbb{#1}}

\newcommand{\Z}{\mathcal{Z}}



\newcolumntype{F}[1]{%
    >{\raggedright\arraybackslash\hspace{0pt}}p{#1}}%
\newcolumntype{T}[1]{%
    >{\centering\arraybackslash\hspace{0pt}}p{#1}}%

\TheoremsNumberedThrough
\EquationsNumberedThrough

\def\R{\mathbb{R}}
\def\Z{\mathbb{Z}}

\def\calF{\mathcal{F}}

\def\calN{\mathcal{N}}
\def\calO{\mathcal{O}}
\def\calP{\mathcal{P}}

\def\calS{\mathcal{S}}

\def\calX{\mathcal{X}}

\def\bA{\boldsymbol{A}}
\def\bB{\boldsymbol{B}}

\def\bD{\boldsymbol{D}}
\def\bE{\boldsymbol{E}}

\def\bH{\boldsymbol{H}}
\def\bI{\boldsymbol{I}}

\def\bM{\boldsymbol{M}}
\def\bN{\boldsymbol{N}}

\def\bQ{\boldsymbol{Q}}
\def\bR{\boldsymbol{R}}
\def\bS{\boldsymbol{S}}

\def\bU{\boldsymbol{U}}
\def\bV{\boldsymbol{V}}
\def\bW{\boldsymbol{W}}
\def\bX{\boldsymbol{X}}

\def\ba{\boldsymbol{a}}
\def\bb{\boldsymbol{b}}

\def\bd{\boldsymbol{d}}

\def\bi{\boldsymbol{i}}

\def\bu{\boldsymbol{u}}
\def\bv{\boldsymbol{v}}
\def\bw{\boldsymbol{w}}
\def\bx{\boldsymbol{x}}
\def\by{\boldsymbol{y}}
\def\bz{\boldsymbol{z}}

\def\bo{\boldsymbol{0}}

\def\balpha{\boldsymbol{\alpha}}
\def\bbeta{\boldsymbol{\beta}}
\def\bgamma{\boldsymbol{\gamma}}

\def\blambda{\boldsymbol{\lambda}}

\def\bmu{\boldsymbol{\mu}}

\def\bpi{\boldsymbol{\pi}}

\def\st{\text{s.t.}}

\definecolor{mygreen}    {RGB}{0,90,0}
\definecolor{myblue}     {RGB}{0,51,140}
\definecolor{myorange}   {RGB}{238,118,0}
\definecolor{myred}      {RGB}{126,0,0}
\definecolor{mygray}     {RGB}{100,100,105}
\definecolor{mygrayblue} {RGB}{0,128,128}
\definecolor{mygraygreen}{RGB}{128,128,0}
\definecolor{DarkPurple}     {RGB}{142, 36, 170}
\definecolor{LightPurple}    {RGB}{57, 130, 7}

\begin{document}

\MANUSCRIPTNO{}

\ARTICLEAUTHORS{
	\AUTHOR{Alexandre Jacquillat, Michael Lingzhi Li, Martin Ramé, Kai Wang}
	\AFF{Sloan School of Management, Massachusetts Institute of Technology,  Cambridge, MA\\
	Harvard Business School, Harvard University, Cambridge, MA\\
	Operations Research Center, Massachusetts Institute of Technology,  Cambridge, MA\\
	School of Vehicle and Mobility, Tsinghua University, Beijing, China
	}
}	
\RUNAUTHOR{Jacquillat, Li, Ramé and Wang}

\RUNTITLE{Branch-and-price for prescriptive contagion analytics}

\TITLE{Branch-and-price for prescriptive contagion analytics}

\ABSTRACT{Predictive contagion models are ubiquitous in epidemiology, social sciences, engineering, and management. This paper formulates a prescriptive contagion analytics model where a decision-maker allocates shared resources across multiple segments of a population, each governed by continuous-time dynamics. We define four real-world problems under this umbrella: vaccine distribution, vaccination centers deployment, content promotion, and congestion mitigation. These problems feature a large-scale mixed-integer non-convex optimization structure with constraints governed by ordinary differential equations, combining the challenges of discrete optimization, non-linear optimization, and continuous-time system dynamics. This paper develops a branch-and-price methodology for prescriptive contagion analytics based on: (i) a set partitioning reformulation; (ii) a column generation decomposition; (iii) a state-clustering algorithm for discrete-decision continuous-state dynamic programming; and (iv) a tri-partite branching scheme to circumvent non-linearities. Extensive experiments show that the algorithm scales to very large and otherwise-intractable instances, outperforming state-of-the-art benchmarks. Our methodology provides practical benefits in contagion systems; in particular, it can increase the effectiveness of a vaccination campaign by an estimated 12-70\%, resulting in 7,000 to 12,000 extra saved lives over a three-month horizon mirroring the COVID-19 pandemic. We provide an open-source implementation of the methodology in an online repository to enable replication.}


\maketitle

\vspace{-12 pt}

\vspace{-12pt}
\section{Introduction}

Epidemiological models have played a central role throughout the COVID-19 pandemic. In the United States for instance, the Center for Disease Control maintained an ensemble forecast based mainly on compartmental contagion models \citep{ray2020ensemble}. At their core, these models rely on susceptible-infected (SI) dynamics, which express the number of new infections proportionally to infected and susceptible individuals. These dynamics have been embedded into more complex models to capture immunization upon recovery (SIR models), time lags from exposure to infection (SEIR models), immunization from vaccinations, asymptomatic cases, quarantine, hospitalization, mortality, etc. Collectively, these models have been instrumental to guide governmental and societal response to the COVID-19 pandemic \citep[see][among many others]{adam2020special,hsiang2020effect,dehning2020inferring,walker2020impact,li2022forecasting,bennouna2022covid}.

Contagion models have a long history that predates the COVID-19 pandemic. Starting from \cite{kermack1927contribution}, SI models have been used to model infectious epidemics, such as influenza \citep{casagrandi2006sirc} and Ebola \citep{berge2017simple}. In marketing, the seminal model of product adoption from \cite{bass1969new} relies on similar dynamics to capture network externalities between existing ``infected'' users and potential ``susceptible'' adopters. Recent applications of this model span multi-generational products \citep{li2013population}, box office revenues \citep{chance2008pricing}, online content diffusion \citep{susarla2012social}, etc. Altogether, SI models---and broader dynamical systems---are ubiquitous in natural sciences, social sciences, engineering and management, with applications to drug addiction \citep{behrens2000optimal}, urban congestion \citep{saberi2020simple}, online rumors \citep{shah2016finding}, climate change mitigation \citep{sterman2012climate}, firm performance \citep{rahmandad2018making}, employee compensation \citep{rahmandad2020if}, etc.

Motivated by the success of \textit{predictive} contagion models, this paper tackles a \textit{prescriptive} contagion analytics problem to optimize spatial-temporal resource allocation over dynamical systems. Specifically, a centralized decision-maker allocates shared resources across multiple segments of a population, each governed by continuous-time dynamics. This paper focuses on a deterministic setting; we show the robustness of the solution under model misspecification but leave prescriptive contagion analytics problems under uncertainty beyond the scope of this work. To demonstrate the broad applicability of the model, we define four problems under its umbrella. The first one involves distributing a vaccine stockpile to combat an epidemic, inspired by the COVID-19 pandemic. The second one adds a discrete facility location structure to deploy mass vaccination centers \citep{bertsimas2022locate}. The third one optimizes online content promotion to maximize product adoption \citep{lin2021content}. The last one deploys emergency vehicles to mitigate urban congestion, based on a new contagion model of traffic congestion developed in this paper using real-world data from Singapore.

Prescriptive dynamical systems---in particular, prescriptive contagion models---involve large-scale mixed-integer non-convex optimization with constraints governed by ordinary differential equations (ODE). As such, they are challenging to even formulate, let alone to solve to optimality. The spatial-temporal resource allocation component, by itself, involves a mixed-integer optimization structure, but the problem is further complicated by complexities of dynamical systems:
\begin{enumerate}[itemsep=0pt,topsep=0pt]
    \item Continuous time dynamics: dynamical systems are governed by ODEs of the form $\frac{d\bM_i(t)}{dt}=f_i(\bM_i(t),\bx_i(t))$, where $\bM_i(t)$ denotes a multi-dimensional state variable in segment $i=1,\cdots,n$; $\bx_i(t)$ denotes a multi-dimensional control variable; and $f_i(\cdot,\cdot)$ refers to a continuous-time transition function. An easy workaround involves time discretization to approximate the dynamics by $\bM_i(t+\Delta t)-\bM_i(t)=f_i(\bM_i(t),\bx_i(t))\cdot\Delta t$. This approach, however, can lead to extensive computational requirements if the time increment $\Delta t$ is too small, or to large approximation errors if it is too large---or both. This challenge is particularly salient in contagion systems, which are highly non-linear and therefore sensitive to even small perturbations in initial conditions. A well-known incarnation of these non-linearities lies in a disease spreading in a population if the basic reproduction number $R_0$ satisfies $R_0>1$ but not if $R_0<1$.
    \item Non-linear interactions: contagion systems are driven by non-convex bilinear interactions of the form $\frac{dS_i(t)}{dt}=-\alpha S_i(t)I_i(t)$, where $S_i(t)$ and $I_i(t)$ refer to the susceptible and infected populations. More broadly, dynamical systems may involve a non-convex transition function $f_i(\cdot,\bx_i(t))$ for a given control variable $\bx_i(t)$. So, even with time discretization, the optimization problem would exhibit a mixed-integer non-convex structure. Accordingly, prescriptive contagion models have typically been solved via approximations and heuristics (see Section~\ref{sec:literature}), but no exact method has been devised for resource allocation in non-linear dynamical systems.
    \item Other non-linearities: system dynamics are endogenous to resource allocation decisions. This dependency can be linear: in epidemiology, for example, immunizations can be safely assumed to be proportional to vaccination rates. But interventions can also induce a non-linear transition function $f_i(\bM_i(t),\cdot)$ for a state variable $\bM_i(t)$ \citep[see, e.g.,][in drug addiction]{behrens2000optimal}. Similarly, dynamical systems may involve non-convex cost functions. Again, these non-linearities create significant complexities in the resulting optimization problem.
\end{enumerate}

\subsubsection*{Contributions and outline.}

This paper develops a branch-and-price methodology for prescriptive contagion analytics. This approach separates the coupled discrete resource allocation decisions in a master problem formulated via mixed-integer linear optimization, and the segment-specific system dynamics in a pricing problem formulated via continuous-state dynamic programming. To our knowledge, it provides the first exact methodology to solve spatial-temporal resource allocation problems in contagion systems---and, more generally, in non-linear dynamical systems.

Specifically, we define a spatial-temporal resource allocation model with endogenous continuous-time non-linear dynamics (Section~\ref{sec:model}). We consider a finite-horizon setting with multiple segments of a population, each constituting a dynamical system governed by non-linear ODEs. The decision-maker optimizes the discrete allocation of a shared resource in each period, which impacts the system dynamics in each segment. This model allows a variety of objective functions and constraints in each segment---even non-linear ones---along with coupling polyhedral constraints across segments (e.g., a shared budget). We apply this formulation to define our four problems: vaccine distribution, deployment of vaccination centers, content promotion, and congestion mitigation. In particular, a byproduct of this research is a new data-driven contagion model of traffic congestion that yields significant improvements in predictive performance against state-of-the-art benchmarks.

Next, this paper develops a branch-and-price algorithm to solve the mixed-integer non-linear optimization problems with ODE constraints (Section~\ref{sec:branchandprice}). The algorithm relies on four components:
\begin{enumerate}
    \item \textit{A set partitioning reformulation.} This formulation uses composite variables to select a resource allocation plan in each segment for the full planning horizon, as opposed to optimizing natural resource allocation decisions for each period. By pre-processing the system dynamics into plan-based variables, this formulation eliminates the three complexities of prescriptive contagion analytics: continuous time dynamics, non-linear interactions, and other non-linearities. However, it comes at the cost of an exponential number of plan-based variables.
    \item \textit{A scalable column generation scheme to solve its linear relaxation.} A master problem solves a coupled resource allocation problem based on a subset of plan-based variables, via mixed-integer linear optimization. A pricing problem adds new plans of negative reduced cost or proves that none exists. Thanks to the problem's structure, the pricing problem can be decomposed into independent segment-specific dynamic programming models. Yet, the system dynamics lead to a continuous state space---a notorious challenge in dynamic programming.
    \item \textit{A state-clustering algorithm for discrete-decision continuous-state dynamic programming.} Forward-enumeration backward-induction algorithms can solve the pricing problem but remain intractable in even small instances. Thus, we propose a linear-time clustering algorithm that exploits the concentration of states in dynamical systems without using the value function or the policy function (which varies from iteration to iteration in column generation). The algorithm provides guarantees on the $\ell_\infty$-diameter of each cluster, which, we prove, controls the global approximation error. The reduction in the size of the state space considerably enhances the scalability of the pricing problem at limited costs in terms of approximation errors.
    \item \textit{A novel tri-partite branching scheme to circumvent the non-linearities of the system.} We embed the column generation procedure into a branch-and-price structure to restore the integrality of resource allocation decisions. We adopt the typical approach of branching on natural resource allocation variables, as opposed to composite plan-based variables. Because of the non-linear dynamics, however, integral resource allocation decisions may not map into equivalent integral plan-based variables. We therefore develop a tri-partite branching scheme that retains a natural branching structure, while guaranteeing finite convergence and optimality.
\end{enumerate}

Finally, this paper demonstrates the scalability of our methodology to otherwise-intractable prescriptive contagion analytics problems (Section~\ref{sec:results}). For resource allocation problems (vaccine allocation, content promotion, and congestion mitigation), the algorithm returns provably optimal or near-optimal solutions in manageable computational times, significantly outperforming state-of-the-art benchmarks. For instance, our algorithm solves vaccine allocation instances involving 21 decisions in each of 51 regions and 12 weeks, with $\calO(21^{612})$ possible decisions. The vaccinations centers problem features an even more challenging facility location structure with linking constraints; still, the branch-and-price methodology generates high-quality solutions and strong bounds in manageable computational times. From a technical standpoint, the methodology significantly outperforms state-of-the-art optimization benchmarks,. Specifically, off-the-shelf implementation based on mixed-integer quadratic optimization or discretization-based linearization do not scale to even small instances; moreover, our methodology provides significant benefits as compared to a tailored coordinate descent heuristic used to circumvent the bilinearities of the problem. From a practical standpoint, the methodology can have a significant impact on the management of contagion-based systems. For instance, the optimized solution can increase the effectiveness of a vaccination campaign by 12-70\% compared to epidemiological benchmarks, resulting in 7,000 to 12,000 extra saved lives over a three-month horizon mirroring the COVID-19 pandemic. These benefits are found to be consistent across contagion analytics problems and robust to model misspecification. Ultimately, our prescriptive contagion analytics approach can deliver significant practical benefits in a variety of domains, by fine-tuning resource allocations based on spatial-temporal system dynamics.
\section{Literature review}\label{sec:literature}

\subsubsection*{Prescriptive contagion analytics.}

The prevalence of contagion models has motivated prescriptive methods to optimize interventions in contagion-based systems. In epidemiology for instance, \cite{goldman2002cost} showed the benefits of social planning interventions against user equilibrium behaviors. \cite{rowthorn2009optimal} found, using a two-region SIS model, that medical interventions should target regions with fewer infections, whereas \cite{ndeffo2011resource} prescribed, using a SIRS model, to prioritize regions with more infections. \cite{yamin2013incentives} designed incentivization schemes for influenza vaccinations by embedding SI-based dynamics into an equilibrium model of vaccination behaviors. In the context of COVID-19, several models traded off health impacts and economic costs to design differentiated lockdowns across age tranches \citep{acemoglu2021optimal}, and to optimize the timing, duration and intensity of lockdowns \citep{caulkins2020long,alvarez2021simple,balderrama2022optimal}. Similarly, the Bass model has been used to guide operational and marketing decisions surrounding product innovation, such as sales planning and time-to-market \citep{ho2002managing}; pricing, production and inventory \citep{shen2014optimal}; multi-product pricing \citep{li2020optimal}; free access and premium subscriptions \citep{mai2022optimizing}; dynamic pricing \citep{cosguner2022dynamic,zhang2022data,agrawal2021dynamic}; etc. The Bass model has also been used to design drug prevention and treatment policies \citep{behrens2000optimal}, advertising campaigns \citep{krishnan2006optimal}, crowdfunding campaigns \citep{zhang2022revenue}, etc. Methodologically, this literature relies on dynamic programming and control to characterize optimal interventions in a single dynamical system (or two systems).

In contrast, our paper considers spatial-temporal resource allocation decisions with coupling constraints across multiple contagion systems. A seemingly related problem is the ventilator sharing problem from \cite{mehrotra2020model} and \cite{bertsimas2021predictions}. However, ventilator availability does not impact the dynamics of the pandemic, so this problem could be formulated via mixed-integer linear optimization by decoupling upstream epidemiological predictions from downstream ventilator allocation decisions. In sharp contrast, vaccinations impact contagion dynamics, which require to integrate epidemiological dynamics into prescriptive resource allocation models.

Thus, our problem involves discrete optimization with endogenous contagion dynamics, combining the mixed-integer optimization difficulties of resource allocation and the non-convex dynamics of contagion models. Existing methods for this class of problems rely on heuristics and approximations. In product adoption, \cite{alban2022resource} used a heuristic to optimize the deployment of mobile healthcare units; and \cite{lin2021content} designed an approximation algorithm with a $1-1/e$ guarantee. In epidemiology, \cite{long2018spatial} used myopic linear optimization and approximate dynamic programming to allocate Ebola treatment units. \cite{bertsimas2022locate} optimized where to open COVID-19 mass vaccination facilities and how to allocate vaccines in the United States, using a coordinate descent heuristic. \cite{fu2021robust} solved a robust vaccine allocation problem with uncertainty in epidemiological forecasts, using linear approximations based on discretized numbers of infections or McCormick reformulations \citep{fu2021robust}. Our paper tackles a deterministic prescriptive contagion analytics problem, and contributes an exact optimization approach that does not rely on finite difference approximations of continuous-time ODE dynamics and that does not involve approximations and heuristics to handle the non-convexities of the problem.

\subsubsection*{Mixed-integer non-linear optimization.}

Our problem exhibits a mixed-integer non-linear optimization (MINLO) structure with ODE constraints,. It therefore remains highly challenging even with finite difference approximations of ODEs \citep[see][for a detailed review of MINLO]{burer2012non}. Bilinear SI-based models fall into mixed-integer quadratic optimization. Methods include the reformulation-linearization technique \citep{sherali2013reformulation}, semi-definite relaxations \citep{fujie1997semidefinite,anstreicher2009semidefinite}, disjunctive programming \citep{saxena2010convex,saxena2011convex}, perspective reformulations \citep{gunluk2010perspective,anstreicher2021quadratic}, augmented Lagrangian duality \citep{gu2020exact}, etc. General-purpose MINLO methods include spatial branch-and-bound \citep{lee2001global}, branch-and-reduce \citep{ryoo1996branch,tawarmalani2004global}, and $\alpha$-branch-and-bound \citep{androulakis1995alphabb}. Recent versions of \cite{gurobi} tackle mixed-integer quadratic optimization problems by combining RLT and spatial branch-and-bound, which we use as a benchmark in Section~\ref{sec:results}.

The ODE constraints also link to the mixed-integer optimal control literature. When the feasible region can be enumerated, partial outer approximation and sum-up rounding can converge to the optimal solution as time discretization becomes infinitely granular \citep{sager2012integer,hante2013relaxation,manns2020improved}. In practice, however, very granular time discretization may be required, and optimality is not guaranteed in the presence of control costs or combinatorial constraints. \cite{jung2015lagrangian} proposed a continuous optimal control problem and a combinatorial integral approximation problem to restore integer feasibility; \cite{gottlich2021penalty} separated a mixed-integer optimal control problem from time-coupling combinatorial constraints; and \cite{bestehorn2021mixed} restored integral feasibility via a shortest path formulation. Our paper departs from this literature by considering large-scale problems with spatial coupling across multiple dynamical systems, and by proposing a new branch-and-price methodology.
 
Since its introduction by \citet{barnhart1998branch}, branch-and-price has been widely applied to mixed-integer linear optimization. In non-linear optimization, \cite{andreas2008branch} solved a reliable $h$-paths problem with non-linearities stemming from failure probabilities. \cite{nowak2018decomposition} proposed an inner-and-outer-approximation of MINLO problems by combining column generation with non-linear optimization. \cite{allman2021branch} synthesized a generic branch-and-price algorithm for non-convex mixed-integer optimization. Our paper proposes a tailored branch-and-price decomposition for spatial-temporal resource allocation problems over continuous-time dynamical systems---including, notably, a novel tri-partite branching disjunction to handle non-linearities.

Finally, one of the main bottlenecks of our methodology lies in the discrete-decision continuous-state dynamic programming structure of the pricing problem \citep[see, e.g.,][for reviews of approximate dynamic programming and reinforcement learning]{bertsekas2015dynamic,powell2022reinforcement}. Continuous-state problems are typically handled via state discretization. \cite{bertsekas1975convergence} showed that static state discretization converges to the optimum as the grid becomes increasingly granular. To manage the growth in the number of states, adaptive discretization methods build a finer grid where the cost function changes rapidly \citep{grune2004using,borraz2011minimizing}. Several studies combined adaptive discretization with reinforcement learning to aggregate states with similar cost-to-go functions \citep[see, e.g.][]{bertsekas1988adaptive,pyeatt2001decision,lolos2017adaptive,sinclair2022adaptive}. \cite{bennouna2021learning} used transition data to learn a discrete partition of a continuous state space. In our problem, however, these methods are not readily applicable because the cost function changes at each column generation iteration. Instead, our state-clustering approach exploits stability in the state space---as opposed to stability in the cost-to-go function.
\section{Prescriptive contagion analytics: definition and formulation}\label{sec:model}

\subsection{Model formulation}
\label{subsection:Formulation}

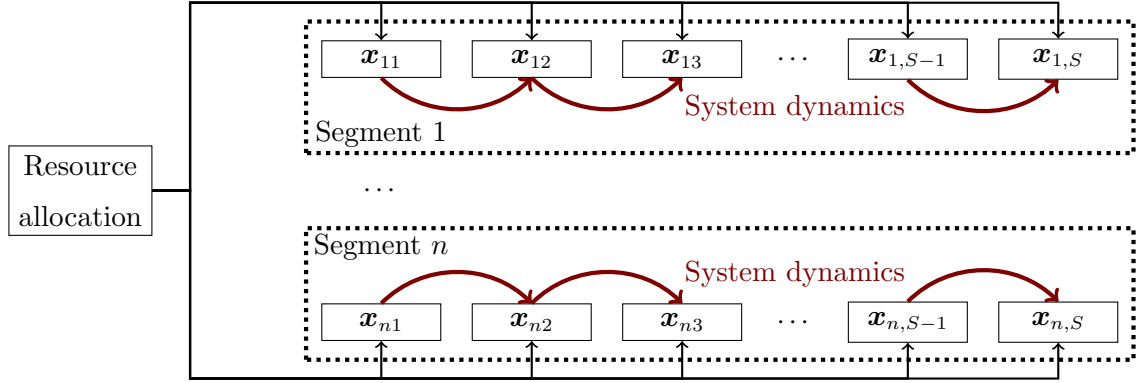
\begin{figure}[h!]
\centering
    \begin{tikzpicture}[every text node part/.style={align=center}]
        \node[draw=black,minimum width=45pt] (RA) at (-2,0){Resource\\allocation};
        \node[draw=black,minimum width=45pt] (x11) at (2,1.75) {$\bx_{11}$};
        \node[draw=black,minimum width=45pt] (x12) at (4,1.75) {$\bx_{12}$};
        \node[draw=black,minimum width=45pt] (x13) at (6,1.75) {$\bx_{13}$};
        \node[] () at (7.5,1.75) {$\cdots$};
        \node[draw=black,minimum width=45pt] (x1S1) at (9,1.75) {$\bx_{1,S-1}$};
        \node[draw=black,minimum width=45pt] (x1S) at (11,1.75) {$\bx_{1,S}$};
        \draw[draw=black,ultra thick,dotted] (1,.5) rectangle ++(11,1.75);
        \node[] () at (2,.75) {Segment $1$};
        \draw[->, color=myred, ultra thick] (x11.south) to [bend right=45] (x12.south);
        \draw[->, color=myred, ultra thick] (x12.south) to [bend right=45] (x13.south);
        \draw[->, color=myred, ultra thick] (x1S1.south) to [bend right=45] (x1S.south);
        \node[color=myred] () at (7.5,1.1) {System dynamics};
        \draw[->, color=black,thick] (RA.east) -- +(0.5,0) -- +(0.5,2.5) -| (x11.north);
        \draw[->, color=black,thick] (RA.east) -- +(0.5,0) -- +(0.5,2.5) -| (x12.north);
        \draw[->, color=black,thick] (RA.east) -- +(0.5,0) -- +(0.5,2.5) -| (x13.north);
        \draw[->, color=black,thick] (RA.east) -- +(0.5,0) -- +(0.5,2.5) -| (x1S1.north);
        \draw[->, color=black,thick] (RA.east) -- +(0.5,0) -- +(0.5,2.5) -| (x1S.north);
        \node[] () at (2,0) {\textbf{$\cdots$}};
        \node[draw=black,minimum width=45pt] (xn1) at (2,-1.75) {$\bx_{n1}$};
        \node[draw=black,minimum width=45pt] (xn2) at (4,-1.75) {$\bx_{n2}$};
        \node[draw=black,minimum width=45pt] (xn3) at (6,-1.75) {$\bx_{n3}$};
        \node[] () at (7.5,-1.75) {$\cdots$};
        \node[draw=black,minimum width=45pt] (xnS1) at (9,-1.75) {$\bx_{n,S-1}$};
        \node[draw=black,minimum width=45pt] (xnS) at (11,-1.75) {$\bx_{n,S}$};
        \draw[draw=black,ultra thick,dotted] (1,-2.25) rectangle ++(11,1.75);
        \node[] () at (2,-.75) {Segment $n$};
        \draw[->, color=myred, ultra thick] (xn1.north) to [bend left=45] (xn2.north);
        \draw[->, color=myred, ultra thick] (xn2.north) to [bend left=45] (xn3.north);
        \draw[->, color=myred, ultra thick] (xnS1.north) to [bend left=45] (xnS.north);
        \node[color=myred] () at (7.5,-1.1) {System dynamics};
        \draw[->, color=black,thick] (RA.east) -- +(0.5,0) -- +(0.5,-2.5) -| (xn1.south);
        \draw[->, color=black,thick] (RA.east) -- +(0.5,0) -- +(0.5,-2.5) -| (xn2.south);
        \draw[->, color=black,thick] (RA.east) -- +(0.5,0) -- +(0.5,-2.5) -| (xn3.south);
        \draw[->, color=black,thick] (RA.east) -- +(0.5,0) -- +(0.5,-2.5) -| (xnS1.south);
        \draw[->, color=black,thick] (RA.east) -- +(0.5,0) -- +(0.5,-2.5) -| (xnS.south);
    \end{tikzpicture}
\caption{Schematic representation of the spatial-temporal resource allocation problem in dynamical systems.}
\label{fig:problem}
\vspace{-12pt}
\end{figure}

We consider a general problem of spatial-temporal resource allocation in dynamical systems, depicted in Figure~\ref{fig:problem}. The ``spatial'' component refers to resource allocation across $n$ segments of a population (e.g., across regions in vaccine allocation, across products in content promotion). The temporal component refers to $S$ decision epochs throughout a planning horizon of length $T$. We denote by $\tau_s$ the time stamp of epoch $s=1,\cdots,S$, with $\tau_1=0$ and $\tau_{S+1}=T$.

We characterize centralized resource allocation decisions via variables $\bx_{is}\in\calF_{is}\subseteq\R^{d_{is}}$ for each segment $i=1,\cdots,n$ and each epoch $s=1,\cdots,S$. These variables are defined as $d_{is}$-dimensional vectors to allow for the allocation of multiple resources (e.g., treatment vs. prevention vehicles in congestion mitigation). We also introduce a variable $\by\in\Z^q\times\R^r$ to capture other decisions (e.g., facility location in our vaccination centers problem). Let $\Delta_s=\sum_{i=1}^nd_{is}$; let $\bX_s\in\R^{\Delta_s}$ denote the concatenated resource allocation variable at epoch $s=1,\cdots,S$. Resource allocation in segment $i=1,\cdots,n$ at epoch $s=1,\cdots,S$ comes at a cost $\Gamma_{is}(\bx_{is})$, and other decisions come at a linear cost $\bd^\top\by$. We make no restriction on the feasible regions $\calF_{is}$ and the cost functions $\Gamma_{is}(\cdot)$ governing segment-specific resource allocations, thus allowing non-linear decision problems. Importantly, we focus on discrete resource allocation decisions, where each region $\calF_{is}$ has finite cardinality $D_{is}=|\calF_{is}|<\infty$. In many cases, resources are naturally restricted to discrete quantities (e.g., emergency vehicles, in our congestion mitigation problem). Even when the decision space is infinite, operational constraints often lead to a discretized set of possible decisions. For example, \cite{moderna} shipped COVID-19 vaccines in pallets of around 20,000 vaccines. Still, the resource allocation problem remains high-dimensional due to coupled spatial-temporal resource allocations---if $D_{is}=\calO(D)$ for all $i$, $s$, the decision space grows in $\calO(D^{n\times S})$, which quickly becomes very large.

Besides segment-specific constraints captured in $\calF_{is}$, global constraints are characterized by a polyhedral set $\left\{(\bX,\by)\in\R^{\Delta_s+q+r}:\bU_s\bX+\bV_s\by\geq\bw_s\right\}$, where $\bU_{s}\in\R^{m_s\times\Delta_s}$, $\bV_s\in\R^{m_s\times(q+r)}$ and $\bw_s\in\R^{m_s}$. In other words, resource allocation decisions are subject to the following linear constraints, where $\bu_{sji}\in\R^{d_{is}}$ denotes the vector corresponding to the $j^{\text{th}}$ constraint and the $i^{\text{th}}$ segment in $\bU_s$, and $\bv_{sj}\in\R^{q+r}$ denotes the vector corresponding to the $j^{\text{th}}$ constraint in $\bV_s$:
\vspace{-6pt}
$$\left(\sum_{i=1}^n\bu_{sji}^\top\bx_{is}\right)+\bv_{sj}^\top\by\geq w_{sj},\ \forall s=1,\cdots,S,\ \forall j=1,\cdots,m_s$$

Next, each population segment $i=1,\cdots,n$ constitutes a dynamical system, governed by a continuous-time state variable $\bM_i(t)\in\R^{r_i}$. We assume that the dynamics are independent across the $n$ segments. In the epidemiological context, for instance, this captures intra-region interactions (e.g., within each state or each country) but not inter-region interactions \citep[see, e.g.,][for similar assumptions]{hsiang2020effect,walker2020impact,li2022forecasting,bennouna2022covid}. Specifically, in each segment $i=1,\cdots,n$, the state variable $\bM_i(t)$ is determined by an initial condition $\bM^{0}_i$, and varies according to the following system of ODEs, where $f_i(\cdot,\cdot)$ denotes the transition function:
$$\frac{d\bM_i(t)}{dt}=f_i(\bM_i(t),\bx_{is}),\ \forall s=1,\cdots,S,\ \forall t\in\left[\tau_s,\tau_{s+1}\right].$$
The system is associated with a continuous-time cost function $g_{it}(\bM_i(t))$ and a terminal cost function $h_i(\bM_i(T))$. We impose no restriction on the transition function $f_i(\cdot,\cdot)$ and the cost functions $g_{it}(\cdot)$ and $h_i(\cdot)$ so, by design, our model encompasses non-convex dynamical systems.

The spatial-temporal resource allocation problem, referred to as Problem $(\calP)$, minimizes total costs subject to resource allocation constraints and the system's dynamics. It is written as follows:
\begin{align}
(\calP)\quad\min \quad & \sum_{i=1}^n\left(\int_0^Tg_{it}(\bM_i(t))dt+h_i(\bM_i(T))\right)+\sum_{i=1}^n\sum_{s=1}^S\Gamma_{is}(\bx_{is})+\bd^\top\by \label{eq:cost}\\
\st \quad & \left(\sum_{i=1}^n\bu_{sji}^\top\bx_{is}\right)+\bv_{sj}^\top\by\geq w_{sj},\ \forall s=1,\cdots,S,\ \forall j=1,\cdots,m_s \label{eq:constraintsC} \\
& \frac{d\bM_i(t)}{dt}=f_i(\bM_i(t),\bx_{is}),\ \forall i=1,\cdots,n,\ \forall s=1,\cdots,S,\ \forall t\in\left[\tau_s,\tau_{s+1}\right] \label{eq:ODE}\\
& \bM_i(\tau_1)=\bM^{0}_i,\ \forall i=1,\cdots,n \label{eq:initial}\\
& \bx_{is}\in\calF_{is},\ \forall i=1,\cdots,n,\ \forall s=1,\cdots,S \label{eq:constraints_i} \\
& \by\in\Z^q\times\R^r\label{eq:domain_y}
\end{align}

\vspace{-6pt}

This formulation, however, is intractable due to continuous time dynamics (Equation~\eqref{eq:ODE}), non-convex system dynamics (functions $f_i(\cdot,\cdot)$, $g_{it}(\cdot)$ and $h_i(\cdot)$), and non-convex segment-wise decisions (functions $\Gamma_{is}(\cdot)$, regions $\calF_{is}$). One workaround would be to model $(\calP)$ via dynamic programming, by defining a state variable that encompasses \textit{all} variables $\bM_1(\tau_s),\cdots,\bM_n(\tau_s)$ as well as all required information to enforce resource allocation constraints (Equation~\eqref{eq:constraintsC}). This approach, however, would scale in $\calO(D^{n\times S})$ if $D_{is}=\calO(D)$ for all $i$, $s$, quickly growing into the curse of dimensionality. Instead, we propose in Section~\ref{sec:branchandprice} a branch-and-price methodology that separates coupled resource allocation decisions across segments and continuous-time dynamics in each segment. We still use dynamic programming in the pricing problem, but we leverage segment-specific state variables $\bM_i(t)\in\R^{r_i}$ rather than a higher-dimensional state variable of the form $(\bM_1(\tau_s),\cdots,\bM_n(\tau_s))\in\R^{r_1}\times\cdots\times\R^{r_n}$. Thus, the pricing problem scales in $\calO(D^{S})$, instead of $\calO(D^{n\times S})$.

Finally, let us underscore that, following \cite{powell2022reinforcement} we define the state variable $\bM_i(t)$ as a necessary and sufficient function of history to compute the cost function, the constraints, and the transition function in segment $i=1,\cdots,n$. For example, an inter-temporal budget constraint of the form $\sum_{s=1}^S\bx_{is}\leq \bgamma_i$ would require state augmentation to capture the budget used up to epoch $s$. In contrast, our model can accommodate budget constraints across segments of the form $\sum_{i=1}^n\bx_{is}\leq\bbeta_s$. Although the problem is not directly modeled as a dynamic program, this structure enables segment-wise dynamic programming decomposition in our algorithm. 

\subsection{Prescriptive contagion analytics model, and applications}

The special case of contagion systems is derived from the general formulation by defining, in each segment $i=1,\cdots,n$, a susceptible state $S_i(t)\in\R$, an infected state $I_i(t)\in\R$, and other states represented via a vector $\bR_i(t)\in\R^V$. We capture bilinear interactions between susceptible and infected populations as follows, with infection rate $\alpha$ and transitions $f^S_i(\cdot)$, $f^I_i(\cdot)$, $f^R_i(\cdot)$:
\begin{align}
&\frac{dS_i(t)}{dt} = -\alpha S_i(t)I_i(t) + f^S_i(S_i(t),I_i(t),\bR_i(t))\label{eq:SIS}\\
&\frac{dI_i(t)}{dt} = +\alpha S_i(t)I_i(t) + f^I_i(S_i(t),I_i(t),\bR_i(t))\label{eq:SII}\\
&\frac{dR_{iv}(t)}{dt} = f^R_i(S_i(t),I_i(t),\bR_i(t)),\quad\forall v=1,\cdots,V \label{eq:SIR}
\end{align}
Similarly, we define functions $g^{SIR}_{it}(\cdot)$ and $h^{SIR}_{i}(\cdot)$ to characterize the costs of the system dynamics in each segment. We then formulate the prescriptive contagion model by minimizing the following objective function, subject to the resource allocation constraints (Equation~\eqref{eq:constraintsC} and Equations~\eqref{eq:constraints_i}--\eqref{eq:domain_y}) and the contagion dynamics (initial conditions, and Equations~\eqref{eq:SIS}--\eqref{eq:SIR}):
\begin{align}
\sum_{i=1}^n\left(\int_0^Tg^{SIR}_{it}(S_i(t),I_i(t),\bR_i(t))dt+h^{SIR}_{i}(S_i(T),I_i(T),\bR_i(T))\right)+\sum_{i=1}^n\sum_{s=1}^S\Gamma_{is}(\bx_{is})+\bd^\top\by\label{eq:costSI}
\end{align}

\begin{figure}[h!]
\centering
\subfloat[Vaccine allocation and vaccination centers \citep{bertsimas2022locate}.]{\label{subfig:vaccines}
    \begin{tikzpicture}[scale=0.5, every node/.style={scale=1}]
        \node[draw=black] (S1) at (0, 4){$S$};
        \node[draw=black] (E1) at (6, 4){$E$};
        \node[draw=black] (I1) at (12, 4) {$I$};
        \node[draw=black] (UD1) at (18, 6) {$U$};
        \node[draw=black] (HD1) at (18, 4) {$H$};
        \node[draw=black] (QD1) at (18, 2) {$Q$};
        \node[draw=black] (D1) at (24,4) {$D$};
        \draw[->,dotted,ultra thick] (S1)--(E1);
        \draw[->] (E1)--(I1);
        \draw[->, draw=black] (I1)--(HD1);
        \draw[->, draw=black] (I1)--(QD1);
        \draw[->, draw=black] (I1)--(UD1);
        \draw[->, draw=black] (UD1)--(D1);
        \draw[->, draw=black] (HD1)--(D1);
        \draw[->, draw=black] (QD1)--(D1);
        \node[draw=black] (SV) at (0, 0){$S’$};
        \node[draw=black] (EV) at (6, 0){$E’$};
        \node[draw=black] (IV) at (12, 0) {$I’$};
        \node[draw=black] (Immune) at (24,0) {$M$};
        \draw[->,dotted,ultra thick] (SV)--(EV);
        \draw[->] (EV)--(IV);
        \draw[->] (IV)--(Immune);
        \draw[->, color=myred, ultra thick] (S1.south) to [bend right=20] node[right,midway]{$\beta x$} (SV.north);
        \draw[draw=black] (13,-1.25) rectangle (11,5);
        \draw[->, draw=black] (11,2.5) to [bend left=10] (S1.south east);
        \draw[->, draw=black] (11,1.5) to [bend right=10] (SV.north east);
    \end{tikzpicture}}\\
\subfloat[Content promotion \citep{lin2021content}.]{\label{subfig:content}
    \begin{tikzpicture}[scale=0.5, every node/.style={scale=1}]
        \node[draw=black] (S) at (0,8){$A$};
        \node[draw=white] (dummy) at (6,-2){};
        \node[draw=white] (dummy) at (-6,-2){};
        \node[draw=black] (B) at (0,0){$B$};
        \draw[->, color=myred,dotted,ultra thick] (S) -- node[right]{$\alpha Ax^2+\beta BA$} (B);
    \end{tikzpicture}}\hspace{0.1cm}
\subfloat[Urban congestion (this paper, Appendix~\ref{app:traffic}).]{\label{subfig:traffic}
    \begin{tikzpicture}[scale=0.5, every node/.style={scale=1}]
        \node[draw=black] (S) at (0, 6){$S$};
        \node[draw=black] (W) at (0, 2){$W$};
        \node[draw=black] (A) at (10, 8){$A$};
        \node[draw=black] (I) at (10, 4){$I$};
        \node[draw=black] (Aprime) at (10, 0){$A'$};
        \node[draw=black] (R) at (16, 4){$R$};
        \draw[->,color=black] (S) to node[above]{$\beta^FS$}(A);
        \draw[<-,color=myred,ultra thick] (S)to[bend left=30] node[above,xshift=-24pt]{$\rho^F(\zeta^F+\psi(x^1))A$}(A);
        \draw[->,dotted,ultra thick] (S)-- node[above,xshift=18pt]{$\alpha^FS(S+A+A')$}(I);
        \draw[->,dotted,ultra thick] (W)-- node[below,xshift=18pt]{$\alpha^WW(S+A+A')$}(I);
        \draw[->,color=black] (W) to node[below]{$\beta^WW$}(Aprime);
        \draw[<-,color=myred,ultra thick] (W) to[bend right=30] node[below,xshift=-24pt]{$\rho^W(\zeta^W+\psi(x^1))A'$}(Aprime);
        \draw[->,color=myred,ultra thick] (A)-- node[right]{$(1-\rho^S)(\zeta^F+\varphi(x^2))$}(I);
        \draw[->,color=myred,ultra thick] (Aprime)-- node[right]{$(1-\rho^W)(\zeta^W+\varphi(x^2))$}(I);
        \draw[->] (I)-- node[above]{$\theta I$}(R);
    \end{tikzpicture}}
\caption{Contagion models for the four problems. Thick red lines indicate transitions that are endogenous to resource allocation. Dotted lines indicate bilinear susceptible-infected interactions.}
\label{fig:SIR}
\vspace{-12pt}
\end{figure}
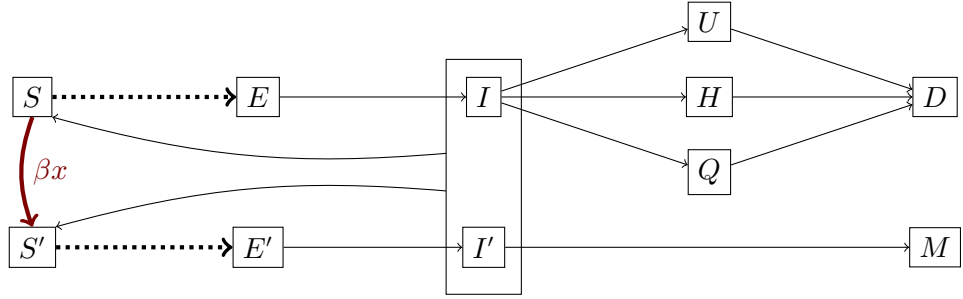
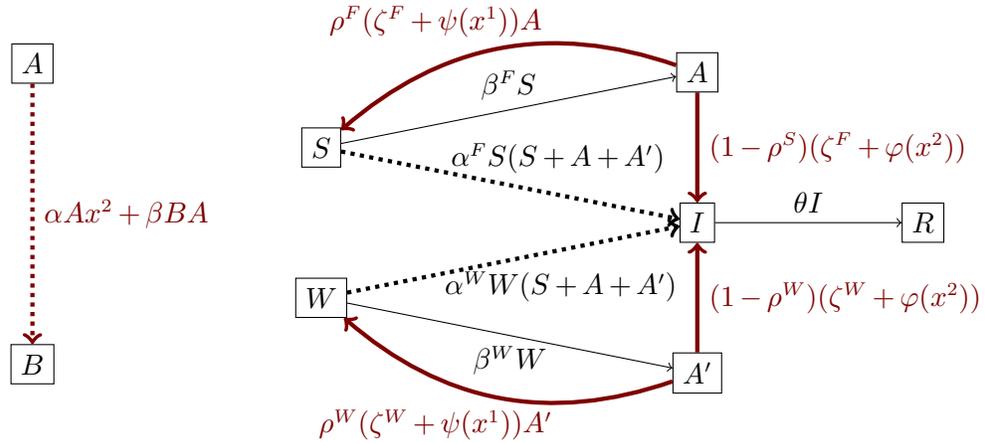

We define our four problems from this generic formulation, all inspired from real-world problems and data-driven contagion models. We defer details and their full formulation to Appendix~\ref{app:formulation}.

\paragraph{Vaccine allocation:}
The problem optimizes the number of vaccines ($x$) received in 51 regions (50 US states plus Washington, DC) to maximize the number of saved lives. We use the DELPHI-V model of COVID-19 from \cite{li2022forecasting} and \cite{bertsimas2022locate}, with 7 non-vaccinated states (susceptible (S), exposed (E), infected (I), undetected (U), hospitalizations (H), quarantine (Q), death (D)) and 4 vaccinated states (susceptible (S’), exposed (E’), infected (I’), immunized (M)). Our experiments consider allocating a weekly stockpile of 2.5 to 7 million vaccines across 51 states and 12 weeks. We discretize decisions into up to 21 decisions per epoch, with $\calO\left(21^{612}\right)$ possible decisions. This setup mirrors the situation from February to April 2021 in the United States.

\paragraph{Vaccination centers.}
The problem extends the vaccine allocation problem by jointly optimizing the location of mass vaccination centers and the subsequent distribution of vaccines within a geographical radius \citep{bertsimas2022locate}. Due to the complexity of the problem, we break down the United States into 10 groups defined by the Centers for Disease Control, with 4-7 states and 20 candidate facilities in each (Appendix~\ref{app:setup}). Each group receives a budget of vaccines proportionally to its share of the US population, discretized into 25,000-vaccine pallets, plus some flexibility buffer. We set a budget of 20 and 30 facilities across the country, again mirroring the situation in 2021. With 20 (resp. 30) facilities, each group needs to select 2 (resp. 3) facilities, and residents can travel up to 150 miles (resp. 100 miles), possibly across states, to access a vaccination center.

\paragraph{Content promotion.}

This problem optimizes product promotion to maximize adoption, subject to a sparsity constraint to promote up to $K$ products per period, a cover constraint, a linking constraint, and contagion dynamics \citep{lin2021content}. The decision variables includes which products should be promoted ($x^1$ variables) and the share of the population to which each product shall be promoted ($x^2$ variables). Each segment corresponds to a product, and adoption follows a Bass model with susceptible users (A) and adopters (B). In our experiments, we consider up to 20 products, 10 decision epochs, 21 decisions per epoch, and 2-6 products promoted per epoch.

\paragraph{Congestion mitigation.}

This problem deploys $B^1$ prevention vehicles and $B^2$ treatment vehicles per hour to mitigate congestion costs. Prevention resources ($x^1$) respond to minor accidents and prevent congestion, whereas treatment resources ($x^2$) respond to major accidents and ease recovery. Contagion models reflect the propagation of congestion from ``infected'' roads to free-flow ``susceptible'' roads  \citep{saberi2020simple}. We leverage new data sources from Singapore on traffic speeds, road work, and traffic accidents to propose a richer six-state model, (Figure~\ref{subfig:traffic}): susceptible roads (S), road work (W), accidents (A), road work and accidents ($A'$), congested (I), and recovered (R). We calibrate this model using a three-peak structure for morning, afternoon, and evening traffic. A byproduct of this paper is a new contagion model of urban congestion, which improves predictive performance against state-of-the-art benchmarks (Appendix~\ref{app:traffic}). In our experiments, we consider a setting with 5 neighborhoods, up to 8 hours and up to 12 vehicles to allocate, developed in collaboration with city officials to capture real-world situations in Singapore.

\paragraph{Discussion.} The vaccine allocation problem has the largest state space. The vaccination centers problem adds a discrete structure with a switching variable $\by$ (facility location) and linking constraints governing $\bx_{is}$ (vaccine allocation). The other two problems include multi-dimensional variables $\bx_{is}$, and the congestion mitigation problem involves a novel contagion model. Together, they cover a range of applications and optimization structures in prescriptive contagion analytics.
\section{The method: branch-and-price for prescriptive contagion analytics}
\label{sec:branchandprice}

We propose a branch-and-price algorithm to solve Problem $(\calP)$, using a set partitioning reformulation (Section~\ref{subsec:SP}), a column generation procedure (Section~\ref{subsec:CG}), a state-clustering algorithm for the continuous-state pricing problem (Section~\ref{subsec:dynamic_program}), and a tri-partite branching scheme to handle non-convexities (Section~\ref{subsec:BP}). Section~\ref{subsec:summary} summarizes the algorithm and establishes its exactness.

\subsection{Set partitioning reformulation}
\label{subsec:SP}

The reformulation optimizes over composite plan-based variables in each segment (which characterize resource allocations over the full planning horizon), as opposed to natural variables (which characterize resource allocation in each period). Specifically, we define the set of feasible plans in each segment $i=1,\cdots,n$ as the combination of feasible resource allocation decisions, as follows:
\vspace{-6pt}
$$\calP_i=\calF_{i1}\times\cdots\times\calF_{iS},\ \forall i=1,\cdots,n$$
In each segment $i=1,\cdots, n$, a plan $p\in\calP_i$ defines a sequence of resource allocations denoted by $(\balpha^p_{i1},\cdots,\balpha^p_{iS})\in\calF_{i1}\times\cdots\times\calF_{iS}$, as well as a cost parameter $C^p_i$ governed by the system's dynamics:
\begin{align}
C^p_i= & \quad\int_0^Tg_{it}(\bM_i(t))dt+h_i(\bM_i(T))+\sum_{s=1}^S\Gamma_{is}(\bx_{is}) \label{SP:cost}\\
  \st & \quad\frac{d\bM_i(t)}{dt}=f_i(\bM_i(t),\balpha^p_{is}),\ \forall s=1,\cdots,S,\ \forall t\in\left[\tau_s,\tau_{s+1}\right] \label{SP:ODE}\\
 & \quad \bM_i(\tau_1)=\bM^0_i \label{SP:initial}
\end{align}

We define the following plan-based decision variables:
\begin{align*}
z^p_i = \begin{cases}
        1 & \text{if plan $p\in\calP_i$ is selected for segment $i=1,\cdots,n$,}\\
        0 & \text{otherwise.}
    \end{cases}
\end{align*}

The set partitioning formulation, referred to as $(\calS\calP)$, minimizes total costs (Equation~\eqref{SP:obj}), while enforcing the coupling constraints (Equation~\eqref{SP:constraintsC}) and ensuring that one plan is selected in each segment (Equation~\eqref{SP:one}). For completeness, we reformulate our four problems in Appendix~\ref{app:SP}.
\begin{align}
(\calS\calP)\qquad\min \quad & \sum_{i=1}^n\sum_{p\in\calP_i}C^p_iz^p_i+\bd^\top\by \label{SP:obj}\\
\st \quad & \left(\sum_{i=1}^n\sum_{p\in\calP_i}\bu_{sji}^\top\balpha^p_{is}z^p_i\right)+\bv_{sj}^\top\by\geq w_{sj},\ \forall s=1,\cdots,S,\ \forall j=1,\cdots,m_s \label{SP:constraintsC} \\
& \sum_{p \in \calP_i}z^p_i = 1 \quad \forall i=1,\cdots,n  \label{SP:one}\\
& z^p_i \in \{0,1\} \quad \forall i=1,\cdots,n,\ \forall p \in \calP_i \label{SP:domain}\\
& \by\in\Z^q\times\R^r\label{SP:domain_y}
\end{align}

\begin{proposition}\label{prop:SP}
The set partitioning formulation $(\calS\calP)$ is equivalent to Problem $(\calP)$.
\end{proposition}

This reformulation eliminates structural complexities of Problem $(\calP)$ by pre-processing the continuous-time dynamics and the non-linear functions into plan-based variables and the corresponding input parameters. The $(\calS\calP)$ formulation therefore exhibits a mixed-integer linear optimization structure, but involves an exponential number of composite plan-based variables. In an instance with $\calO(D)$ decisions in each segment and at each epoch, the number of plans scales in $\calO(D^{n\times S})$. It is therefore intractable to even enumerate the full set of plans, let alone to estimate the cost parameters by running a dynamical model for each one (Equations~\eqref{SP:cost}--\eqref{SP:initial}).

\subsection{Solving the linear optimization relaxation via column generation}
\label{subsec:CG}

To address this challenge, we generate plans iteratively via column generation. A restricted master problem (RMP) solves the linear relaxation of $(\calS\calP)$ over subsets of plans $\calP^0_i\subseteq\calP_i$:
\begin{align}
(RMP)\qquad\min \quad & \sum_{i=1}^n\sum_{p\in\calP^0_i}C^p_iz^p_i+\bd^\top\by \label{MP:obj}\\
\st \quad & \left(\sum_{i=1}^n\sum_{p\in\calP^0_i}\bu_{sji}^\top\balpha^p_{is}z^p_i\right)+\bv_{sj}^\top\by\geq w_{sj},\ \forall s=1,\cdots,S,\ \forall j=1,\cdots,m_s \label{MP:constraintsC} \\
& \sum_{p \in \calP^0_i}z^p_i = 1 \quad \forall i=1,\cdots,n  \label{MP:one}\\
& z^p_i \geq0 \quad \forall i=1,\cdots,n,\ \forall p \in \calP^0_i \label{MP:domain}\\
& \by\in\R^{q+r}\label{MP:domain_y}
\end{align}

Let $\lambda_{sj}\in\R_+$ be the dual variable of the coupling constraint (Equation~\eqref{MP:constraintsC}) and $\mu_i\in\R$ the dual variable of the set partitioning constraint (Equation~\eqref{MP:one}). The pricing problem exploits segment-wise decomposition to generate new plan-based variables or a certificate that none exists. Specifically, it seeks the plan in $\calP_i$ with the minimal reduced cost, for each segment $i=1,\cdots,n$:
\begin{align}
    (PP_i)\qquad\min\quad   &  \left(\int_0^Tg_{it}(\bM_i(t))dt+h_i(\bM_i(T))+\sum_{s=1}^S\Gamma_{is}(\bx_{is})-\sum_{s=1}^S\sum_{j=1}^{m_s}\lambda_{sj}\bu_{sji}^\top\bx_{is} - \mu_i\right)\label{PP:obj}\\
    \st\quad                &   \bx_{is}\in\calF_{is},\ \forall s=1,\cdots,S \label{PP:constraints_i} \\
                            &   \frac{d\bM_i(t)}{dt}=f_i(\bM_i(t),\bx_{is}),\  \forall s=1,\cdots,S,\ \forall t\in\left[\tau_s,\tau_{s+1}\right] \label{PP:ODE}\\
                            &   \bM_i(\tau_1)=\bM^{0}_i \label{PP:initial}
\end{align}

The pricing problem exhibits a discrete non-linear optimization structure with continuous-time dynamics. However, we can leverage temporal decomposition to formulate it via dynamic programming. By definition, the continuous-time state variable $\bM_i(t)$ encapsulates all necessary system history to compute the cost function, the constraints in $\calF_{is}$, and the transition function $f_i(\cdot,\cdot)$, so the system follows Markovian dynamics. At each epoch $s=1,\cdots,S$, the dynamic programming state variable is denoted by $\bN^{i}_s$, equal to $\bM_i(\tau_s)$. The decision is $\bx_{is}\in\calF_{is}$, and the transition function is governed by the continuous-time dynamics between $\tau_s$ and $\tau_{s+1}$, with initial condition determined by the state variable. The problem involves a per-period cost characterizing the continuous-time cost $g_{it}(\cdot)$ accrued between $\tau_s$ and $\tau_{s+1}$, the cost of resource allocation $\Gamma_{is}(\cdot)$, and the dual price of the coupling constraints. In addition, the problem involves a terminal cost characterizing the end-state cost $h_i(\cdot)$ and the dual price of the set partitioning constraints. The Bellman equation is given as follows, using $J_s(\cdot)$ to define the cost-to-go function at epoch $s=1,\cdots,S$:
\begin{align}
    J_s\left(\bN^{i}_s\right)&=\min_{\bx_{is}\in\calF_{is}}\left\{\int_{\tau_s}^{\tau_{s+1}}g_{it}(\bM_i(t))dt+\Gamma_{is}(\bx_{is})-\sum_{j=1}^{m_s}\lambda_{sj}\bu_{sji}^\top\bx_{is}+J_{s+1}\left(\bN^{i}_{s+1}\right)\right\} \label{eq:bellman_1}\\
    &\quad\text{where $\bM_i(\tau_s)\gets\bN^{i}_s$; $\frac{d\bM_i(t)}{dt}=f_i(\bM_i(t),\bx_{is}),\ \forall t\in\left[\tau_s,\tau_{s+1}\right]$; and $\bN^{i}_{s+1}\gets \bM_i(\tau_{s+1})$}\nonumber\\
    J_{S+1}\left(\bN^{i}_{S+1}\right)&=h_i\left(\bN^{i}_{S+1}\right)-\mu_i\label{eq:BellmanT}
\end{align}

The column generation procedure solves the linear relaxation of the $\calS\calP$ formulation to optimality. Starting with initial sets $\calP^0_i$ of plan-based variables, the RMP provides a feasible primal solution at each iteration, along with the dual variables $\lambda_{sj}$ and $\mu_i$. We then solve $(PP_i)$ for each segment $i=1,\cdots,n$; if its optimal value is negative, we expand the set $\calP^0_i$ with the new solution and proceed. Otherwise, the incumbent RMP solution is optimal for the $(\calS\calP)$ relaxation.

This column generation scheme separates the two complexities of the problem: the RMP handles the coupled resource allocation decisions via mixed-integer optimization (Equations~\eqref{MP:obj}--\eqref{MP:domain}), and the pricing problem handles the continuous-time system dynamics via dynamic programming (Equations~\eqref{eq:bellman_1}--\eqref{eq:BellmanT}). However, the pricing problem exhibits a continuous-state dynamic programming structure---a notoriously challenging class of problems. This motivates our state-clustering dynamic programming algorithm to solve it efficiently at each column generation iteration.

\subsection{A state-clustering dynamic programming algorithm for the pricing problem}
\label{subsec:dynamic_program}

Since the pricing problem is separable across segments $i=1,\cdots,n$, we ignore the index $i$ in this section; for instance, $\bN_s$ refers to the state, and $\bx_s\in\calF_s$ to the decision with $D_s=|\calF_s|$.

\subsubsection*{Exact algorithm.} Any finite-horizon discrete-decision continuous-state dynamic programming model can be solved via forward enumeration and backward induction, as follows:
\begin{enumerate}
\item \textit{Forward enumeration}: From the initial state $\bM^0$, evaluate all possible decisions $\bx_1,\cdots,\bx_S$, and store all possible states at each epoch $s=1,\cdots,S+1$ in an exhaustive set $\calN^*_s$.
\item \textit{Backward induction}: Starting from $s=S$, derive the optimal policy $\bpi_s^*(\cdot)$ at epoch $s=S,\cdots,1$ and update the cost-to-go function $J_s(\cdot)$, as follows:
\vspace{-6pt}
\begin{align}
    \bpi_s^*\left(\bN_s\right)&=\argmin_{\bx_s\in\calF_s}\left\{\int_{\tau_s}^{\tau_{s+1}}g_{t}(\bM(t))dt+\Gamma_s(\bx_s)-\sum_{j=1}^{m_s}\lambda_{sj}\bu_{sj}^\top\bx_s+J_{s+1}\left(\bN_{s+1}\right)\right\}\label{policy}\\
    J_s\left(\bN_s\right)&=\int_{\tau_s}^{\tau_{s+1}}g_{t}(\bM(t))dt+\Gamma_s\left(\bpi_s^*\left(\bN_s\right)\right)-\sum_{j=1}^{m_s}\lambda_{sj}\bu_{sj}^\top\bpi_s^*\left(\bN_s\right)+J_{s+1}\left(\bN_{s+1}\right) \label{eq:costtogo}\\
    &\quad\text{where $\bM(\tau_s)\gets\bN_s$; $\frac{d\bM(t)}{dt}=f(\bM(t),\bpi_s^*\left(\bN_s\right)),\ \forall t\in\left[\tau_s,\tau_{s+1}\right]$; and $\bN_{s+1}\gets \bM(\tau_{s+1})$}\nonumber
\end{align}
\end{enumerate}

This approach is detailed in Algorithm~\ref{alg:dp_full} in Appendix~\ref{app:DP}. Per the Bellman principle of optimality, the optimal path is optimal at every decision point (Proposition~\ref{prop:exactDP}).
\begin{proposition}\label{prop:exactDP}
Algorithm~\ref{alg:dp_full} returns a feasible and optimal solution to Equations~\eqref{eq:bellman_1}--\eqref{eq:BellmanT}.
\end{proposition}

However, Algorithm~\ref{alg:dp_full} is highly computationally expensive due to exhaustive state enumeration in the forward-enumeration and backward-induction loops. Complexity grows in $\calO\left(\prod_{s=1}^S D_s\right)$, or $\calO(D^S)$ if $D_s=\calO(D),\ \forall s$, which can quickly become very large as the planning horizon and the decision space increase. Even for a simple eight-period problem with two decisions that can each take four values, $\prod_{s=1}^SD_s=16^8\sim1$ billion, which severely hinders exact dynamic programming.

As noted earlier, most approximate dynamic programming and reinforcement learning methods depend on learning approximations of the cost-to-go function $J_s(\bN_s)$ and/or the optimal policy function $\bpi_s^*(\bN_s)$. In our problem, however, both functions change from one column generation iteration to the next, due to the updates in the dual values $\lambda_{sj}$ from the RMP. Thus, standard approaches would need be to re-learned at every iteration, resulting in significant inefficiencies. 

Instead, we propose a state-clustering algorithm that exploits the concentration of states, but not the cost function. This approach proceeds by aggregation to reduce the state space during the forward-enumeration step of dynamic programming. Note that forward enumeration does not depend on the cost function, so the state clustering algorithm needs to be applied only once at the beginning of the column generation algorithm. Then,  we leverage the backward-induction algorithm with the clustered state space to speed up the pricing problem at each column generation iteration.

\subsubsection*{A state-clustering acceleration.}

Many dynamical systems exhibit a natural concentration of states. In vaccine allocation, for example, there are only mild differences between a no-vaccine baseline and an allocation of just a few vaccines; similarly, two allocations that merely swap decisions between epochs $s$ and $s+1$ may yield similar future outcomes. We leverage this observation to approximate the full set of states $\calN_s^*$ by a clustered state space $\calN_{s}$ at each epoch $s=1,\cdots,S$.

Specifically, at epoch $s=1,\cdots,S$, we enumerate all pairs of states in $\calN_s$ and all decisions in $\calF_s$; for each one, we compute the next state at period $s+1$ and the corresponding cost accrued between times $\tau_s$ and $\tau_{s+1}$. This yields $Q=|\calN_s| \times |\calF_s|$ state-decision pairs $(\bN_s^1,\bx_s^1),\cdots,(\bN_s^Q,\bx_s^Q)$ with $Q$ subsequent states $\bN_{s+1}^1,\cdots,\bN_{s+1}^Q$ and corresponding cost parameters $\chi_s^1,\cdots,\chi_s^Q$. These are defined as follows:
\begin{align*}
    & \bM^q(\tau_s)\gets\bN_s^q,\quad\frac{d\bM^q(t)}{dt}=f(\bM^q(t),\bx_s^q),\quad\bN^{q}_{s+1}\gets \bM^q(\tau_{s+1}),\quad\forall q=1,\cdots,Q\\
    & \chi_s^q=\int_{\tau_s}^{\tau_{s+1}}g_{t}(\bM^q(t))dt+\mathbf{1}\{s=S\}(h(\bM^q(T))),\quad\forall q=1,\cdots,Q
\end{align*}

We group states $\bN_{s+1}^1,\cdots,\bN_{s+1}^Q$ into $\Omega_{s+1}$ clusters, and treat the centroids as representative states in $\calN_{s+1}$ at epoch $s+1$. We then apply the backward-induction algorithm in the clustered state space. By design, this procedure can reduce the state space significantly as long as $\Omega_{s+1}\ll |\calN_s| \times |\calF_s|$. Obviously, the resulting dynamic program is an approximation of the full dynamic program, but the optimality loss can be small if the full state space is indeed concentrated.

The next question lies in the design of the clustering algorithm. Algorithms based on global similarity measures have at least quadratic complexity to compute the entire distance matrix. In our setup, $|\calN_s| \times |\calF_s|$ is already very large, so a complexity of $\calO\left(|\calN_s|^2 \times |\calF_s|^2\right)$ would be intractable. Most linear-time clustering algorithms, such as $k$-means, cannot effectively bound the diameter of each cluster, potentially leading to significant error propagation in non-linear dynamical systems. Instead, we propose a linear-time clustering algorithm with guarantees on cluster diameter in an $\ell_\infty$-space. Unlike $k$-means, the algorithm creates clusters dynamically without a pre-specified number of clusters. The benefits of our approach are showed theoretically in Proposition~\ref{prop:approx_error} (namely, the bound on cluster diameter bounds the global approximation error) and numerically in Section~\ref{sec:results} (namely, the smaller global approximation error induces stronger solutions than a $k$-means benchmark).

Specifically, for each cluster $\omega\in\{1,\cdots,\Omega_{s+1}\}$, we store: (i) the number of data points $\eta(\omega)$ and the sum of all states $\bN^\Sigma(\omega)$ (to define the centroid); (ii) the element-wise minimum and maximum states, $\underline{\bN}(\omega)$ and $\overline{\bN}(\omega)$ (to bound the cluster diameter); and (iii) the set of state-decision pairs $\calX(\omega)$ that lead to the cluster and the total cost $c(\omega)$ (to define cost functions in the clustered state space). For a diameter tolerance $\varepsilon$, we cluster the states $\bN_{s+1}^1,\cdots,\bN_{s+1}^Q$ as follows
\begin{enumerate}
\item We assign the first state $\bN_{s+1}^1$ to a new cluster $\omega_1$. We initialize: $\eta(\omega_1)\gets 1;\ \bN^\Sigma(\omega_1) \gets \bN_{s+1}^1;\ \underline{\bN} (\omega_1)\gets \bN_{s+1}^1;\ \overline{\bN}(\omega_1)\gets \bN_{s+1}^1;\ \calX(\omega_1)\gets \{(\bN_s^1,\bx_s^1)\};\ c^\Sigma(\omega_1)\gets \chi_s^1$
\item We assign $\bN_{s+1}^q$ to a cluster based on its $\ell_\infty$ distance to the minimum and maximum states:
\begin{itemize}
    \item[--] If $\min_{\ell=1,\cdots,k}\max(\|\bN_{s+1}^q-\underline{\bN}(\omega_\ell)\|_\infty, \|\bN_{s+1}^q-\overline{\bN}(\omega_\ell)\|_\infty)\leq\varepsilon$, then $\bN^q$ is assigned to $\omega_{l^*}$ where $\ell^*$ denotes the index of the cluster that attains the minimum:
    \vspace{-6pt}
\begin{align*}
    &\eta(\omega_{l^*}) \gets \eta(\omega_{l^*}) + 1;\ \bN^\Sigma(\omega_{l^*}) \gets \bN^\Sigma(\omega_{l^*}) + \bN_{s+1}^q;\ \underline{\bN} (\omega_{l^*})\gets \min(\underline{\bN}(\omega_{l^*}),\bN_{s+1}^q)\\
    &\overline{\bN}(\omega_{l^*})\gets \max(\overline{\bN}(\omega_{l^*}),\bN_{s+1}^q);\ \calX(\omega_{l^*})\gets \calX(\omega_{l^*}) \cup \{(\bN_s^q, \bx_s^q)\};\  c^\Sigma(\omega_{l^*})\gets c^\Sigma(\omega_{l^*}) +  \chi_s^q,
\end{align*}
\item[--] Otherwise, $\bN_{s+1}^q$ is assigned to a new cluster $\omega_{k+1}$: $\eta(\omega_{k+1}) \gets 1;\ \bN^\Sigma(\omega_{k+1})\gets \bN_{s+1}^q;\  \underline{\bN} (\omega_{k+1})\gets \bN_{s+1}^q;\  \overline{\bN}(\omega_{k+1})\gets \bN_{s+1}^q;\ \calX(\omega_{k+1})\gets  \{(\bN_{s}^q, \bx_{s}^q)\};\  c^\Sigma(\omega_{k+1})\gets  \chi_s^q$
\end{itemize}
\item We retrieve each cluster's centroid $\bN_{s+1}(\omega_\ell)=\bN^\Sigma(\omega_\ell)/\eta(\omega_\ell)$ and define the cost function in the clustered state space, as follows: $c(\bN_s, \bx_{s})=\frac{c^\Sigma(\omega_\ell)}{\eta(\omega_\ell)}+ \Gamma_{s}(\bx_{s})$ for all $(\bN_s, \bx_s) \in \calX(\omega_\ell)$.
\end{enumerate}

We solve a pricing problem approximation in the clustered state space, via backward induction:
\begin{align}
    \bpi_s^*\left(\bN_s\right)&=\argmin_{\bx_s\in\calF_s}\left\{c(\bN_s, \bx_{s})-\sum_{j=1}^{m_s}\lambda_{sj}\bu_{sj}^\top\bx_s+J_{s+1}\left(\bN_{s+1}\right)\right\}\\
    J_s\left(\bN_s\right)&=c(\bN_s, \bpi_s^*\left(\bN_s\right))-\sum_{j=1}^{m_s}\lambda_{sj}\bu_{sj}^\top\bpi_s^*\left(\bN_s\right)+J_{s+1}\left(\bN_{s+1}\right) \\
    &\quad\text{where $\bM(\tau_s)\gets\bN_s$; $\frac{d\bM(t)}{dt}=f(\bM(t),\bpi_s^*\left(\bN_s\right)),\ \forall t\in\left[\tau_s,\tau_{s+1}\right]$; and $\bN_{s+1}\gets \bM(\tau_{s+1})$}\nonumber
\end{align}

Pseudo-code is given in Algorithm~\ref{alg:dp_clust} in Appendix~\ref{app:DP}. The algorithm makes a single pass through all the states, thus terminating in linear time (as long as the number of clusters is much smaller than the number of states). Still, the algorithm controls the maximum pair-wise distance within each cluster, guaranteeing an $\ell_\infty$-diameter within $\varepsilon$. This is formalized in Proposition~\ref{prop:cluster}.
\begin{proposition}\label{prop:cluster}
At each decision epoch $s=1,\cdots,S+1$, each cluster $\omega_1,\cdots,\omega_k$ satisfies:
\vspace{-6pt}
\begin{equation}
\|\overline{\bN}(\omega_\ell)-\underline{\bN}(\omega_\ell)\|_\infty\leq \varepsilon.
\end{equation}
\end{proposition}
This guarantee allows us to bound the global approximation error proportionally to $\varepsilon$, in Proposition \ref{prop:approx_error}. Note that the error grows exponentially in time, which is typical in dynamical systems \citep{sager2012integer}. Nonetheless, by controlling the cluster diameter, the state-clustering algorithm can control the global error. Our results in Section~\ref{subsec:DPresults} show that, in practice, our state-clustering algorithm induces moderate errors in small instances and can scale to large instances. Moreover, by controlling the cluster diameter, it leads to stronger optimization solutions than $k$-means.

\begin{proposition}\label{prop:approx_error}
Assume that the transition function $f(., \bx_{s})$ is $L_{s}$-Lipschitz in the state variables in the $l_\infty$ metric, for $s \in \{1,\cdots, S+1\}$. For each true state $\bN^*_{s} \in \calN^*_{s}$ (from Algorithm~\ref{alg:dp_full}), there exists a clustered state $\bN_{s} \in \calN_s$ (from Algorithm~\ref{alg:dp_clust}) such that:
\vspace{-6pt}
\begin{equation}
\|\bN^*_s - \bN_s\|_\infty\leq 
\begin{cases}
    0 & \text{if $s =1$} \\
\varepsilon\sum\limits_{\sigma=3}^{s+1} \exp\left(\sum\limits_{\nu=\sigma}^{s} L_{\nu-1}(\tau_{\nu} - \tau_{\nu-1})\right) ,&  \text{for all}\ s\geq 2
\end{cases}
\end{equation}
\end{proposition}

\subsection{A tri-partite branching disjunction}
\label{subsec:BP}

We embed the column generation procedure into a branch-and-bound structure to retrieve an integral solution to the $(\calS\calP)$ formulation. We proceed via multi-phase branching: we first branch on the mixed-integer variables $\by$; we then branch on the natural resource allocation variables $\bx_{is}$; and we finally restore integrality of the plan-based variables $z^p_i$ via a tri-partite branching disjunction.

\paragraph{Branching on mixed-integer variables.} In the case where the formulation involves discrete variables $\by$, we first branch on those variables whenever one of its components is fractional:
\vspace{-6pt}
\begin{equation}\label{eq:branchingY}
\underbrace{(y_\ell\leq\left\lfloor\widehat{y}_\ell\right\rfloor)}_{\text{left branch}}\lor\underbrace{(y_\ell\geq\left\lceil \widehat{y}_\ell\right\rceil)}_{\text{right branch}},\ \text{where $\widehat{y}_\ell$ denotes the $\ell^{\text{th}}$ component of the incumbent solution}
\end{equation}

\paragraph{Branching on the natural variables.} To avoid building deep and one-sided trees and to maintain the structure of the pricing problem, a typical approach in branch-and-price involves branching on the natural variables (i.e., the resource allocation decisions $\bx_{is}$ in our case) as opposed to branching on the composite plan-based variables (i.e., the variables $z^p_i$). Let $\bx_{is}(\bz)\in\R^{d_{is}}$ denote the resource allocation variable in segment $i=1,\cdots,n$ and epoch $s=1,\cdots,S$ for a plan-based solution $\bz$; and let $x^k_{is}(\bz)\in\R$ be its $k^{\text{th}}$ component, for $k=1,\cdots,d_{is}$. A column generation solution $\widehat{\bz}$ leads to infeasible resource allocations if $\bx_{is}(\widehat\bz)\notin\calF_{is}$. In that case, we exploit the discreteness of the feasible region $\calF_{is}$ to create a valid disjunction. Since $\calF_{is}$ does not necessarily comprise contiguous integers, we introduce dedicated floor and ceiling functions: $\lfloor \ba\rfloor^k_{is}=\max\{\beta\leq a^k:\exists\bb\in\calF_{is},b^k=\beta\}$ and $\lceil\ba\rceil^k_{is}=\min\{\beta\geq a^k:\exists\bb\in\calF_{is}:b^k=\beta\}$ (where $a^k$ and $b^k$ are the $k^{\text{th}}$ components of $\ba$ and $\bb$).

Armed with these notations, we can then define the following valid branching disjunction:
\vspace{-6pt}
\begin{equation}\label{eq:branching}
\underbrace{\bigg(x^k_{is}(\bz)\leq\left\lfloor\bx_{is}(\widehat{\bz})\right\rfloor^k_{is}\bigg)}_{\text{left branch}}\lor\underbrace{\bigg(x^k_{is}(\bz)\geq\left\lceil\bx_{is}(\widehat{\bz})\right\rceil^k_{is}\bigg)}_{\text{right branch}},\ \text{with}\ \bx_{is}(\widehat\bz)=\sum_{p\in\calP_i}\balpha^p_{is}\widehat{z}^p_i
\end{equation}

Out of these disjunctions, we select one that corresponds to a variable $x^k_{is}(\widehat\bz)$ with the largest value of $\min\left\{x^k_{is}(\widehat\bz)-\left\lfloor\bx_{is}(\widehat\bz)\right\rfloor^k_{is},\left\lceil\bx_{is}({\widehat\bz})\right\rceil^k_{is}-x^k_{is}(\widehat\bz)\right\}$. This branching strategy is an analog to the \textit{most fractional} branching strategy in integer optimization. Importantly, this disjunction preserves the structure of the pricing problem, by merely restricting the search from the full feasible regions $\calF_{is}$ to their subregions defined by the corresponding lower-bound and upper-bound constraints.

\paragraph{Tri-partite branching.}

By design, the above branching schemes enforce integrality of the $\by$ variables and ensure that the plan-based variables define feasible resource allocation decisions $\bx_{is}$. However, they may not guarantee integral plan-based variables $z^p_i$; for instance, we may obtain $\alpha^{p_1}_{is}=4$, $\alpha^{p_2}_{is}=8$, $z^{p_1}_i=z^{p_2}_i=0.5$, and $6\in\calF_{is}$. In linear problems, this solution can be brought into an equivalent feasible solution by considering the ``average plan'' $p^*$ with $\alpha^{p^*}_{is}=6$. In our problem, however, the non-convex system dynamics break this equivalence because $C^{p^*}_i\neq0.5C^{p_1}_i+0.5C^{p_2}_i$ in general. Thus, the plan $p^*$ is no longer guaranteed to form an optimal solution to Problem $(\calP)$.

A direct way to enforce the integrality of the plan-based variables would be to create disjunctions of the form $(z^p_i=0)\lor(z^p_i=1)$. However, as noted earlier, this disjunction can lead to weak and imbalanced tree structures. Moreover, it breaks the structure of the pricing problem---notably, it is difficult to seek the ``second-best'' solution in the branch associated with the $z^p_i=0$ disjunction.

Instead, we devise a novel tri-partite branching disjunction to handle the non-linearities. Consider a node with a solution $\widehat{\bz}$. Assume that $\bx_{is}(\widehat\bz)\in\calF_{is}$ for all $i,s$ but that there exists a fractional plan-based variable $z^p_i\in(0,1)$. Let $\alpha^{pk}_{is}(\bz)\in\R$ denote the $k^{\text{th}}$ component of $\balpha^p_{is}$ for $k=1,\cdots d_{is}$. There must exist a segment $i$, an epoch $s$, a component $k$ and a plan $p_0$ such that:
\vspace{-6pt}
\begin{equation}\label{eq:tripartiteINF}
\bx_{is}(\widehat{\bz})=\sum_{p\in\calP_i}\balpha^p_{is}\widehat{z}^p_i\in\calF_{is}\ \text{and}\ \alpha^{k,p_0}_{is}\neq x^k_{is}(\widehat{\bz})\ \text{and}\ \widehat{z}^{p_0}_i>0.
\end{equation}

Intuitively, $\bx_{is}(\widehat{\bz})$ defines a ``promising'' resource allocation decision for segment $i$ at epoch $s$. Accordingly, we create a corresponding branch in the tree. To retain a mutually exclusive and collectively exhaustive tree, we create two additional branches with lower and higher resource allocations. Formally, our tri-partite branching disjunction is defined via the following disjunction:
\begin{equation}\label{eq:tripartite}
\underbrace{\bigg(x^k_{is}(\bz)<x^k_{is}(\widehat{\bz})\bigg)}_{\text{left branch}}\lor\underbrace{\bigg(x^k_{is}(\bz)=x^k_{is}(\widehat{\bz})\bigg)}_{\text{middle branch}}\lor\underbrace{\bigg(x^k_{is}(\bz)>x^k_{is}(\widehat{\bz})\bigg)}_{\text{right branch}},\ \text{with}\ \bx_{is}(\widehat\bz)=\sum_{p\in\calP_i}\balpha^p_{is}\widehat{z}^p_i
\end{equation}

A final question lies in variable selection to determine which $x^k_{is}(\widehat{\bz})$ to branch upon. We extend the \textit{most fractional} principle, by selecting the variable with the largest weighted difference between the plan-based allocations and the resource allocations from the column generation solution. Formally, we select $x^k_{is}(\widehat{\bz})$ such that the $(i,k,s)$ tuple maximizes $\sum_{p \in \calP_i}\widehat{z}^p_i| \alpha^{pk}_{is}-x^k_{is}(\widehat{\bz})|$.

Again, this branching disjunction does not impact the structure of the pricing problem, which can easily accommodate lower-bounding and upper-bounding constraints on the decision variables. As we shall establish in Theorem~\ref{thm:exact}, this procedure eliminates any infeasible solution satisfying the conditions in Equation~\eqref{eq:tripartiteINF} and yields integer plan-based variables $z_i^p$ upon termination.

\subsection{Summary of the branch-and-price algorithm}\label{subsec:summary}

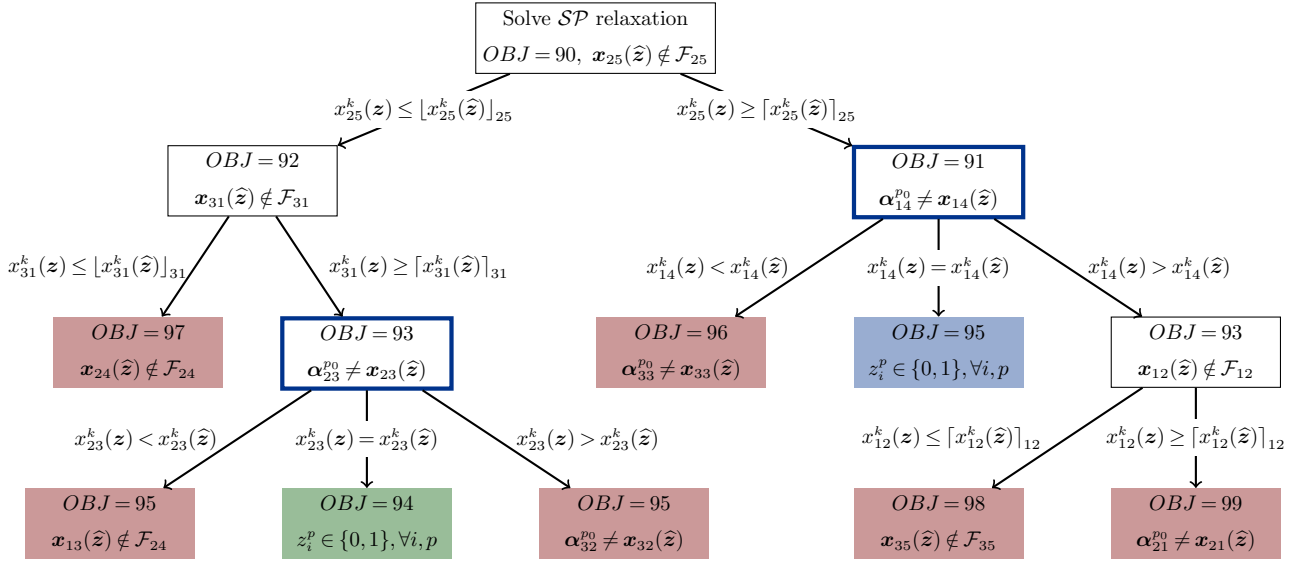
\begin{figure}[h!]
\centering\small
    \begin{tikzpicture}[scale=0.38, every node/.style={scale=0.8},every text node part/.style={align=center}]
        \tikzset{reg/.style={rectangle,draw, minimum width = 80pt, minimum height=22pt}}
        \tikzset{opt/.style={rectangle,fill=mygreen!40, minimum width = 80pt, minimum height=22pt}}
        \tikzset{bound/.style={rectangle,fill=myred!40, minimum width = 80pt, minimum height=22pt}}
        \tikzset{intx/.style={rectangle,draw=myblue,ultra thick, minimum width = 80pt, minimum height=22pt}}
        \tikzset{int/.style={rectangle,fill=myblue!40, minimum width = 80pt, minimum height=22pt}}
        
        \node[reg,anchor=north] (parent) at (0,17){Solve $\calS\calP$ relaxation\\$OBJ=90,\ \bx_{25}(\widehat{\bz})\notin\calF_{25}$};
        
        \node[reg,anchor=north] (subproblem1) at (-12,12){$OBJ=92$\\$\bx_{31}(\widehat{\bz})\notin\calF_{31}$};
        \draw[->,thick] (parent) -- node[fill=white,midway]{$x^k_{25}(\bz)\leq\left\lfloor x^k_{25}(\widehat{\bz})\right\rfloor_{25}$} (subproblem1);
        
        \node[intx,anchor=north] (subproblem2) at (12,12){$OBJ=91$\\$\balpha^{p_0}_{14}\neq \bx_{14}(\widehat{\bz})$};
        \draw[->,thick] (parent) -- node[fill=white,midway]{$x^k_{25}(\bz)\geq\left\lceil x^k_{25}(\widehat{\bz})\right\rceil_{25}$} (subproblem2);
        
        \node[bound,anchor=north] (subproblem11) at (-16,6){$OBJ=97$\\$\bx_{24}(\widehat{\bz})\notin\calF_{24}$};
        \draw[->,thick] (subproblem1) -- node[left,midway]{$x^k_{31}(\bz)\leq\left\lfloor x^k_{31}(\widehat{\bz})\right\rfloor_{31}$} (subproblem11);
        
        \node[intx,anchor=north] (subproblem12) at (-8,6){$OBJ=93$\\$\balpha^{p_0}_{23}\neq \bx_{23}(\widehat{\bz})$};
        \draw[->,thick] (subproblem1) -- node[right,midway,xshift=0.2 cm]{$x^k_{31}(\bz)\geq\left\lceil x^k_{31}(\widehat{\bz})\right\rceil_{31}$} (subproblem12);

        \node[bound,anchor=north] (subproblem21) at (3,6){$OBJ=96$\\$\balpha^{p_0}_{33}\neq \bx_{33}(\widehat{\bz})$};
        \draw[->,thick] (subproblem2) -- node[left,midway,xshift=-0.2 cm]{$x^k_{14}(\bz)<x^k_{14}(\widehat{\bz})$} (subproblem21);
        
        \node[int,anchor=north] (subproblem22) at (12,6){$OBJ=95$\\$z^p_i\in\{0,1\},\forall i,p$};
        \draw[->,thick] (subproblem2) -- node[fill=white,midway]{$x^k_{14}(\bz)=x^k_{14}(\widehat{\bz})$} (subproblem22);
        
        \node[reg,anchor=north] (subproblem23) at (21,6){$OBJ=93$\\$\bx_{12}(\widehat{\bz})\notin\calF_{12}$};
        \draw[->,thick] (subproblem2) -- node[right,midway,xshift=0.2 cm]{$x^k_{14}(\bz)>x^k_{14}(\widehat{\bz})$} (subproblem23);
        
        \node[bound,anchor=north] (subproblem121) at (-17,0){$OBJ=95$\\$\bx_{13}(\widehat{\bz})\notin\calF_{24}$};
        \draw[->,thick] (subproblem12) -- node[left,midway,xshift=-0.2 cm]{$x^k_{23}(\bz)<x^k_{23}(\widehat{\bz})$} (subproblem121);
        
        \node[opt,anchor=north] (subproblem122) at (-8,0){$OBJ=94$\\$z^p_i\in\{0,1\},\forall i,p$};
        \draw[->,thick] (subproblem12) -- node[fill=white,midway]{$x^k_{23}(\bz)=x^k_{23}(\widehat{\bz})$} (subproblem122);
        
        \node[bound,anchor=north] (subproblem123) at (1,0){$OBJ=95$\\$\balpha^{p_0}_{32}\neq \bx_{32}(\widehat{\bz})$};
        \draw[->,thick] (subproblem12) -- node[right,midway,xshift=0.2 cm]{$x^k_{23}(\bz)>x^k_{23}(\widehat{\bz})$} (subproblem123);
        
        \node[bound,anchor=north] (subproblem231) at (12,0){$OBJ=98$\\$\bx_{35}(\widehat{\bz})\notin\calF_{35}$};
        \draw[->,thick] (subproblem23) -- node[left,midway,xshift=-0.3cm]{$x^k_{12}(\bz)\leq\left\lceil x^k_{12}(\widehat{\bz})\right\rceil_{12}$} (subproblem231);

        \node[bound,anchor=north] (subproblem232) at (21,0){$OBJ=99$\\$\balpha^{p_0}_{21}\neq \bx_{21}(\widehat{\bz})$};
        \draw[->,thick] (subproblem23) -- node[fill=white,midway]{$x^k_{12}(\bz)\geq\left\lceil x^k_{12}(\widehat{\bz})\right\rceil_{12}$} (subproblem232);
        
    \end{tikzpicture}
\caption{Illustration of the branching structure, combining bi-partite branching and tri-partite branching. Red nodes are pruned by bound; blue nodes are pruned by feasibility; white nodes trigger bi-partite branching; blue squares trigger tri-partite branching; the green node indicates the optimal solution.}
\label{fig:branching}
\end{figure}

Algorithm~\ref{alg:exact} summarizes our solution approach. In each node, the algorithm solves the linear relaxation of $\calS\calP$ via column generation (Step~1), iterating between the RMP (Step~1.1) and the dynamic programming algorithm for the pricing problem (Step~1.2), until convergence. We also apply in Step~2 an optional upper-bounding scheme for acceleration: in our resource allocation problems for instance, we solve a restricted master problem with integrality constraints; in the vaccination centers problem, this becomes much more cumbersome, so we fix facility construction variables based on the linear relaxation solution and solve the subsequent vaccine allocation problem. The algorithm then proceeds to bi-partite branching when the variables $\by$ are not integral or when the natural resource allocation variables $\bx$ are infeasible (Steps 3 and 4), and to tri-partite branching to restore integrality in plan-based variables (Step~5). This overall branching strategy is illustrated in Figure~\ref{fig:branching}. Step~6 corresponds to the node selection step---we implemented breadth-first and depth-first search strategies, and found comparable performance. The algorithm uses bounding and feasibility rules to prune leaves (Step~2) and terminates when no remaining leaf is active (Step~6).

\begin{algorithm}[h!]
\floatname{algorithm}{Algorithm}
\caption{Branch-and-price algorithm for Problem ($\calP$).}
\renewcommand{\arraystretch}{1.0}\small
\label{alg:exact}
\begin{algorithmic}[1] 
\State \textbf{Initialization:} node list $\Phi=\{0\}$ of the branch-and-bound tree; parent node $cur=0$; initial sets of plans $\calP^0_1,\cdots,\calP^0_n$ and corresponding RMP \eqref{MP:obj}--\eqref{MP:domain}; upper bound $UB=+\infty$; incumbent solution $\bz^*=\emptyset$.

\State \textbf{\textit{Step~1.}} Solve parent node $cur$ via column generation (Section~\ref{subsec:CG}); store objective $OBJ$, solution $\widehat{\bz}$, and corresponding resource allocation decisions $\bx_{is}(\bz)=\sum_{p\in\calP_i}\balpha^p_{is}\widehat{z}^p_i$; remove $cur$ from node list $\Phi$.

\State ~~\textbf{\textit{Step~1.1.}} Solve RMP relaxation with $\calP^0_1,\cdots,\calP^0_n$; store solution $\widehat{\bz}$, objective $OBJ$, and dual costs $\blambda,\bmu$.

\State ~~\textbf{\textit{Step~1.2.}} Solve $(PP_i)$ for each segment $i=1,\cdots,n$, using the exact dynamic programming algorithm (Algorithm~\ref{alg:dp_full}) or the state-clustering dynamic programming algorithm (Algorithm~\ref{alg:dp_clust}). If the objective is negative for some segment $i$, add optimal plan to $\calP^0_i$, and go to Step~1.1.

\State \textbf{\textit{Step~2.}} Apply  upper-bounding heuristic (e.g.,  restricted master problem with integrality constraints, variables fixing); store objective $\overline{OBJ}$ and solution $\overline{z}$. If $\overline{OBJ}<UB$, update $UB\gets\overline{OBJ}$, $\bz^*\gets\overline{\bz}$.

\State \hspace{12pt}-- If $OBJ<UB$ and $\widehat y_\ell\notin\Z$ for some $\ell=1,\cdots,q$; go to Step~2.

\State \hspace{12pt}-- If $OBJ<UB$ and $\bx_{is}(\widehat\bz)\notin\calF_{is}$ for some $i,s$; go to Step~3.

\State \hspace{12pt}-- If $OBJ<UB$ and $\bx_{is}(\widehat\bz)\in\calF_{is},\ \forall i,s$ and $z^p_i$ is fractional for some $i,p$; go to Step~4.

\State \hspace{12pt}-- If $OBJ<UB$ and $\bx_{is}(\widehat\bz)\in\calF_{is},\ \forall i,s$ and  $\widehat{\bz}$ is integral; update $UB\gets OBJ$, $\bz^*\gets\widehat{\bz}$; go to Step~5.

\State \hspace{12pt}-- If $OBJ\ge UB$; go to Step~5.

\State  \textbf{\textit{Step~3:} Branching on $\by$ variables.} Add two children nodes to $\Phi$ (Equation \eqref{eq:branchingY}); go to Step~5.

\State  \textbf{\textit{Step~4:} Bi-partite branching.} Add two children nodes to $\Phi$ (Equation \eqref{eq:branching}); go to Step~5.

\State  \textbf{\textit{Step~5:} Tri-partite branching.}  Add three children nodes to $\Phi$ (Equation~\eqref{eq:tripartite}); go to Step~5.

\State \textbf{\textit{Step~6.}} If $\Phi = \emptyset$, \texttt{STOP}; return solution $\bx^*$ and objective $UB$. Otherwise, choose $cur\in\Phi$; go to Step~1.
\end{algorithmic}
\end{algorithm}

Theorem~\ref{thm:exact} establishes the convergence and exactness of the algorithm, as long as the pricing problem is solved to optimality (Step~1.2). Exactness stems from the facts that the algorithm maintains valid lower and upper bounds. Finiteness follows from the fact that the branching disjunction eliminates any infeasible solution and from the finiteness of the decision space. Note that the pricing problem only needs to be solved to optimality upon convergence; therefore, an exact acceleration involves applying the state-clustering dynamic programming approximation in initial iterations and then turning to the exact algorithm. In our experiments, we merely apply the state-clustering approximation to retain tractability, but we evaluate the performance of the optimized solution in view of the full dynamical system---not the clustered approximation (Section~\ref{subsec:impact}).

\begin{theorem}\label{thm:exact}
The branch-and-price algorithm with the tri-partite branching disjunction converges in a finite number of iterations and returns an optimal solution to the $\calS\calP$ formulation.
\end{theorem}
\section{Experimental results}
\label{sec:results}

We evaluate our methodology on the four problems described in Section~\ref{sec:model}. All experiments were run on Julia v1.5.2 with the JuMP package \citep{dunning2017jump} on 10-core i9-10900k CPU, with a 3-hour limit and a 0.1\% tolerance. All instances and code are made available online.\footnote{\hyperlink{https://github.com/martinrame24/Prescriptive-Contagion-Analytics}{https://github.com/martinrame24/Prescriptive-Contagion-Analytics}}

\subsection{Benefits of the state-clustering dynamic programming algorithm}\label{subsec:DPresults}

Our methodology hinges on the ability to solve the pricing problem efficiently. We first compare our state-clustering dynamic programming algorithm (Algorithm~\ref{alg:dp_clust}) to the exact dynamic programming algorithm based on state space enumeration (Algorithm~\ref{alg:dp_full}). Table~\ref{tab:DP} reports computational times, broken down into initialization times to create the state space (``init.'') and dynamic programming times (``DP'). For all instances solved by the exact algorithm, we report the Median Average Percentage Error (MAPE) and Median Absolute Error (MAE) between clustered and true states.

\begin{table}[h!]
\centering
\footnotesize\renewcommand{\arraystretch}{1.0}
\caption{State-clustering dynamic programming against exhaustive enumeration (vaccine allocation problem).}
\label{tab:DP}
\begin{tabular}{cccllccccccccccc}
\toprule
 & & & & & & \multicolumn{2}{c}{\textbf{Time (sec.)}} & \multicolumn{4}{c}{\textbf{MAPE (\%)}} & \multicolumn{4}{c}{\textbf{MAE ($\times 10^{-4}$)}} \\ \cmidrule(lr){7-8}\cmidrule(lr){9-12}\cmidrule(lr){13-16}
$n$ & $S$ & $D$ & \multicolumn{1}{c}{Method} & \multicolumn{1}{c}{Tolerance} & $|S|$ & Init. & DP & $s=1$ & $s=2$ & $s=3$ & $s=4$ & $s=1$ & $s=2$ & $s=3$ & $s=4$ \\\hline 
51 & 4 & 6 & Enumeration &\multicolumn{1}{c}{|}& 79,305 & 20.01 &  1.33 &|&|&|&|&|&|&|&| \\
&& & Clustering & $\varepsilon=0.002$ &2,350 & 3.33 &  0.08 & 0.00&0.00&0.48&0.69 & 0.00&0.00&0.29&0.27  \\
&& & Clustering & $\varepsilon=0.005$ &1,458 & 2.20 &  0.04 & 0.00&0.00&1.32&2.08 & 0.00&0.00&1.04&0.95  \\
&& & Clustering & $\varepsilon=0.01$ & 982 & 1.57 &  0.03 & 0.00&0.01&7.22&9.38 & 0.00&0.01&3.72&4.02 \\\hline
51 & 4 & 11 & Enumeration &\multicolumn{1}{c}{|}& 821K & 194 &  17.1 &|&|&|&|&|&|&| \\
&& & Clustering & $\varepsilon=0.002$ &4,971 & 13.3 &  0.33 & 0.00&0.08&0.57&0.86 & 0.00&0.08&0.37&0.36  \\
&& & Clustering & $\varepsilon=0.005$ & 2,533 & 6.74 &  0.19 & 0.00&0.12&2.86&3.67 & 0.00&0.29&1.59&1.42  \\
&& & Clustering & $\varepsilon=0.01$ & 1,521& 4.29 &  0.10 & 0.00&0.11&6.27&10.02 & 0.00&0.29&3.58&4.23 \\\hline
51 & 4 & 21 & Enumeration &\multicolumn{1}{c}{|}& 10.4M & 2,950 &  234 &|&|&|&|&|&|&|&| \\
&& & Clustering & $\varepsilon=0.002$ &10,615& 74.6&  1.09 & 0.00&0.08&0.87&1.08 & 0.00&0.08&0.52&0.41 \\
&& & Clustering & $\varepsilon=0.005$ & 4,477& 26.9 &  0.58& 0.00 & 0.09 & 3.51 & 5.12 & 0.00 &0.01&1.86&1.57 \\
&& & Clustering & $\varepsilon=0.01$ & 2,473 & 14.8  &  0.29 & 0.00 & 0.22 & 7.17 &  13.64 & 0.00 &0.28&3.78&4.91 \\\toprule
51 & 6 & 21 & Enumeration &\multicolumn{1}{c}{|}& 4.59B & n/a &  n/a &|&|&|&|&|&|&|&| \\
&& & Clustering & $\varepsilon=0.002$ & 29,020 & 370&  7.41 & n/a & n/a & n/a &  n/a & n/a &  n/a & n/a &  n/a  \\
&& & Clustering & $\varepsilon=0.005$ & 9,212 & 75.0 &  1.48 & n/a & n/a & n/a &  n/a & n/a &  n/a & n/a &  n/a  \\
&& & Clustering & $\varepsilon=0.01$ & 4,328 & 31.3 &  0.55 & n/a &  n/a & n/a &  n/a & n/a &  n/a & n/a &  n/a \\\hline
51 & 8 & 21 & Enumeration &\multicolumn{1}{c}{|}& 2.06T & n/a &  n/a &|&|&|&|&|&|&|&| \\
&& & Clustering & $\varepsilon=0.002$ & 56,635 & 1,057 &  13.6  & n/a & n/a & n/a &  n/a & n/a &  n/a & n/a &  n/a  \\
&& & Clustering & $\varepsilon=0.005$ & 15,076 & 147 &  3.50 & n/a & n/a & n/a &  n/a & n/a &  n/a & n/a &  n/a  \\
&& & Clustering & $\varepsilon=0.01$ & 6,341 & 52.4 &  1.47 & n/a & n/a & n/a &  n/a & n/a &  n/a & n/a &  n/a \\\hline
51 & 10 & 21 & Enumeration &\multicolumn{1}{c}{|}& 893T & n/a &  n/a &|&|&|&|&|&|&|&| \\
&& & Clustering & $\varepsilon=0.002$ & 89,685 & 2,586 &  32.8 & n/a & n/a & n/a &  n/a & n/a &  n/a & n/a &  n/a  \\
&& & Clustering & $\varepsilon=0.005$ & 21,308 & 238 &  5.50 & n/a & n/a & n/a &  n/a & n/a &  n/a & n/a &  n/a  \\
&& & Clustering & $\varepsilon=0.01$ & 8,353 & 77.9 &  2.13 & n/a & n/a & n/a &  n/a & n/a &  n/a & n/a &  n/a \\\bottomrule
\end{tabular}
\begin{tablenotes}
    \vspace{-6pt}
    \item ``n/a'' means that the exact dynamic programming algorithm does not terminate due to memory limitations.
\end{tablenotes}
    \vspace{-12pt}
\end{table}

Note, first, that the number of states grows very quickly, reflecting the ``curse of dimensionality'' in dynamic programming \citep{powell2022reinforcement}. Even by exploiting the segment-based decomposition in the pricing problem, the state space scales in $\calO(n\times D^{S})$, with trillions of states in our largest instances. Exhaustive enumeration does not scale to even small instances, requiring minutes to converge with 4 segments and 4 decision epochs. The state-clustering algorithm considerably reduces the state space---by up to a factor of $10^{10}$ to $10^{11}$. Thus, it accelerates convergence by two orders of magnitude in small instances (seconds versus minutes), and scales to the largest instance with 51 states, 10 epochs and 21 decisions in minutes when the exact algorithm fails to terminate.

Moreover, the state-clustering algorithm provides a high-quality approximation of the full state space. With $\varepsilon=0.01$, the algorithm results in a median absolute error up to 4.9$\times10^{-4}$ and a median absolute percentage error up to 13.6\%. A tolerance of $\varepsilon=0.002$ further reduces the MAE within 6$\times10^{-5}$ and the MAPE within 1.1\%. As expected, the quality of the approximation deteriorates over time due to the propagation of errors. Whereas we cannot estimate the error over longer time horizons because of the limitations of the enumerative algorithm, this observation could motivate dynamic implementations of our methodology, for instance, by re-evaluating the system's dynamics and re-optimizing resource allocations at each epoch. Nonetheless, these results suggest that the drift in state approximation remains moderate, underscoring the critical role of the state-clustering algorithm to enhance the tractability of the pricing problem at small costs in terms of accuracy.

Finally, we provide in Appendix~\ref{app:clustering} a detailed comparison of our state-clustering algorithm (Algorithm~\ref{alg:dp_clust}) and a $k$-means benchmark. By design, the $k$-means algorithm leads to a smaller average distance but a larger worst-case distance within each cluster. Per Proposition~\ref{prop:approx_error}, this worst-case error can propagate over time in dynamical systems. In turn, our algorithm results in much higher fidelity across the time horizon, with significant benefits in the downstream optimization---achieving a 1.4\% to 25\% cost reduction against the $k$-means benchmark for the vaccine allocation problem.

\subsection{Benefits of the branch-and-price algorithm}
\label{subsec:Benefits}

Armed with the state-clustering algorithm, we now solve our four prescriptive contagion analytics problems. Table~\ref{tab:BP} reports computational times, broken down into the state-clustering initialization, the restricted master problems, and the pricing problems. It also reports the solution and the lower bound obtained via column generation, via bi-partite branching (i.e., Algorithm~\ref{alg:exact} without Step 5), and with the full branch-and-price algorithm with tri-partite branching. The column generation and bi-partite branching solutions are obtained by solving the restricted master problem with integrality constraints upon convergence. For brevity, the table reports results for ``easy'', ``medium'' and ``hard'' instances. Full results in Appendix~\ref{app:results} show the robustness of these findings.

\begin{table}[h!]
\centering
\renewcommand{\arraystretch}{1.0}
\footnotesize\caption{Performance of the branch-and-price algorithm.}
\label{tab:BP}
\begin{adjustbox}{width=\textwidth}
\begin{tabular}{lllcccccccc}
\toprule
  & & & & \multicolumn{4}{c}{CPU times (s)}  & \multicolumn{3}{c}{Solution quality} \\ \cmidrule(lr){5-8} \cmidrule(lr){9-11}
Problem & Instance & \multicolumn{1}{c}{Method} & Nodes & Init. & RMP & PP & Total & Upper bound & Solution & Gap \\\hline
Vaccine & Easy &  Column Generation &          1 &     21 &    0.02 &     0.22 &   0.29 &   93,695 &   92,261 &   1.53\% \\
allocation&  & Bi-partite B\&P &       28 &     21 &     0.27 &     1.43 &     7.11 &   93,379 &   93,266 &     0.12\% \\
&  &  Tri-partite B\&P &       46 &     21 &     0.42 &     2.2 &     9.6 &   93,353 &   93,263 &     0.09\% \\\cmidrule{2-11}
 &  Medium &  Column Generation &          1 &     43 &    0.022 &    0.81 &   0.92 &  156,465 &  155,624 &   0.54\% \\
&  &  Bi-partite B\&P &     4,778 &     43 &   53 &   426 &  1,064 &  156,183 &  156,026 &     0.1\% \\
&  &  Tri-partite B\&P &     4,778 &     43 &    55 &   420 &  1,058 &  156,183 &  156,027 &     0.1\% \\\cmidrule{2-11}
 &  Hard &  Column Generation &          1 &    144 &     0.03 &    4.07 &   4.32 &  225,323 &  222,592 &   1.21\% \\
&  &  Bi-partite B\&P &     9,990 &    144 &  101 &   1,765 &  3,475 &  225,048 &  224,495 &     0.24\% \\
& &  Tri-partite B\&P &     9,990 &    144 &  103 &  1,808 &  3,535 &  225,048 &  224,495 &     0.24\% \\\toprule
Vaccination & Easy, $K=2$ &  Column Generation &        1 &       36 &       0.05 &       0.03 &    0.46 &    786 &     652 &     29.15\% \\
centers & & Bi-partite B\&P &     2,919 &   36 &      77 &      11 &     608 &   764 &   763 &        0.09\% \\
& & Tri-partite B\&P &     2,919 &   36 &      75 &      11 &     604 &   764 &   763 &        0.09\% \\\cmidrule{2-11}
 & Easy, $K=3$ &  Column Generation &        1 &       43 &       0.02 &       0.03 &    0.57 &    678 &     652 &      3.75\% \\
 & & Bi-partite B\&P &    56,273 &  43 &    1,272 &     142 &   10,809 &   675 &   674 &        0.23\% \\
& & Tri-partite B\&P &    58,496 &  43 &    1,382 &     137 &   10,801 &   675 &   674 &        0.23\% \\\cmidrule{2-11}
& Hard, $K=2$ &  Column Generation &        1 &       45 &       0.19 &       0.43 &    1.11 &   8200 &    6660 &     29.91\% \\
& & Bi-partite B\&P &    11,203 &   45 &     216 &     147 &   10,808 &  8,134 &  7,912 &        2.73\%\\
& & Tri-partite B\&P &    11,899 &   45 &     224 &     178 &   10,808 &  8,132 &  7,917 &        2.64\%\\\cmidrule{2-11}
& Hard, $K=3$ &  Column Generation &        1 &       47 &       0.27 &       0.31 &    1.11 &   8,200 &    6,822 &     19.07\% \\
& &  Bi-partite B\&P &    15,982 &  47 &     356 &     207 &   10,885 &  8,048 &  7,714 &        4.14\% \\
& &  Tri-partite B\&P &    16,358 &  47 &     363 &     209 &   10,955 &  8,048 &  7,716 &        4.13\% \\\toprule
Content & Easy & Column Generation &         1 &       79 &      0.04 &      4.6 &    4.8 &  1.013 &   0.98 &    3.86\% \\
promotion& & Bi-Partite B\&P &      696 &           79 &  9.3 &  198 &  327 &  1.013 &  1.009 &      0.53\% \\
& & Tri-partite B\&P &      142 &           79 &  1.6 &   52 &   74 &  1.013 &  1.0125 &     0.09\% \\\cmidrule{2-11}
& Medium & Column Generation &         1 &       58 &      0.04 &      5.01 &   5.4 &  1.52 &   1.44 &   5.66\% \\
& & Bi-Partite B\&P &    17,960 &           58 &  294 &  1,949 &  7,201 &  1.52 &  1.51 &     0.83\% \\
& & Tri-partite B\&P &    17,819 &           58 &  3,380 &  2,0570 &  7,201 &  1.52 &  1.517 &     0.31\% \\\cmidrule{2-11}
& Hard & Column Generation &         1 &      134 &      0.08 &       15 &    16 &  2.02 &   1.87 &   7.08\% \\
& & Bi-Partite B\&P &    12,606 &          134 &  229 &  3,221 &  10,801 &  2.009 &  1.99 &     1.04\% \\
& & Tri-partite B\&P &    11,864 &          134 &  229 &  3,279 &  10,801 &  2.009 &  1.99 &     1.04\% \\\toprule
Congestion &  Easy &  Column Generation &          1 &    10 &    0.012 &    0.08 &   0.11 &   0.82 &     0.81 &   1.23\% \\
mitigation& & Bi-partite B\&P &        2 &    10 &   0.018 &    0.126 &    0.56 &  0.82 &  0.81 &     0.16\% \\
&  & Tri-partite B\&P &        8 &    10 &   0.053 &     0.23 &      0.8 &  0.82 &  0.82 &     0\% \\\cmidrule{2-11}
 & Medium  &  Column Generation &          1 &   185 &    0.024 &    6.2 &   6.27 &   2.20 &     2.05 &   6.61\% \\
&  & Bi-partite B\&P &       48 &   185 &   0.34 &   32 &   34 &  2.19 &  2.19 &     0.21\% \\
&  & Tri-partite B\&P &       57 &   185 &   0.43 &   39 &    42 &  2.19 &  2.19 &     0.06\% \\\cmidrule{2-11}
 &  Hard &  Column Generation &          1 &   459 &    0.044 &   20 &  20 &   3.32 &     2.93 &  11.9\% \\
& & Bi-partite B\&P &      206 &   459 &   4.14 &  350 &  389 &  3.31 &  3.31 &     0.2\% \\
&  & Tri-partite B\&P &      243 &   459 &   4.9 &  396 &   443 &  3.31 &  3.31 &     0.09\% \\
\bottomrule
\end{tabular}
\end{adjustbox}
\begin{tablenotes}
    \vspace{-6pt}
    \item Vaccine allocation: $n=51$; $B=2.5$M; $S=8$, $D=6$ (``easy''); $S=10$, $D=11$ (``medium''); $S=12$, $D=21$ (``hard'').\vspace{-3pt}
    \item Vaccination centers: Group H (``easy'') and Group A (``hard''), with $K=2$ or $K=3$ facilities per group, $S=6$.\vspace{-3pt}
    \item Content promotion: $n=20$; $K=2$ ; $S=6$, $D=21$,  (``easy''); $S=8$, $D=11$,  (``medium''); $S=10$, $D=21$ (``hard'').\vspace{-3pt}
    \item Congestion mitigation: $n=5$; $B^1=6$ and $B^2=4$ ($D=35$); $S=2$ (``easy''); $S=4$ (``medium''); $S=6$ (``hard'').\vspace{-3pt}
\end{tablenotes}
\vspace{-10pt}
\end{table}

\begin{table}[h!]
\centering
\footnotesize\renewcommand{\arraystretch}{1.0}
\caption{Detailed results for the vaccination centers problem (full United States, 3 facilities per group)}
\label{tab:BP_full_country}
\begin{tabular}{ccccclcccccccc}
\toprule
  & & & & & & & \multicolumn{4}{c}{CPU times (s)}  & \multicolumn{3}{c}{Solution quality}  \\ \cmidrule(lr){8-11} \cmidrule(lr){12-14}
Group & n & S & D & F & \multicolumn{1}{c}{Method} & Nodes & Init. & RMP & PP & Total & Upper bound & Solution & Gap \\\hline
A &  7 &  6 &  12 &  3  &  Column Generation &        1 &       47 &       0.27 &       0.31 &    1.11 &   8,200 &    6,822 &     19.07\% \\
& & & & &  Tri-partite B\&P &    16,358 &  47 &     363 &     209 &   10,955 &  8,048 &  7,716 &        4.13\% \\\hline
B &  7 &  6 &  10 &  3  &  Column Generation &        1 &       43 &       0.06 &       0.34 &    0.88 &   3,659 &    3,322 &      9.20\% \\
& & & & &  Tri-partite B\&P &    11,965 &  43 &     329 &     134 &   10,806 &  3,569 &  3,467 &        2.87\% \\\hline
C &  4 &  6 &  13 &  3 &  Column Generation &        1 &       38 &       0.04 &       0.11 &    0.55 &   3,880 &    3,171 &     18.28\% \\
& & & & &  Tri-partite B\&P &    10,101 &  38 &     369 &      98 &   10,812 &  3,880 &  3,692 &        4.84\% \\\hline
D &  4 &  6 &   5 &  3  &  Column Generation &        1 &       38 &       0.04 &       0.04 &    0.48 &   1,056 &     733 &     30.65\% \\
& & & & &  Tri-partite B\&P &    22,600 &  38 &     539 &      45 &   10,824 &  1,052 &  1,000 &        4.95\% \\\hline
E &  6 &  6 &  16 &  3  &  Column Generation &        1 &       55 &       0.08 &       0.51 &    1.02 &   8,331 &    6,927 &     16.86\% \\
& & & & & Tri-partite B\&P &     7,824 &  55 &     214 &     199 &   10,832 &  8,316 &  8,053 &        3.16\% \\\hline
F &  5 &  6 &  13 &  3 &  Column Generation &        1 &       47 &       0.07 &       0.25 &    0.83 &   3,593 &    2,466 &     45.25\% \\
& & & & & Tri-partite B\&P &    16,882 &  47 &     698 &     181 &   10,822 &  3,530 &  3,328 &         5.74\%\\\hline
G &  4 &  6 &   4 &  3 &  Column Generation &        1 &       42 &       0.02 &       0.02 &    0.50 &   1,985 &    1,679 &     15.39\% \\
& & & & &  Tri-partite B\&P &    39,651 &  42 &     609 &      72 &   10,802 &  1,955 &  1,921 &        1.76\% \\\hline
H &  6 &  6 &   4 &  3  &  Column Generation &        1 &       43 &       0.02 &       0.03 &    0.57 &    678 &     652 &      3.75\% \\
& & & & & Tri-partite B\&P &    58,496 &  43 &    1,382 &     137 &   10,801 &   675 &   674 &        0.23\% \\\hline
I &  4 &  6 &  20 &  3 &  Column Generation &        1 &       39 &       0.07 &       0.23 &    0.68 &   3,573 &    2,587 &     34.89\% \\
& & & & & Tri-partite B\&P &    55,864 &  39 &    4,357 &     507 &   10,871 &  3,423 &  3,399 &        0.72\% \\\hline
J &  4 &  6 &   5 &  3 &  Column Generation &        1 &       36 &       0.03 &       0.03 &    0.47 &    626 &     565 &      9.73\% \\
& & & & &  Tri-partite B\&P &    25,016 &  36 &     327 &      41 &   10,803 &   624 &   565 &        9.46\% \\\hline
Full &  51 &  6 &   --- &  3 &  Column Generation &       --- &          --- &         --- &         --- &         --- &  35,581 &  28,924 &      18.7\% \\
& & & & & Tri-partite B\&P &       --- &          --- &            --- &           --- &             --- &  35,072 &  33,815 &        3.58\% \\
\bottomrule
\end{tabular}
\end{table}

Note, first, that the column generation algorithm is highly scalable for the three resource allocation problems: vaccine allocation, content promotion, and congestion mitigation. Indeed, column generation solves the set partitioning relaxation in seconds even in the harder instances. The restricted master problem is virtually instantaneous but the pricing problem is more time-consuming due to the non-linear system dynamics---reinforcing the need for an efficient dynamic programming algorithm (Section~\ref{subsec:DPresults}). In fact, the state-clustering algorithm shifts the complexity away from the dynamic programming algorithm itself---enabling the online column generation algorithm to terminate in seconds---toward the offline state space clustering---which can take up to minutes.

Moreover, the column generation algorithm can derive high-quality solutions, thanks to the tight set partitioning formulation for these problems. In the vaccine allocation problem in particular, the column generation algorithm yields solutions within a 1-2\% optimality gap. For the content promotion and congestion mitigation problems, column generation generally performs well but can also leave a larger optimality gap---due to the sparsity constraint in content promotion and to the coordination of prevention and treatment resources in congestion mitigation. This limitation motivates the branch-and-price algorithm to further improve solutions and tighten optimality gaps.

Turning to our main observation, the branch-and-price algorithm solves every instance of the vaccine allocation, content promotion and congestion mitigation problems to near-optimality. In fact, it reaches optimality within a 0.1\% tolerance in most instances, and leaves a small optimality gap otherwise. The number of nodes generally remains moderate but can grow larger for harder instances, leading to higher computational requirements. As a result, the branch-and-price algorithm terminates in minutes to hours, but can still lead to strong solution improvements from column generation. These benefits can go up to 6-12\% in the content promotion and congestion mitigation cases. Even in vaccine allocation, where column generation leaves a small gap, the branch-and-price solution yields a 1\% improvement, amounting to an extra 2,000 lives saved over three months.

Finally, note that the standard bi-partite branching scheme is not sufficient to guarantee convergence to an optimal solution. Notably, it leaves a 1\% optimality gap in content promotion instances and a 0.2\% gap in congestion mitigation instances. This result underscores the impact of system non-linearities on the branch-and-price algorithm. Instead, the tri-partite branching scheme developed in this paper can be instrumental to ensure convergence to an optimal solution, while retaining an efficient branching tree and an efficient pricing problem structure.

Next, recall that the vaccination centers problem features a joint facility location and resource allocation structure, leading to a looser set-partitioning relaxation. As a result, the column generation heuristic leaves large optimality gaps up to 30\%. Restoring integrality via branch-and-price is much more challenging. Interestingly, the restricted master problem ends up being more time-consuming than the pricing problem, thus highlighting the joint complexities from coupled resource allocation and from non-linear system dynamics. Yet, the branch-and-price algorithm leads to significant solution improvements and gap reductions: it returns an optimal solution in the ``easy'' two-facility instance, a near-optimal solution in the ``hard'' two-facility instances, and moderate optimality gaps in the three-facility cases (0-4\% versus 3-30\% with column generation).

To shed more light into these findings, Table~\ref{tab:BP_facility} reports results for the all ten groups constituting the full United States. Despite variations in the problem's complexity, our core observations are highly robust: the column generation algorithm leaves consistently large optimality gaps, and the branch-and-price algorithm significantly improves the solution and tightens the gap. Ultimately, the branch-and-price algorithm results in an estimated 5,000 extra lives saved across the country over a 6-week period---a 17\% improvement over the column generation solution.

\subsection{Practical impact of the methodology}\label{subsec:impact}

Tables~\ref{tab:vaccine0} and~\ref{tab:others} evaluate the practical benefits of our branch-and-price methodology, as compared to a do-nothing baseline, practical benchmarks (i.e., easily-implementable heuristics) and optimization benchmarks (i.e., existing methods from the literature). To provide a fair assessment, we evaluate all solutions with the full continuous-state contagion models, as opposed to the state-clustering approximation used in our algorithm---thus ensuring an apples-to-apples comparison.

We define two practical benchmarks for the vaccine allocation and the congestion mitigation problems, as well as slight variants for the content promotion problem due to the sparsity constraint.
\begin{itemize}[topsep=0pt,itemsep=0pt]
    \item[--] \textbf{uniform allocation}: each segment receives the same amount of resources of at each epoch (in content promotion, we randomize the $K$ products at each epoch, repeated 100 times).
    \item[--] \textbf{cost-based allocation}: a domain-based benchmark in which each segment receives a constant share of resources proportionally to the total cost under the do-nothing baseline (in content promotion, we select the $K$ products for which promotions have the strongest impact).
\end{itemize}

For the vaccination centers problem, we consider a demographic-based ``\textbf{top-K}'' benchmark that selects the $K$ facilities that can serve the most people within access distance restrictions. We test it both with uniform vaccine allocations across the resulting $K$ facilities, and with optimized vaccine allocations (obtained by fixing the facility variables and optimizing subsequent resource allocation).

We also define three optimization benchmarks for the vaccine allocation and vaccination centers problems, from the recent prescriptive contagion analytics literature in epidemiology applications:
\begin{itemize}[topsep=0pt,itemsep=0pt]
    \item[--] \textbf{MIQO implementation}: mixed-integer bilinear optimization implementation of Problem $(\calP)$ in Gurobi 9.5 \citep{gurobiquadratic}, based on a time-discretization approximation of the ODEs.
    \item[--] \textbf{discretization} to approximate Problem $(\calP)$ via mixed-integer linear optimization, using time discretization to eliminate continuous-time dynamics and a staircase approximation of the infected population to handle bilinearities. This benchmark mirrors the approach from \cite{fu2021robust} in a robust optimization setting. To optimize its performance, we divide the $[0,2\%]$ interval into sub-intervals of length $\delta$ (since the infected population never exceeded 2\%). We consider a coarse and a granular discretization, with $\delta = 0.001$ and $\delta = 0.002$.
    \item[--] \textbf{coordinate descent heuristic}, which circumvents bilinearities by iterating between optimizing vaccine allocations over discretized time increments with fixed numbers of infections (via mixed-integer linear optimization) and re-estimating infections \citep{bertsimas2022locate}.
\end{itemize}

\begin{table}[h!]
\footnotesize\renewcommand{\arraystretch}{1.0}
\caption{Death toll comparison for the vaccine allocation problem (full country, $D=21$).}
\label{tab:vaccine0}
\begin{center}
\begin{tabular}{clcccccc}
\toprule
&&\multicolumn{2}{c}{\textbf{$S=6$}}& \multicolumn{2}{c}{\textbf{$S=8$}} & \multicolumn{2}{c}{\textbf{$S=12$}} \\ \cmidrule(lr){3-4}\cmidrule(lr){5-6}\cmidrule(lr){7-8}
Budget & \multicolumn{1}{c}{Method} & Time (sec.) & Deaths & Time (sec.) & Deaths & Time (sec.) & Deaths \\\toprule 
| & Do nothing & | & 573.37K& | & 604.22K & | & 649.75K \\\toprule
2.5M& Uniform allocation & \multicolumn{1}{c}{|} & -5.91K& \multicolumn{1}{c}{|} & -11.02K & \multicolumn{1}{c}{|} & -21.42K \\
& Cost-based allocation & \multicolumn{1}{c}{|} & -6.59K& \multicolumn{1}{c}{|} &  -12.92K & \multicolumn{1}{c}{|} & -27.64K \\\cmidrule{2-8}
& MIQO implementation & n/a & n/a& n/a & n/a & n/a & n/a \\
& Discretization ($\delta=0.002$) & 1,000$^*$ & -5.03K& n/a & n/a & n/a & n/a\\
& Discretization ($\delta=0.001$) & 1,000$^*$ & -5.56K& n/a & n/a & n/a & n/a\\
& Coordinate Descent & 5.43 & \textbf{-11.36K}& 4.73 & -20.02K & 6.03 & -39.01K \\\cmidrule{2-8}
&Branch-and-price ($\varepsilon=0.002$) & 426.2 & \textbf{-11.28K}& 1234.2 & \textbf{-20.43K} & 3671.4 & \textbf{-39.41K}  \\\toprule
7M&Uniform allocation& \multicolumn{1}{c}{|} & -13.22K& \multicolumn{1}{c}{|} & -23.67K & \multicolumn{1}{c}{|} & -43.60K\\
&Cost-based allocation& \multicolumn{1}{c}{|} & -17.44K& \multicolumn{1}{c}{|} & -31.77K  & \multicolumn{1}{c}{|}& -58.71K \\\cmidrule{2-8}
&MIQO implementation& n/a & n/a& n/a & n/a & n/a & n/a \\ 
&Discretization ($\delta=0.002$) & 1,000$^*$ & -5.88K& n/a & n/a & n/a & n/a \\ 
&Discretization ($\delta=0.001$) & 1,000$^*$ & -2.88K& n/a & n/a & n/a & n/a \\ 
&Coordinate Descent  & 5.46 & -20.61K& 12.7 & -36.22K & 118.5 & -62.55K \\\cmidrule{2-8}
&Branch-and-price ($\varepsilon=0.002$) & 67.5 & \textbf{-21.46K}& 156.8 & \textbf{-37.38K} & 708.7 & \textbf{-65.65K} \\\bottomrule
\end{tabular}
\end{center}
\begin{tablenotes}
    \vspace{-6pt}
    \item $*$ and ``n/a'': no optimal and feasible solution, respectively. Bold font:  solutions within $1\%$ of the best-found solution.
\end{tablenotes}
\end{table}

\begin{table}[h!]
\footnotesize\renewcommand{\arraystretch}{1.0}
\caption{Performance comparison: vaccination centers, congestion mitigation and content promotion problems.}
\label{tab:others}
\centering
\begin{adjustbox}{width=\textwidth}
\begin{tabular}{lcclcclcc}
\toprule
\multicolumn{3}{c}{Vaccination Centers} & \multicolumn{3}{c}{Content Promotion} & \multicolumn{3}{c}{Congestion Mitigation} \\ \cmidrule(lr){1-3}\cmidrule(lr){4-6}\cmidrule(lr){7-9}
Method & Time (sec.) & Deaths & Method & Time (sec.) & Market share & Method & Time (sec.) & Cost\\\hline 
 Do nothing & | & 573K &  Do nothing & | & 2.1 p.p.  &  Do nothing & | & 18.98  \\\hline
 Top-K, Uniform & \multicolumn{1}{c}{|} & -6.98K & Uniform (random) & \multicolumn{1}{c}{|} & \multicolumn{1}{c}{+5.1 p.p.} & Uniform & \multicolumn{1}{c}{|} & -1.80\% \\
  Top-K, Optimized & \multicolumn{1}{c}{|} & -7K & Cost-based & \multicolumn{1}{c}{|} & \multicolumn{1}{c}{+7.3 p.p.} 
  & Cost-based & \multicolumn{1}{c}{|} & -1.79\%\\\hline
 MIQO & n/a & n/a&  \multicolumn{1}{c}{|} & \multicolumn{1}{c}{|} & 
 \multicolumn{1}{c}{|}&  \multicolumn{1}{c}{|} & \multicolumn{1}{c}{|} & \multicolumn{1}{c}{|}\\
Discretization  & 10,275$^*$ & -4.3K & \multicolumn{1}{c}{|} & \multicolumn{1}{c}{|} & \multicolumn{1}{c}{|}  & \multicolumn{1}{c}{|} & \multicolumn{1}{c}{|} & \multicolumn{1}{c}{|}\\
Coordinate descent & 254 & -7.7K &\multicolumn{1}{c}{|} & \multicolumn{1}{c}{|} & \multicolumn{1}{c}{|} &\multicolumn{1}{c}{|} & \multicolumn{1}{c}{|} & \multicolumn{1}{c}{|} \\\hline
Branch-and-price & 10,978 & \textbf{-8.6K} & Branch-and-price & 11,864 & \textbf{+10.0 p.p.} & Branch-and-price & 416 & \textbf{-3.32\%}  \\\bottomrule
\end{tabular}
\end{adjustbox}
\begin{tablenotes}
    \vspace{-6pt}
    \item $*$ and ``n/a'': no optimal and feasible solution, respectively. Bold font:  solutions within $1\%$ of the best-found solution.\vspace{-3pt}
    \item Vaccination centers: $S=6$, $K=3$, results for the full country across 10 groups.\vspace{-3pt}
    \item Content promotion: $n=20$ $K=2$, $S=10$, $D=21$ (``hard'' instance); market share averaged across 20 products.\vspace{-3pt}
    \item Congestion mitigation: $n=5$; $B^1=6$; $B^2=4$; $S=6$ (``hard'' instance).
\end{tablenotes}
\vspace{-6pt}
\end{table}

Results show that our solution provides significant benefits against all benchmarks. Table~\ref{tab:vaccine0} reports results for the vaccine allocation problem, with weekly budgets of 2.5 million and 7 million vaccines---corresponding to the supplies originally planned and actually available in the United States in 2021, respectively. Without vaccinations, the pandemic would lead to an estimated 650,000 fatalities over three months (accounting for undetected deaths). Under a uniform allocation strategy, the vaccination campaign can save around 21,000 of these fatalities (or 3.3\%) with the smaller vaccine budget, and 43,000 fatalities (or 6.6\%) with the larger budget. The cost-based benchmark saves up to 30\% of additional lives, by capturing epidemiological information. But then, our solution saves 20-70\% additional lives over six weeks and 12-40\% additional lives over three months, as compared to the cost-based benchmark. In other words, optimized stockpile management can increase the effectiveness of the vaccination campaign by a factor of 1.12 to 1.7 without increasing vaccine capacity or vaccine efficacy. In absolute terms, these improvements represent 7,000 to 12,000 extra lives saved over three months, demonstrating the importance of vaccine distribution and the edge of optimization to guide resource allocation in complex contagion systems.

We obtain similar findings for the other three problems (Table~\ref{tab:others}). In the content promotion and congestion mitigation problems, the relative gains of the optimized solution are even more significant, estimated at 84-96\% against uniform allocation and at 37-85\% against cost-based allocation. Simple benchmarks perform comparatively worse due to the sparsity constraint in the content promotion problem and the interdependencies between prevention and treatment vehicles in the congestion mitigation problem---leading to stronger benefits of optimization. In the vaccination centers problem, the optimized solution also provides significant improvements against the simple demographic-based top-K benchmark. Specifically, the benefit of optimization amounts to 23-31\% under uniform vaccine allocations and to 3-23\% under optimized vaccine allocations. Exactly like epidemiological proxies are insufficient to guide tactical resource allocation, simple demographic proxies are insufficient to guide strategic planning in dynamical and non-linear contagion systems.

Finally, our branch-and-price methodology is instrumental to reap the benefits of optimization, with significant impact against simpler optimization benchmarks. First, Tables~\ref{tab:vaccine0} and~\ref{tab:others} show the strong limitations of simple workarounds to circumvent the non-convex contagion dynamics in the optimization problem. Direct implementation using mixed-integer quadratic optimization solvers fails to even provide a feasible solution in small six-week instances. Similarly, mixed-integer linear optimization implementation (using a discretization-based staircase approximation of the number of infections) performs worse than even uniform allocation in small instances and does not return a feasible solution in larger instances. In comparison, our algorithm achieves a Pareto improvement: faster computational times, stronger scalability, and much higher-quality solutions. Our methodology can also outperform the tailored coordinate descent heuristic, albeit in longer computational times. For the vaccine allocation problem, our solution yields a 1-5\% improvement, resulting in up to 3,000 extra lives saved over a three-month horizon. In fact, the benefits of our methodology increases as we increase the number of available vaccines, because the interaction between the vaccine allocation decisions and the infection dynamics become stronger, favoring our global optimization approach. In the vaccination centers problem, the benefits can also be significant when coordinate descent leads to local optima that induce a different set of facilities. In the $F=3$ instance, our methodology results in an extra 900 lives saved over six weeks---a 12\% improvement.

\begin{figure}[h!]
    \centering
    \includegraphics[scale=0.45]{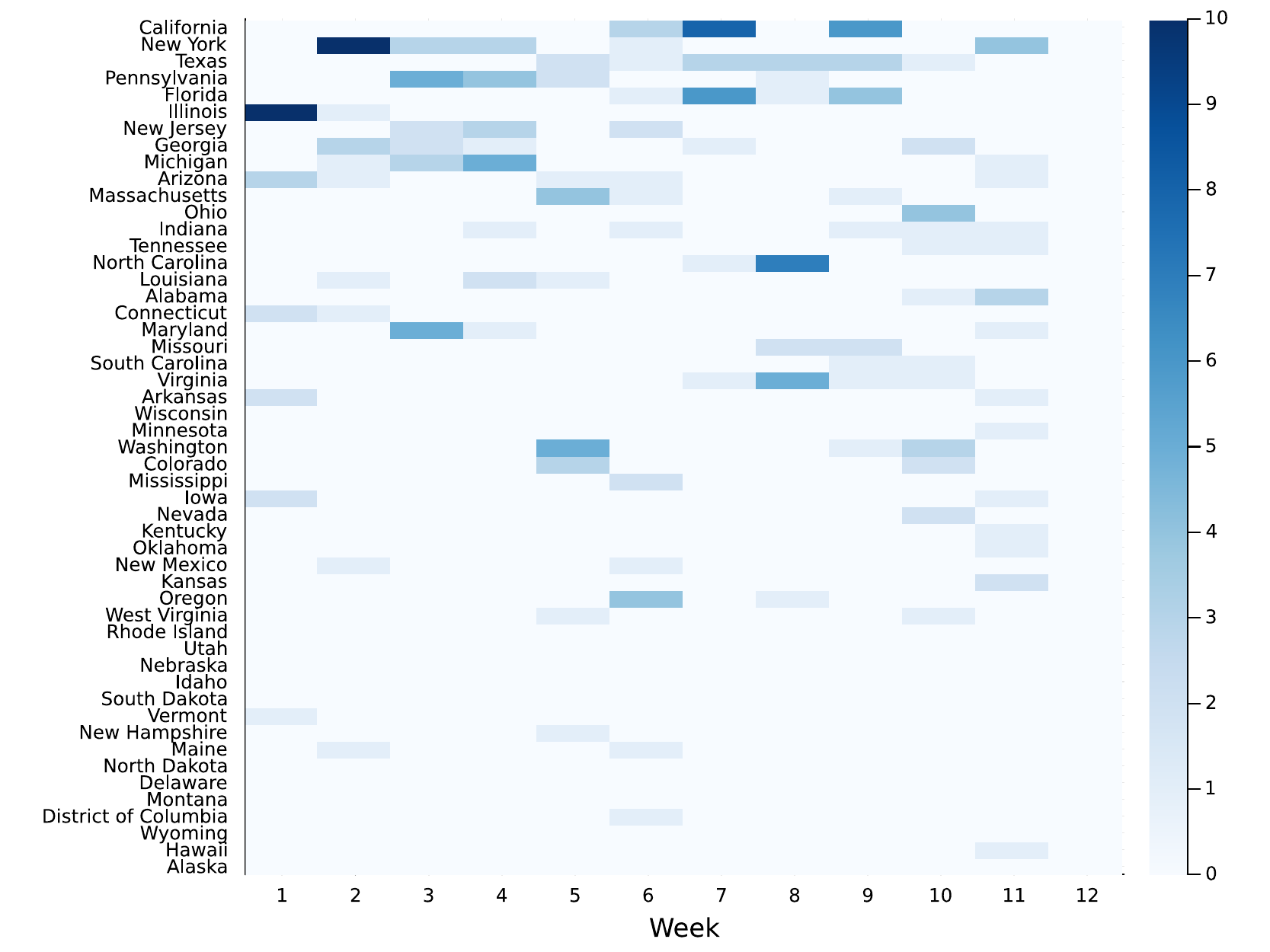}
    \caption{Spatial-temporal vaccine allocation over a 12 weeks. States are ordered from top to bottom in decreasing order of the total cost under do-nothing---a proxy for the prevalence of the pandemic.}
    \label{fig:vaccine_heatmap_2}
    \vspace{-24pt}
\end{figure}

\begin{figure}[h!]
\centering
\subfloat[Treatment vehicles, 8AM]{\includegraphics[width=0.33\textwidth]{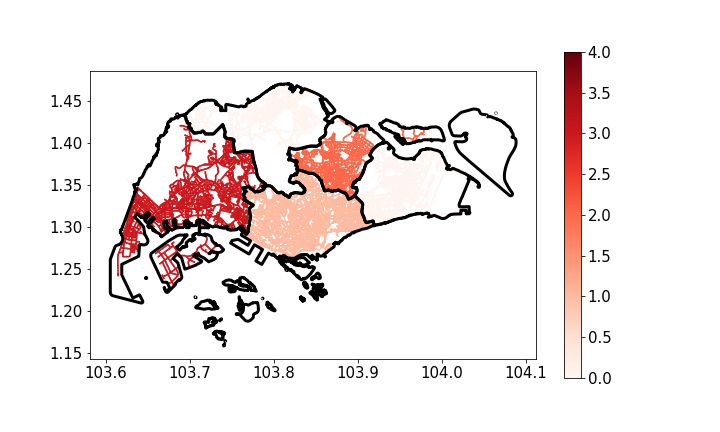}} 
\subfloat[Treatment vehicles, 9AM]{\includegraphics[width=0.33\textwidth]{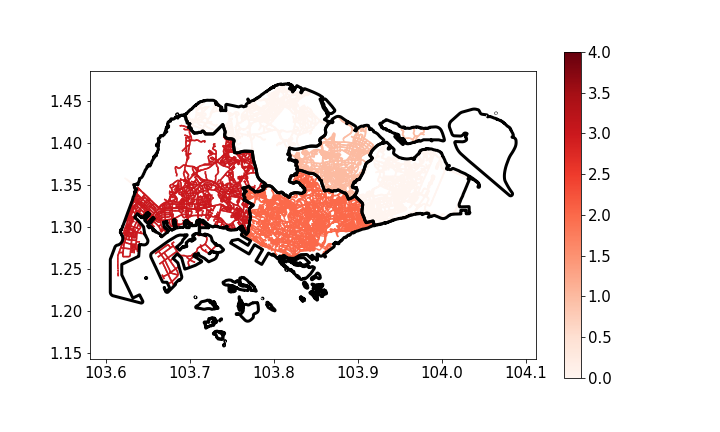}} 
\subfloat[Treatment vehicles, 10AM]{\includegraphics[width=0.33\textwidth]{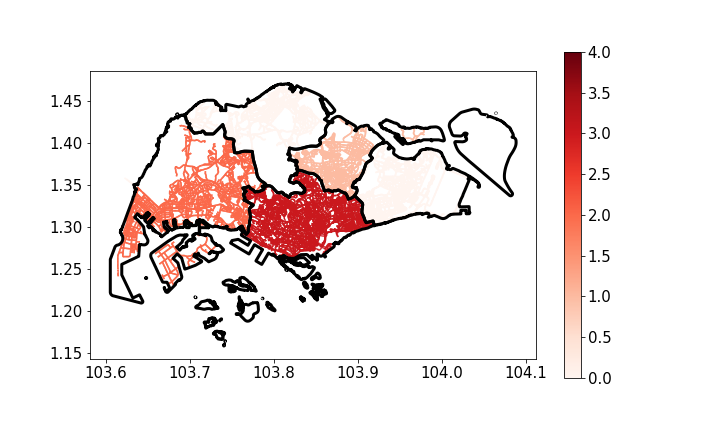}}
\hspace{0.1cm}
\subfloat[Prevention vehicles, 8AM]{\includegraphics[width=0.33\textwidth]{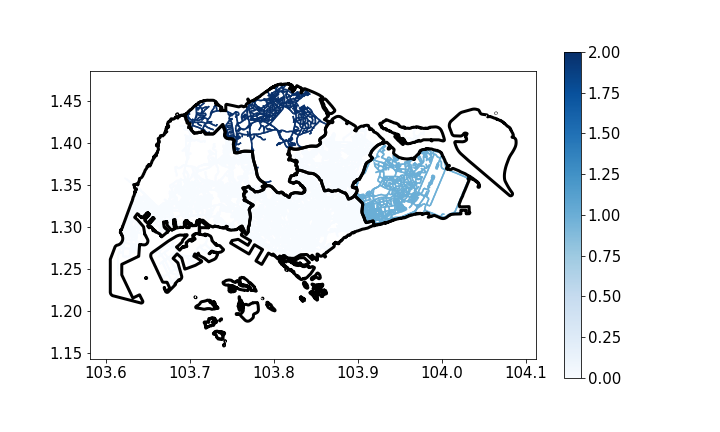}} 
\subfloat[Prevention vehicles, 8AM]{\includegraphics[width=0.33\textwidth]{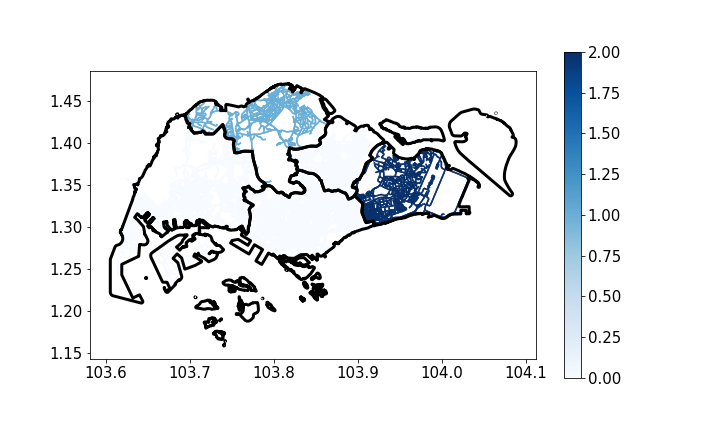}} 
\subfloat[Prevention vehicles, 8AM]{\includegraphics[width=0.33\textwidth]{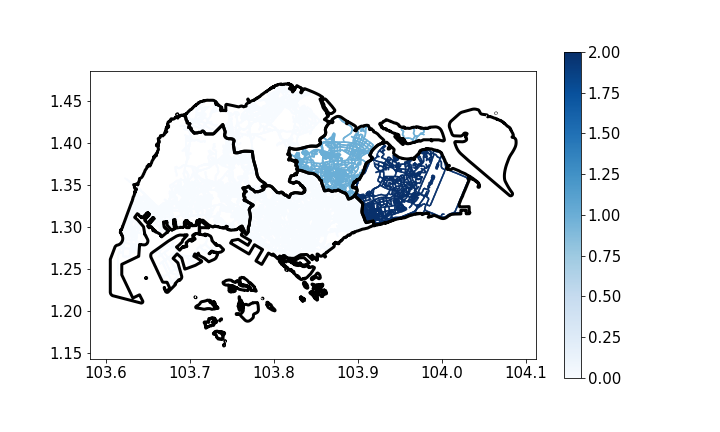}}
\caption{Allocation of treatment vehicles (top) and prevention vehicles (bottom) in Singapore over three hours.}
\label{fig:traffic_heatmap}
\vspace{-18pt}
\end{figure}

We conclude by illustrating vaccine allocations in Figure~\ref{fig:vaccine_heatmap_2} and vehicle allocations in Figure~\ref{fig:traffic_heatmap}. In vaccine allocation, the states where the pandemic is more costly (e.g., those with a larger population, more infections) generally receive more vaccines, but the relationship is not monotonic. Moreover, vaccine allocations exhibit different temporal patterns; for instance, New York and Illinois receive most vaccines early on; California and Texas receive most of their vaccines later on; and Florida stands in between. In the congestion mitigation problem, treatment vehicles are mainly sent to dense urban centers and near the airport, to clear large accidents and ease recovery; in contrast, prevention vehicles are primarily sent to suburban regions to prevent the spread of congestion. In both cases, there is no one-size-fits-all allocation strategy. Instead, the benefits of optimization stem from fine-tuning resource allocations to mitigate near-term impact via treatment interventions and to manage network propagation via prevention interventions.

\vspace{-6pt}

In summary, our branch-and-price methodology generates consistently high-quality solutions in large-scale practical instances of our prescriptive contagion analytics problems, with significant benefits as compared to all practical and optimization benchmarks. In Appendix~\ref{app:robust}, we report additional results showing the robustness of these findings to parameter estimation errors.
\section{Conclusion}
\label{sec:conclusion}

Predictive contagion systems have shown considerable success in epidemiology and other domains of science, engineering and management. This paper developed a prescriptive algorithm to optimize spatial-temporal resource allocation decisions in systems governed by contagion dynamics---and in more general dynamical systems. By combining the difficulties of combinatorial optimization and those of dynamical systems, however, this class of problems exhibits a large-scale and complex mixed-integer non-convex optimization structure with continuous-time ODE constraints. In response, this paper has developed a branch-and-price algorithm, using (i) a set partitioning reformulation to eliminate non-linearities; (ii) column generation to separate combinatorial optimization in a master problem from non-linear dynamics in a pricing problem; (iii) a novel state-clustering algorithm for discrete-decision continuous-state dynamic programming to solve the pricing problem; and (iv) a novel tri-partite branching scheme on natural variables to circumvent the non-linearities.

We implemented the methodology on four prescriptive analytics contagion problems: vaccine allocation, mass vaccination centers deployment, content promotion, and congestion mitigation. The branch-and-price algorithm scales to realistic and otherwise-intractable instances, with up to 21 decisions, 51 regions and 12 decision epochs in vaccine allocation, hence $\calO(21^{612})$ possible decisions overall. From a practical standpoint, the methodology significantly outperforms easily-implementable benchmarks based on demographic or epidemiological proxies, as well as state-of-the-art optimization algorithms. In the vaccine allocation example, the methodology can improve the effectiveness of a vaccination campaign by 12-70\%, resulting in 7,000 to 12,000 extra lives saves in a situation mirroring the midst of the COVID-19 pandemic in the United States. These results are robust across problem instances and robust to parameter estimation errors. Ultimately, our prescriptive contagion analytics methodology can deliver significant benefits across contagion-based domains, by fine-tuning resource allocations based on spatial-temporal system dynamics.

The promising results reported in this paper also motivate extensions of our methodology to tackle prescriptive contagion analytics problems with different decision-making structures (e.g., non-discrete resource allocation decisions, dynamic resource allocations) and more complex system dynamics (e.g., inter-population mixing). Perhaps the main limitation of this paper is the focus on deterministic resource allocation problems; although our robustness tests in Appendix~\ref{app:robust} showed the benefits of our solutions under model misspecification, our optimization methodology could be combined in future work with the robust epidemiological approach from \cite{fu2021robust}. Yet, this paper provides methodological foundations to support complex resource allocation decisions in a broad class of contagion-based systems, and of dynamical systems more generally.

\bibliographystyle{informs2014}
\bibliography{refs}

\newpage
\begin{APPENDICES}
\section{Details on the prescriptive contagion models}
\label{app:models}

In this appendix, we first formulate the four prescriptive contagion analytics problems as subcases of the general problem presented in Section~\ref{sec:model}. We then provide details on the contagion model for the congestion mitigation problem, which we developed and calibrated using real-world data from Singapore. Finally, we apply the set partitioning formulation proposed in Section~\ref{subsec:SP} to each problem.

\subsection{Formulation of prescriptive contagion models}
\label{app:formulation}

\subsubsection*{Vaccine allocation.}

A centralized decision-maker distributes $B_s$ vaccine doses each week. Each decision epoch $s=1,\cdots,S$ represents a week. Each segment represents one of 51 US regions (50 states plus Washington, DC). The dynamics in each region are characterized by a contagion model, shown in Figure~\ref{subfig:vaccines}. This model comprises seven states among the non-vaccinated population: susceptible (S), exposed (E), infected (I), undetected (U), hospitalizations (H), quarantine (Q), and death (D). In addition, the model includes four vaccinated states: susceptible (S’), exposed (E’), infected (I’), and immunized (M). The core model, called DELPHI, was fitted from historical data at the state level in the United States; it was incorporated into the ensemble forecast from the U.S. Center for Disease Control, with strong predictive performance across the various waves of the pandemic \citep{li2022forecasting}. It was then augmented into the DELPHI-V model to reflect the impacts of vaccinations by \cite{bertsimas2022locate}; we use this model in this paper, with the exception that we restrict ourselves to a single age group for tractability purposes.

The decision variable $x_{is}$ denotes the number of vaccines allocated to region $i=1,\cdots,n$ at epoch $s=1,\cdots,S$. Because vaccines are sent by pallets, this amount must be a multiple of the number of vaccines per pallet, denoted by $L$.  For convenience, let us denote by $\overline{x}_i(t)$ a piece-wise constant function such that $\overline{x}_i(t)=\frac{x_{is}}{\tau_{s+1}-\tau_s}$ for all $\tau_s\leq t\leq\tau_{s+1}$. Vaccinations induce a transfer of the susceptible population to the alternative susceptible state, proportionally to an effectiveness parameter $\beta$. This model assumes conservatively that vaccinations do not prevent infections but immunize patients, thus preventing downstream mortality.

The system dynamics are then governed by the following ordinary differential equations:
\begin{align}
\frac{dS_i(t)}{dt} &=-\alpha \gamma(t)\left(S_i(t)-\beta\overline{x}_i(t)\right)I_i(t)-\beta \overline{x}_i(t) \label{V:S}\\
\frac{dE_i(t)}{dt} &=\alpha \gamma(t)\left(S_i(t)-\beta\overline{x}_i(t)\right)I_i(t)-r^{I} E_i(t) \label{V:E}\\
\frac{dI_i(t)}{dt} &=r^{I} E_i(t)-r^{d} I_i(t) \label{V:I}\\
\frac{dU_i(t)}{dt} &=r_{k}^{U}(t) I_i(t)-r^{D} U_i(t) \label{V:U}\\
\frac{dH_i(t)}{dt} &=r_{k}^{H}(t) I_i(t)-r^{D} H_i(t) \label{V:H}\\
\frac{dQ_i(t)}{dt} &=r_{k}^{Q}(t) I_i(t)-r^{D} Q_i(t) \label{V:Q}\\
\frac{dD_i(t)}{dt} &=r^{D}\left(U_i(t)+H_i(t)+Q_i(t)\right)\label{V:D}\\
\frac{dM_{i}(t)}{dt}  &=\beta \overline{x}_i(t)\label{V:M}
\end{align}

The vaccine allocation problem is defined as follows. Equation~\eqref{eq:Vobj} minimizes the number of deaths (an absorbing state) along with the terminal number of hospitalized, quarantined and exposed people. Constraint~\eqref{eq:Vbudget} applies the budget of vaccines in each period. The remaining constraints define the contagion dynamics and apply restrictions on the domain of definition of the variables.
\begin{align}
\min \quad & \sum_{i=1}^n \left(c_DD_i(T) + c_{H}H_i(T) + c_{Q}Q_i(T) + c_{E}E_i(T)\right)\label{eq:Vobj}\\
\st \quad & \sum_{i=1}^n x_{is} \leq B_s \quad \forall s=1,\cdots,S\label{eq:Vbudget}\\
& \text{Equations~\eqref{V:S}--\eqref{V:M}}\\
& \text{Initial conditions}\\
& \bS, \bE, \bI, \bU, \bH, \bQ, \bD, \bM \geq \bo\\
& x_{is} \in \{0,L,2L,\cdots\},\ \forall i=1,\cdots,n,\ \forall s=1,\cdots,S\label{eq:Vdomain}
\end{align}

In fact, we maximize the number of lives saved by the vaccination campaign, rather than minimizing the death toll. Indeed, system dynamics may lead to ``unavoidable'' costs due to the initial conditions or to limited mitigation resources. To capture the effects of the vaccination campaign, we optimize the savings from a do-nothing baseline, with no vaccines. Obviously, the two formulations are equivalent, but maximizing the cost savings provides better estimates of the optimality gaps. Specifically, we denote by $\overline{D}_i$, $\overline{H}_i$, $\overline{Q}_i$ and $\overline{E}_i$ the states under the do-nothing baseline, and formulate the problem as follows:
\begin{align}
\max \quad & \sum_{i=1}^n \left(c_D\overline{D}_i(T) + c_{H}\overline{H}_i(T) + c_{Q}\overline{Q}_i(T) + c_{E}\overline{E}_i(T)-c_DD_i(T) - c_{H}H_i(T) - c_{Q}Q_i(T) - c_{E}E_i(T)\right)\label{eq:Vobj2}\\
\st \quad
& \text{Equations~\eqref{eq:Vbudget}--\eqref{eq:Vdomain}}
\end{align}

\subsubsection*{Vaccination centers.}

We consider the same system dynamics as in the vaccine allocation problem, but optimize the deployment of mass vaccination centers. Specifically, we select a subset of $K$ vaccination centers out of a set of $F$ candidates. At each epoch $s=1,\cdots,S$, a budget of $B_s$ vaccines needs to be allocated across all regions. Each facility $j=1,\cdots,F$ has access to a maximum capacity of $Cap_{js}$ vaccines at decision epoch $s=1,\cdots,S$---this capacity is not necessarily equal to $B_s/K$ to create flexibility in vaccine allocations. We impose a threshold on the maximum distance that a resident can travel to reach a vaccination center. Let $P_{ij}$ denote the number of residents in region $i$ than can be served by vaccination center $j$. For example, if vaccination center $j$ is located in New York City and $i$ refers to New Jersey, then $P_{ij}$ captures the number of New Jersey residents living within the distance threshold from New York City.

We define the following variables:
\begin{align*}
    y_{j}&=\begin{cases}1&\text{if vaccination center $j=1,\cdots,F$ is constructed}\\0&\text{otherwise}\end{cases}\\
    x_{ijs}&=\text{number of vaccines sent to region $i=1,\cdots,n$ from vaccination center $j=1,\cdots,F$ at epoch $s=1,\cdots,S$}
\end{align*}

The problem minimizes the number of deaths (Equation~\eqref{eq:VFobj}). Constraint~\eqref{eq:VFbudget_f} ensures that $K$ vaccination centers are selected. Constraint~\eqref{eq:VFBudget} applies the total budget of vaccines. Constraint~\eqref{eq:VFCap} applies the vaccine capacity in each center, and constraint~\eqref{eq:VFPopulation_Tot} ensures that vaccine allocation decisions are consistent with the coverage of each vaccination center. Finally, constraint~\eqref{eq:linking} couples facility location decisions and subsequent vaccine allocation decisions. This problem mirrors the problem studied by \cite{bertsimas2022locate} in the midst of the COVID-19 pandemic in the United States, again with the restriction to a single age group.
\begin{align}
\min \quad & \sum_{i=1}^n \left(c_DD_i(T) + c_{H}H_i(T) + c_{Q}Q_i(T) + c_{E}E_i(T)\right)\label{eq:VFobj}\\
\st \quad & \sum_{j=1}^F y_j = K \label{eq:VFbudget_f}\\
& \sum_{i=i}^n\sum_{j=1}^Fx_{ijs} \leq B_s \quad \forall s=1,\cdots,S \label{eq:VFBudget} \\
& \sum_{i=1}^nx_{ijs} \leq Cap_{js} \quad \forall s=1,\cdots,S, \quad \forall j=1,\cdots,F \label{eq:VFCap}\\
& \sum_{s=1}^Sx_{ijs} \leq P_{ij}\quad\forall i=1,\cdots,n,\quad\forall j=1,\cdots,F \label{eq:VFPopulation_Tot}\\
& x_{ijs} \leq \min(Cap_{js},P_{ij})y_{j}\quad\forall i=1,\cdots,n,\quad\forall j=1,\cdots,F,\quad\forall s=1,\cdots,S \label{eq:linking}\\
& \text{Equations~\eqref{V:S}--\eqref{V:M}}\\
& \text{Initial conditions}\\
& \bS, \bE, \bI, \bU, \bH, \bQ, \bD, \bM \geq \bo\\
& x_{ijs} \in \{0,L,2L,\cdots\},\ \forall i=1,\cdots,n,\ \forall j=1,\cdots,F,\ \forall s=1,\cdots,S\\
& \by \in \{0,1\}^F\label{eq:Vdomain_f}
\end{align}

Again, we actually maximize the number of lives saved, as follows:
\begin{align}
\max \quad & \sum_{i=1}^n \left(c_D\overline{D}_i(T) + c_{H}\overline{H}_i(T) + c_{Q}\overline{Q}_i(T) + c_{E}\overline{E}_i(T)-c_DD_i(T) - c_{H}H_i(T) - c_{Q}Q_i(T) - c_{E}E_i(T)\right)\label{eq:VFobj2}\\
\st \quad
& \text{Equations~\eqref{eq:VFbudget_f}--\eqref{eq:Vdomain_f}}
\end{align}

\subsubsection*{Content promotion.}

A platform needs to promote a portfolio of $n$ products. For implementation reasons, we assume that the platform can promote each product to a fraction of the population by increments of a constant fraction $L$ and can only promote up to $K$ products at a time. We define the following variables:
\begin{align*}
    x^1_{is}&=\begin{cases}1&\text{if product $i=1,\cdots,n$ is promoted at epoch $s=1,\cdots,S$}\\0&\text{otherwise}\end{cases}\\
    x^2_{is}&=\text{population share to which product $i=1,\cdots,n$ is promoted at epoch $s=1,\cdots,S$}
\end{align*}

Again, we define a piece-wise constant function $\overline{x}^2_{i}(t)=\frac{x^2_{is}}{\tau_{s+1}-\tau_s}$ for all $\tau_s\leq t\leq\tau_{s+1}$.

The contagion model captures product adoption dynamics via a Bass model, comprising susceptible users (A) and adopters (B). Product promotion increases adoption from external sources, which, in turn, increases positive externalities between adopters and susceptible users (Figure~\ref{subfig:content}). This model follows the one from \cite{lin2021content}, developed using data from a social media platform. The dynamics are defined as follows:
\begin{align}
\frac{dA_{i}(t)}{dt} &= -\alpha \overline{x}^2_{i}(t) A_i(t) - \frac{\beta}{m}S_i(t)B_i(t) \label{C:S}\\
\frac{dB_{i}(t)}{dt} &= \alpha \overline{x}^2_{i}(t) A_i(t) + \frac{\beta}{m}S_i(t)B_i(t)\label{C:B}
\end{align}

The problem maximizes total adoption (Equation~\eqref{eq:Cobj}), subject to a sparsity constraint ensuring that the platform promotes at most $K$ products per period (Constraint~\eqref{eq:Cbudget}), to a cover constraint ensuring that every user is shown a product (Constraint~\eqref{eq:Ccover}), to a linking constraint ensuring consistency between the $\bx^1$ and $\bx^2$ variables (Constraint~\eqref{eq:Cconsistency}), and to the contagion dynamics.
\begin{align}
\max \quad & \sum_{i=1}^n B_{i}(T)\label{eq:Cobj}\\
\st \quad & \sum_{i=1}^n x^1_{is} \leq K \quad \forall s=1,\cdots,S \label{eq:Cbudget}\\
& \sum_{i=1}^n x^2_{is} = 1 \quad \forall s=1,\cdots,S \label{eq:Ccover}\\
& x^2_{is} \leq x^1_{is}\quad \forall s=1,\cdots,S, \ \forall i=1,\cdots,n \label{eq:Cconsistency}\\
& \text{Equations~\eqref{C:S}--\eqref{C:B}}\\
& \text{Initial conditions}\\
& \bS, \bB  \geq \bo\\
& x^2_{is} \in \{0,L,2L,\cdots,1\},\ \forall i=1,\cdots,n,\ \forall s=1,\cdots,S\\
& \bx^1\in\{0,1\}^{n\times S}\label{eq:domainB}
\end{align}

Again, we maximize the improvements from a do-nothing baseline, with no product promotion. We denote by $\overline{B}_i$ the adoption state under the do-nothing baseline, and formulate the problem as follows:
\begin{align}
\max \quad & \sum_{i=1}^n (B_{i}(T)-\overline{B}_{i}(T))\label{eq:Cobj2}\\
\st \quad
& \text{Equations~\eqref{eq:Cbudget}--\eqref{eq:domainB}}
\end{align}

\subsubsection*{Congestion mitigation.}

We consider $n$ neighborhoods, and each population corresponds to the set of road segments in each neighborhood. The contagion model comprises six states, shown in Figure~\ref{subfig:traffic}: susceptible free-flow roads (S), road work (W), roads with accidents (A), roads with road work and accidents ($A'$), congested roads (I), and recovered roads (R). Details on this model are provided in Appendix~\ref{app:traffic}.

The centralized decision-maker can deploy $B^1$ emergency vehicles for treatment purposes (namely, to respond to major accidents in order to clear the roads faster) and $B^2$ emergency vehicles for prevention purposes (namely, to respond to minor accidents and prevent the formation of congestion). Accordingly, we define the following variables:
\begin{align*}
    x^1_{is}&=\text{number of treatment vehicles allocated to neighborhood $i=1,\cdots,n$ at epoch $s=1,\cdots,S$}\\
    x^2_{is}&=\text{number of prevention vehicles allocated to neighborhood $i=1,\cdots,n$ at epoch $s=1,\cdots,S$}
\end{align*}

Again, we define piece-wise constant functions $\overline{x}^1_{i}(t)=\frac{x^1_{is}}{\tau_{s+1}-\tau_s}$ and $\overline{x}^2_{i}(t)=\frac{x^2_{is}}{\tau_{s+1}-\tau_s}$ for all $\tau_s\leq t\leq\tau_{s+1}$.

The system in each neighborhood is governed by the following equations, where $\phi(\cdot,\cdot)$ and $\psi(\cdot,\cdot)$ denote generic functions reflecting the effects of the prevention and treatment interventions:
\begin{align}
&\frac{dS_i(t)}{dt} = -\alpha^FS_i(t)(I_i(t)+A_i(t)+A'_i(t)) - \beta^FS_i(t) + \rho^F(\zeta^F+\psi(\overline{x}^1_i(t),\zeta^F))A_i(t)\label{T:S}\\
&\frac{dA_i(t)}{dt} =  \beta^FS_i(t) - \rho^F(\zeta^F+\psi(\overline{x}^1_i(t),\zeta^F))A_i(t) - (1-\rho^F)(\zeta^F+\phi(\overline{x}^2_i(t),\zeta^F))A_i(t)\label{T:A}\\
&\frac{dW_i(t)}{dt} = -\alpha^WW_i(t)(I_i(t)+A_i(t)+A'_i(t)) - \beta^WW_i(t) + \rho^W(\zeta^W+\psi(\overline{x}^1_i(t),\zeta^W))A'_i(t)\label{T:W}\\
&\frac{dA'_i(t)}{dt} =  \beta^WW_i(t) - \rho^W(\zeta^W+\psi(\overline{x}^1_i(t),\zeta^W))A'_i(t) - (1-\rho^W)(\zeta^W+\phi(\overline{x}^2_i(t),\zeta^W))A'_i(t)\label{T:Aprime}\\
&\frac{dR_i(t)}{dt} = \theta I_i(t)\label{T:R}\\
&\frac{dI_i(t)}{dt} = -\left(\frac{dS_i(t)}{dt} + \frac{dW_i(t)}{dt} + \frac{dA_i(t)}{dt} + \frac{dA'_i(t)}{dt} + \frac{dR_i(t)}{dt}\right)\label{T:I}
\end{align}
In our experiments, we consider linear functions $\phi(\cdot,\cdot)$ and $\psi(\cdot,\cdot)$ defined as follows, although our model can accommodate more complex and non-linear functions.
\begin{align}
\phi(\overline{x}^2_i(t),\zeta) =  \frac{n\zeta}{2D_y}\left(\overline{x}^2_i(t)-\frac{B^2}{n}\right)\\
\psi(\overline{x}^1_i(t),\zeta) = \frac{n\zeta}{2D_x}\left(\overline{x}^1_i(t)-\frac{B^1}{n}\right),
\end{align}
The problem minimizes the costs of congestion and accidents, subject to budget constraints and congestion dynamics governed by the contagion model:
\begin{align}
\min \quad & \sum_{i=1}^n\int_0^T \left( c_I I_i(t) + c_{A}A_i(t) + c_{A'}A'_i(t) \right)dt\\
\st \quad & \sum_{i=1}^n x^1_{is} \leq B^1 \quad \forall s=1,\cdots,S \label{eq:congestionbudget}\\
& \sum_{i=1}^n x^2_{is} \leq B^2 \quad \forall s=1,\cdots,S\\
& \text{Equations~\eqref{T:S}--\eqref{T:I}}\\
& \text{Initial conditions}\\
& \bS, \bW, \bA, \bA', \bI, \bR \geq \bo\\
& \bx^1,\bx^2\in\Z_+^{n\times T}\label{eq:congestiondomain}
\end{align}

Again, we maximize the improvements from a do-nothing baseline, with no emergency vehicle. We denote by $\overline{I}_i$, $\overline{A}_i$ and $\overline{A}'_i$ the states under the do-nothing baseline, and formulate the problem as follows:
\begin{align}
\max \quad & \sum_{i=1}^n\int_0^T \left( c_I \overline{I}_i(t) + c_{A}\overline{A}_i(t) + c_{A'}\overline{A}'_i(t)-c_I I_i(t) - c_{A}A_i(t) - c_{A'}A'_i(t) \right)dt\\
\st \quad
& \text{Equations~\eqref{eq:congestionbudget}--\eqref{eq:congestiondomain}}
\end{align}

\subsection{A data-driven contagion model of traffic congestion}
\label{app:traffic}

We collected data in Singapore that tracks the average speed on each road every five minutes, as well as the log of incidents and road works in the city. We divide the city into five neighborhoods, and the population in each neighborhood characterizes all the road segments. In each five-minute interval, a road segment is labeled as congested if the average speed recorded is less than 30\% of the average speed on that road.

We considered three contagion models: a basic Susceptible-Infected-Recovered (SIR) model inspired from \cite{saberi2020simple}, an SAIR model with an ``accident'' state, and the full model depicted in Figure~\ref{subfig:traffic} where we separate road segments with and without road work (referred to as 2-SAIR). We divide each day into a morning period (7am to 1pm), an afternoon period (1pm to 5:30pm) and an evening period (5:30pm onward). We trained our contagion models on those three periods separately, using historical data from March 2021. We then tested the models on data from April 2021. The results are reported in Figure~\ref{fig:fitting} and Table~\ref{tab:rmse}.

\begin{figure}[h]
    \centering
    \subfloat[In sample (March 2021)]{\includegraphics[width=0.5\textwidth]{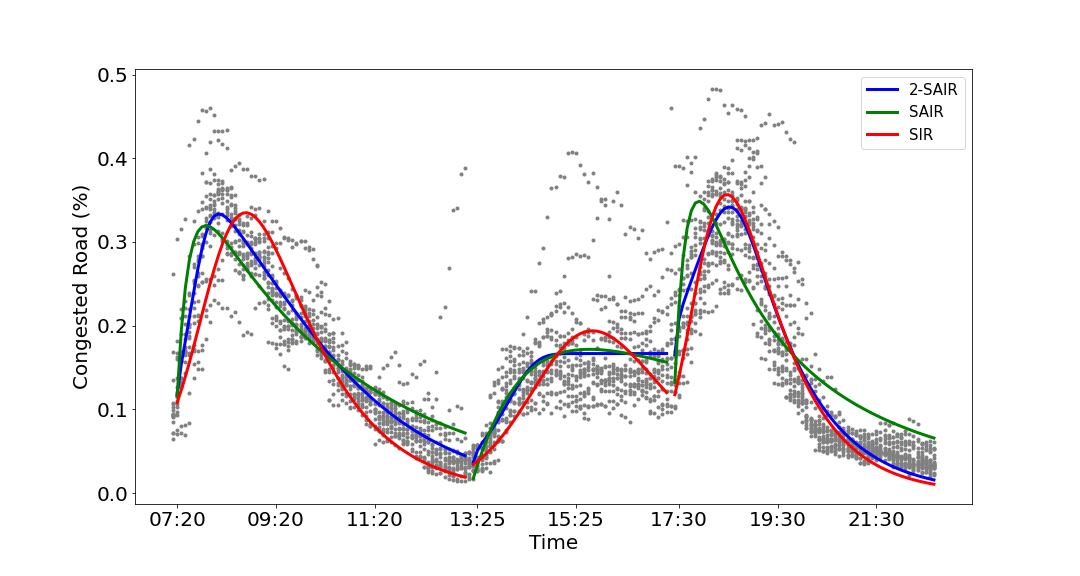}} 
    \subfloat[Out of sample (April 01, 2021)]{\includegraphics[width=0.5\textwidth]{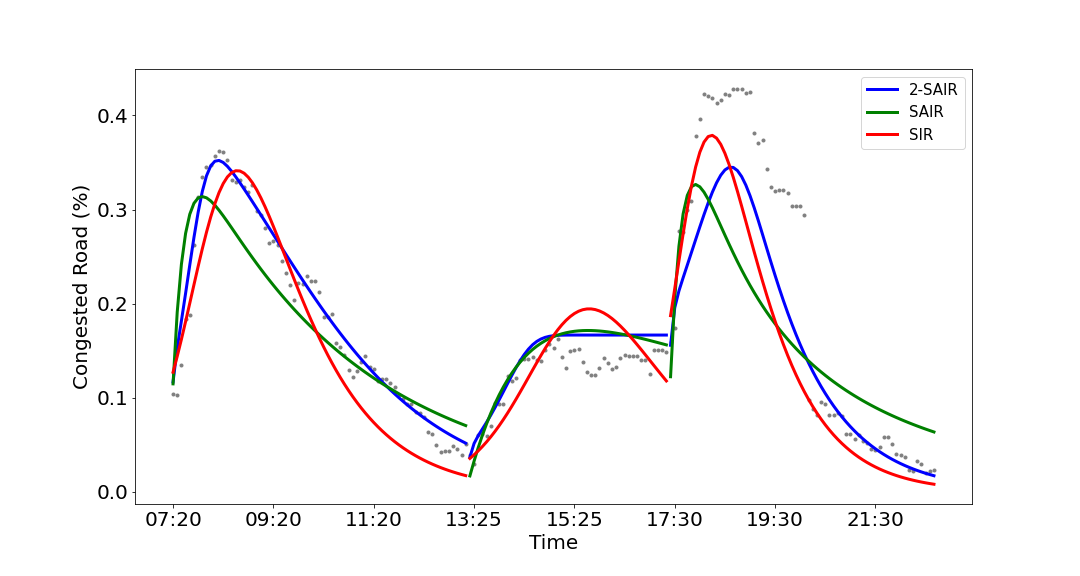}}
    \caption{Comparison of congestion levels predicted by the three contagion models.}
    \label{fig:fitting}
\end{figure}

\begin{table}[h!]
        \caption{Root-mean-square error (RMSE) for each model on the three periods of the day.}
    \label{tab:rmse}\small\centering\renewcommand{\arraystretch}{1.0}
    \begin{tabular}{lcccccc} \toprule
    &\multicolumn{3}{c}{In sample (March 2021}&\multicolumn{3}{c}{Out of sample (April 2021)}\\\cmidrule(lr){2-4}\cmidrule(lr){5-7}
    {Model} & {Morning} & {Afternoon} & {Evening} &{Morning} & {Afternoon} & {Evening} \\ \midrule
    {SIR}  & 0.0551 & 0.0546 & 0.0586  &0.0526 & 0.1028 & 0.2196  \\
    {SAIR}  & 0.0506  & 0.0502 & 0.0702   &0.0434  & 0.1002 & 0.205   \\
    {2-SAIR}  & \textbf{0.0389}  & \textbf{0.0501} & \textbf{0.0533}  &\textbf{0.0406}  & \textbf{0.0998} & \textbf{0.1885}  \\
     \bottomrule
    \end{tabular}
\end{table}

Note that the 2-SAIR model consistently outperforms the two benchmarks, both in sample and out of sample---with a 23\% reduction in RMSE during the morning peak and a 14\% reduction during the evening peak. Its parameters, reported in Table~\ref{tab:param}, provide some insights on the underlying dynamics; for instance, the accident rate and the congestion rate are both higher on roads segments with road work than on roads segments without road work. Ultimately, the model provides strong estimates of the congestion dynamics observed in the city, as illustrated in Figure~\ref{fig:fitting} across the training data and for one out-of-sample day.

\begin{table}[h!]
    \caption{Parameters of the 2-SAIR model for every region of the city.}
    \label{tab:param}
    \centering
    \small
    \begin{tabular}{llccccccccc} \toprule
    {Time period} & {Region} & {$\alpha^F$} & {$\beta^F$} & {$\rho^F$} & {$\zeta^F$} & {$\alpha^W$} & {$\beta^W$} & {$\rho^W$} & {$\zeta^W$} & {$\theta$} \\ \midrule
    Morning & West      &    0.531 &   0.578 &  1.126 &  0.569 &    3.966 &   1.585 &  0.411 &  4.835 &  10.354 \\
    & Central   &    0.000 &   0.048 &  0.034 &  0.076 &    1.224 &   0.875 &  1.180 &  1.544 &   2.135 \\
    & Northeast &    0.000 &   0.137 &  0.345 &  0.078 &    1.240 &   0.198 &  1.023 &  2.273 &   4.742 \\
    & East      &    0.052 &   0.041 &  0.432 &  0.086 &    2.745 &   3.086 &  1.201 &  1.934 &   3.950 \\
    & North     &    0.492 &   0.562 &  1.450 &  2.043 &    1.396 &   0.896 &  2.162 &  1.878 &   1.113 \\ \midrule
    Afternoon & West &    0.897 &   0.180 &  1.913 &  3.180 &    2.852 &   1.321 &  1.208 &  3.444 &   2.512 \\
    & Central   &    0.862 &   0.070 &  3.259 &  1.144 &    1.098 &   1.026 &  0.338 &  1.831 &   3.477 \\
    & Northeast &    0.782 &   0.024 &  5.786 &  0.283 &    1.215 &   2.360 &  0.002 &  3.555 &   1.381 \\
    & East      &    0.838 &   0.115 &  2.143 &  1.260 &    0.000 &   0.288 &  0.000 &  0.137 &   1.283 \\
    & North     &    1.538 &   0.196 &  1.856 &  1.656 &    0.987 &   1.085 &  0.996 &  0.930 &   1.111 \\ \midrule
    Evening & West      &    0.091 &   0.448 &  1.249 &  1.176 &    0.568 &   0.990 &  2.194 &  1.264 &   3.775 \\
    & Central   &    0.024 &   0.331 &  1.176 &  0.793 &    0.068 &   0.733 &  1.873 &  1.335 &   3.350 \\
    & Northeast &    0.018 &   0.304 &  1.247 &  1.326 &    0.116 &   1.514 &  1.728 &  1.123 &   2.861 \\
    & East      &    0.603 &   0.191 &  2.255 &  1.141 &    0.962 &   0.985 &  1.527 &  1.517 &   6.799 \\
    & North     &    0.001 &   0.580 &  0.905 &  0.698 &    1.148 &   1.110 &  1.136 &  1.142 &   4.264\\
     \bottomrule
    \end{tabular}

\end{table}

\subsection{Set partitioning formulation applied to prescriptive contagion models}
\label{app:SP}

\subsubsection*{Vaccine allocation.}

A plan $p$ is a list of $S$ elements $\alpha_{is}^p$, where $\alpha_{is}^p$ represents the number of vaccines allocated to region $i$ at epoch $s$ in plan $p$. The set partitioning formulation is given as follows:
\begin{align}
\min \quad & \sum_{i=1}^n\sum_{p\in\calP_i}C^p_iz^p_i  \\ 
\text{s.t} \quad & \sum_{i=1}^n \alpha_{is}^pz_i^p \leq B_s, \quad \forall s =1,\cdots,S\label{eq:Vlambda}\\ 
& \sum_{p \in \mathcal{P}_i}z_i^p = 1,\quad\forall i=1,\cdots,n\label{eq:Vmu} \\
& z^p_i \in \{0,1\} \quad \forall p \in \calP_i,\quad\forall i=1,\cdots,n
\end{align}

Let $\lambda_s\in\R_+$ and $\mu_i\in\R$ denote the dual variables associated with Constraints~\eqref{eq:Vlambda} and~\eqref{eq:Vmu}, respectively. We obtain the following pricing problem:
\begin{align}
\min \quad & \left(c_DD_i(T) + c_{H}H_i(T) + c_{Q}Q_i(T) + c_{E}E_i(T)\right) + \sum_{s=1}^S\lambda_s x_{is} - \mu_i\\
\st \quad & x_{is} \leq B_s \quad \forall s=1,\cdots,S\\
& \text{Equations~\eqref{V:S}--\eqref{V:M}}\\
& \text{Initial conditions}\\
& \bS, \bE, \bI, \bU, \bH, \bQ, \bD, \bM \geq \bo\\
& x_{is} \in \{0,L,2L,\cdots\},\ \forall i=1,\cdots,n,\ \forall s=1,\cdots,S
\end{align}

\subsubsection*{Vaccination centers.}

A plan $p$ is defined as above. The set partitioning formulation is given as follows:
\begin{align}
\min \quad & \sum_{i=1}^n\sum_{p\in\calP_i}C^p_iz^p_i  \\ 
\text{s.t} \quad & \sum_{j=1}^Fy_j = K \\
& \sum_{i=i}^n\sum_{j=1}^Fx_{ijs} \leq B_s \quad \forall s=1,\cdots,S \label{eq:SPbudget} \\
& \sum_{i=1}^nx_{ijs} \leq Cap_{js} \quad \forall s=1,\cdots,S, \quad \forall j=1,\cdots,F \label{eq:SPcapacity} \\
& \sum_{s=1}^Sx_{ijs} \leq P_{ij}\quad\forall i=1,\cdots,n,\quad\forall j=1,\cdots,F \label{eq:SPdistance} \\
& x_{ijs} \leq \min(Cap_{js},P_{ij})y_{j}\quad\forall i=1,\cdots,n,\quad\forall j=1,\cdots,F,\quad\forall s=1,\cdots,S \label{eq:SPlinking}\\
& \sum_{p \in \mathcal{P}_i} \alpha_{is}^pz_i^p = \sum_{j=1}^Kx_{i,j,s} \quad \forall s =1,\cdots,S \quad \forall i=1,\cdots,n \label{eq:VFconsistent}\\ 
& \sum_{p \in \mathcal{P}_i}z_i^p = 1\quad\forall i=1,\cdots,n \label{eq:VFunique}\\
& z^p_i \in \{0,1\} \quad \forall p \in \calP_i,\ \forall i=1,\cdots,n\\
& y_j\in\{0,1\}\quad\forall j=1,\cdots,F
\end{align}
 
Note that constraint~\eqref{eq:VFconsistent} ensures consistency between the number of vaccines sent to a region $i$ and the plan-based variables. Let $\lambda_{is}\in\R$ and $\mu_i\in\R$ denote the dual variables associated with Constraints~\eqref{eq:VFconsistent} and~\eqref{eq:VFunique}, respectively. We obtain the following pricing problem:
\begin{align}
\min \quad & \left(c_DD_i(T) + c_{H}H_i(T) + c_{Q}Q_i(T) + c_{E}E_i(T)\right) - \sum_{s=1}^S\lambda_{is} x_{is} - \mu_i\\
\st \quad & x_{is} \leq B_s \quad \forall s=1,\cdots,S\\
& \text{Equations~\eqref{V:S}--\eqref{V:M}}\\
& \text{Initial conditions}\\
& \bS, \bE, \bI, \bU, \bH, \bQ, \bD, \bM \geq \bo\\
& x_{is} \in \{0,L,2L,\cdots\},\ \forall i=1,\cdots,n,\ \forall s=1,\cdots,S
\end{align}

\subsubsection*{Content promotion.} 

A plan $p$ is a list of $S$ tuples $(\alpha^{1p}_{is},\alpha^{2p}_{is})$ where $\alpha^{1p}_{is}$ is equal to 1 if plan $p$ promotes product $i$ at epoch $s$ and $\alpha^{2p}_{is}$ is the share of the population to which product $i$ is shown at epoch $s$ under plan $p$. The set partitioning formulation is given by:
\begin{align}
\min \quad & \sum_{i=1}^n\sum_{p\in\calP_i}C^p_iz^p_i  \\ 
\text{s.t} \quad & \sum_{i=1}^n \alpha^{1p}_{is}z_i^p \leq K, \quad \forall s =1,\cdots,S\label{eq:Clambda1}\\ 
& \sum_{i=1}^n \alpha^{2p}_{is}z_i^p=1, \quad \forall s =1,\cdots,S\label{eq:Clambda2}\\ 
& \sum_{p \in \mathcal{P}_i}z_i^p = 1,\quad\forall i=1,\cdots,n\label{eq:Cmu} \\
& z^p_i \in \{0,1\} \quad \forall p \in \calP_i,\quad\forall i=1,\cdots,n
\end{align}

Let $\lambda^1_s\in\R_+$, $\lambda^2_s\in\R$ and $\mu_i\in\R$ denote the dual variables associated with Constraints~\eqref{eq:Clambda1},~\eqref{eq:Clambda2} and~\eqref{eq:Cmu}, respectively. We obtain the following pricing problem:
\begin{align}
\max \quad & B_{i}(T) + \sum_{s=1}^S\lambda^1_sx^1_{is} - \sum_{s=1}^S\lambda^2_sx^2_{is} - \mu_i\\
& x^2_{is} \leq x^1_{is}\quad \forall s=1,\cdots,S\\
& \text{Equations~\eqref{C:S}--\eqref{C:B}}\\
& x^2_{is} \in \{0,L,2L,\cdots,1\},\ \forall i=1,\cdots,n,\ \forall s=1,\cdots,S\\
& \bx^1\in\{0,1\}^{n\times S}
\end{align}

\subsubsection*{Congestion mitigation.}
A plan $p$ is a list of $S$ tuples $(\alpha^{1p}_{is},\alpha^{2p}_{is})$ where $\alpha^{1p}_{is}$ and $\alpha^{2p}_{is}$ are the number of emergency vehicles deployed for treatment and prevention purposes, respectively, in segment $i$ at epoch $s$ under plan $p$. The set partitioning formulation is given by:
\begin{align}
\min \quad & \sum_{i=1}^n\sum_{p\in\calP_i}C^p_iz^p_i  \\ 
\text{s.t} \quad & \sum_{i=1}^n \alpha^{1p}_{is}z_i^p \leq B^1, \quad \forall s =1,\cdots,S\label{eq:CMlambda1}\\
& \sum_{i=1}^n \alpha^{2p}_{is}z_i^p \leq B^2, \quad \forall s =1,\cdots,S\label{eq:CMlambda2}\\
& \sum_{p \in \mathcal{P}_i}z_i^p = 1,\quad\forall i=1,\cdots,n \label{eq:CMmu}\\
& z^p_i \in \{0,1\} \quad \forall p \in \calP_i,\quad\forall i=1,\cdots,n
\end{align}

Let $\lambda^1_s\in\R_+$, $\lambda^2_s\in\R_+$ and $\mu_i\in\R$ denote the dual variables associated with Constraints~\eqref{eq:CMlambda1},~\eqref{eq:CMlambda2} and~\eqref{eq:CMmu}, respectively. We obtain the following pricing problem:
\begin{align}
\min \quad & \int_0^T \left( c_I I_i(t) + c_{A}A_i(t) + c_{A'}A'_i(t) \right)dt + \sum_{s=1}^S\lambda^1_sx^1_{is}+\lambda^2_sx^2_{is} - \mu_i\\
\st \quad & \sum_{i=1}^nx^1_{is}  \leq B^1 \quad \forall s=1,\cdots,S\\
&  \sum_{i=1}^nx^2_{is} \leq B^2 \quad \forall s=1,\cdots,S\\
& \text{Equations~\eqref{T:S}--\eqref{T:I}}\\
& \bS, \bW, \bA, \bA', \bI, \bR \geq \bo\\ & \bx^1,\bx^2\in\Z_+^{n\times T} \
\end{align}

\section{Technical details}

\subsection{Proof of Proposition~\ref{prop:SP}}

Let $\widehat{\bz},\ \widehat{\by}$ be a feasible solution to $(\calS\calP)$. We define:
\begin{align*}
  \widehat{\bx}_{is} = \sum_{p\in\calP_i}\balpha_{is}^p z_i^p \quad \forall i =1,\dots,n, \ \forall s=1,\dots,S 
\end{align*}
By definition of the costs in $(\calS\calP)$ and in $(\calP)$, these two solutions achieve the same costs and satisfy constraints~\eqref{eq:ODE} and~\eqref{eq:initial}. Moreover, because of constraint (\ref{SP:one}) and by definition of $\balpha$, we know that, for all $i = 1,...,n$, there exists a single plan $p\in \calP_i$ such that $\widehat{z}_i^p=1$ and that $\balpha_{is}^p \in \calF_{is}$ for all $s=1,\dots,S$. Therefore, $\widehat{\bx}_{is} \in \calF_{is}$ for all $i=1,\dots,n,$ and $ s=1,\dots,S$. Finally, using constraint (\ref{SP:constraintsC}), we have: 
\begin{align*}
    \sum_{i=1}^n\bu^\top_{sji}\widehat{\bx}_{is} + \bv_{sj}^\top \widehat{\by} = \sum_{i=1}^n\sum_{p\in\calP_i}\bu_{sji}^\top \balpha_{is}^p\widehat{z}_i^p + \bv_{sj}^\top\widehat{\by} \geq w_{sj} \quad \forall s=1,\dots,S, \  \forall j=1,\dots,m_s
\end{align*}
Therefore, this yields a feasible solution to $(\calP)$ with the exact same cost.

Conversely, let $\widehat{\bx},\ \widehat{\by}$ be a solution to $(\calP)$. We define the plan $\widehat{p}_i$ such that:
$$\balpha_{is}^{\widehat{p}_i} = \widehat{\bx}_{is} \quad \forall i =1,\dots,n, \ \forall s=1,\dots,S$$
Now we consider the following plan-based solution, for all $i=1,\cdots,n$:
\begin{align*}
z^p_i = \begin{cases}
        1 & \text{if } p=\widehat{p}_i\\
        0 & \text{otherwise}
        \end{cases}
\end{align*}
By construction, constraints (\ref{SP:one}) and (\ref{SP:domain}) are satisfied. Moreover:
\begin{align*}
    \sum_{i=1}^n\sum_{p\in\calP_i}\bu_{sji}^\top\balpha_{is}^pz_i^p + \bv_{sj}^\top \widehat{\by} = \sum_{i=1}^n\bu_{sji}^\top\balpha_{is}^{\widehat{p}_i}z_i^{\widehat{p}_i}  + \bv_{sj}^\top \widehat{\by} = \sum_{i=1}^n\bu_{sji}^\top\widehat{\bx}_{is} + \bv_{sj}^\top \widehat{\by} \geq g_j^s \quad \forall s=1,\dots,S, \ \forall j=1,\dots,m_s
\end{align*}
We define the corresponding cost parameter $C^{\widehat{p}_i}_i$ per Equations~\eqref{SP:cost}--\eqref{SP:initial}. By construction, the two solutions achieve the same costs and satisfy constraints~\eqref{eq:ODE} and~\eqref{eq:initial}. This yields a feasible solution to $(\calS\calP)$ with the exact same cost. Therefore, Problem $(\calP)$ and Problem $(\calS\calP)$ are equivalent.\hfill\Halmos

\subsection{Pseudo-code of dynamic programming algorithms}\label{app:DP}

Algorithm~\ref{alg:dp_full} provides pseudo-code for the exact dynamic programming algorithm comprising the forward-enumeration and backward-induction loops. The forward-enumeration loop returns the full state space, denoted by $\calN^*_1,\cdots,\calN^*_{S+1}$ in the paper. The backward-induction loop computes the optimal policy given a state space. Algorithm~\ref{alg:dp_clust} provides pseudo-code for the state-clustering algorithm. It returns the clustered state space $\calN_1,\cdots,\calN_{S+1}$. The approximate dynamic programming algorithm considered in the paper consists of the state-clustering procedure in Algorithm~\ref{alg:dp_clust} and of the backward-induction loop in Algorithm~\ref{alg:dp_full}.

\begin{algorithm}
    \caption{Dynamic programming algorithm for pricing problem}
    \label{alg:dp_full}
    \renewcommand{\arraystretch}{1.0}\footnotesize
    \begin{algorithmic}[1] 
        \Procedure{StateSpaceCreation}{$\calN_1$}
            \For{$s \gets 1$ to $S$} \Comment{Forward enumeration}
            \State $\calN_{s+1} \gets \emptyset$
            \For{$\bN_s\in\calN_s,\ \bx_s\in\calF_s$}
                \State Transition: $\bM(\tau_s)\gets\bN_s$; $\frac{d\bM(t)}{dt}=f(\bM(t),\bx_{s}),\ \forall t\in\left[\tau_s,\tau_{s+1}\right]$; $\calN_{s+1} \gets \calN_{s+1} \cup \left\{\bM(\tau_{s+1})\right\}$
                \State Cost: $c_s(\bN_s,\bx_s) \gets \int_{\tau_s}^{\tau_{s+1}}g_{t}(\bM(t))dt+\Gamma_{s}(\bx_{s})$  
                \If{$s=S$}
                    \State $c_s(\bN_s,\bx_s)\gets c_s(\bN_s,\bx_s)+h(\bN_s)$
                \EndIf
            \EndFor
            \EndFor 
            \State \textbf{return} state space $\calN_1,\cdots,\calN_{S+1}$, cost $c_s(\bN_s,\bx_s)$ for all $\bN_s\in\calN_s,\ \bx_s\in\calF_s,\ s\in\{1,\cdots,S\}$
        \EndProcedure
        \Procedure{DynamicProgram}{$\calN_1,\cdots,\calN_S,\ c_1(\cdot,\cdot),\cdots,c_S(\cdot,\cdot),\lambda_{sj}$}   
            \State Initialization: $c^{opt}(\bN_s)=+\infty;\ \bpi^{opt}_{s}(\bN_s) \gets \emptyset$ $\forall \bN_s$ $\forall s \in \{1,\cdots, S\}$
            \For{$s \gets S$ to $1$} \Comment{Backward induction}
            \For{$\bN_s\in\calN_s$}
            \State $c^{temp} \gets +\infty$; $\bx^{temp} \gets \emptyset$    
            \For{$\bx_s\in\calF_s$}
                \State $\widetilde{c}_s(\bN_s, \bx_s) \gets c_s(\bN_s,\bx_s)-\sum_{j=1}^{m_s}\lambda_{sj}\bu^\top_{sj}\bx_{s}-\mathbf{1}\{s=S\}\mu$
                \State If $\widetilde{c}_s(\bN_s, \bx_s) < c^{temp}$, update $c^{temp}\gets \widetilde{c}_s(\bN_s, \bx_s) $; $\bx^{temp}\gets \bx_s$
            \EndFor
            \State $c^{opt}(\bN_s)=c^{temp};\ \bpi^{opt}_s(\bN_s) \gets\bx^{temp}$
            \EndFor
            \EndFor
            \State Initialization: $\bN^{opt}_1 \gets\bM^0$, $OPT\gets0$
            \For{$s \gets 1$ to $S$} \Comment{Forward evaluation}
                \State Find optimal decision $\bx^{opt}_s \gets \bpi^{opt}_s(\bN^{opt}_s)$ and update cost $OBJ \gets OBJ + c_s(\bN^{opt}_s, \bx^{opt}_s)$ 
                \State Find next state: $\bM(\tau_s)\gets\bN^{opt}_s$; $\frac{d\bM(t)}{dt}=f(\bM(t),\bx_{s}),\ \forall t\in\left[\tau_s,\tau_{s+1}\right]$; $\bN^{opt}_{s+1}\gets\bM(\tau_{s+1})$
            \EndFor
            \State \textbf{return} cost $OBJ$, decisions $\bx_s \; \forall s \in \{1, \cdots, S\}$
        \EndProcedure
    \end{algorithmic}
\end{algorithm}

\begin{algorithm}[htbp]
    \small
    \caption{State-clustering acceleration of dynamic programming algorithm for pricing problem}
    \label{alg:dp_clust}
    \renewcommand{\arraystretch}{1.0}\footnotesize
    \begin{algorithmic}[1]
        \Procedure{StateClustering}{$\calN_1$, $\varepsilon$}
            \For{$s \gets 1$ to $S$} \Comment{Forward pass}
            \State Initialization: $\calN_{s+1} \gets \emptyset$, $k \gets 0$ 
            \For{$\bN_s\in\calN_s,\ \bx_s\in\calF_s$}
                \State Next state: $\bM(\tau_s)\gets\bN_s$; $\frac{d\bM(t)}{dt}=f(\bM(t),\bx_{s}),\ \forall t\in\left[\tau_s,\tau_{s+1}\right]$; $\bN_{s+1}\gets\bM(\tau_{s+1})$
                \State  $l^*  \gets \argmin_{l\in \{0,1,\cdots, k\}} \max(\|\bN_{s+1}-\underline{\bN}(\gamma_l)\|_\infty, \|\bN_{s+1}-\overline{\bN}(\gamma_l)\|_\infty)$
                \If{$\max(\|\bN_{s+1}-\underline{\bN}(\gamma_{l^*})\|_\infty, \|\bN_{s+1}-\overline{\bN}(\gamma_{l^*})\|_\infty)\leq\varepsilon$}
                    \State Update: $\bN^\Sigma(\gamma_{l^*}) \gets \bN^\Sigma(\gamma_{l^*}) + \bN_{s+1}$ , $\eta(\gamma_{l^*}) \gets \eta(\gamma_{l^*}) + 1$
                    \State Update: $\underline{\bN}(\gamma_{l^*})\gets \min(\underline{\bN}(\gamma_{l^*}),\bN_{s+1})$, $\overline{\bN}(\gamma_{l^*})\gets \max(\overline{\bN}(\gamma_{l^*}),\bN_{s+1})$
                    \State Update: $\calX(\gamma_{l^*}) \gets \calX(\gamma_{l^*}) \cup \{(\bN_s, \bx_s)\}$, $c^\Sigma(\gamma_{l^*}) \gets c^\Sigma(\gamma_{l^*}) + \int_{\tau_s}^{\tau_{s+1}}g_{t}(\bM(t))dt+\mathbf{1}\{s=S\}(h(\bM(T)))$
                \Else
                    \State New cluster: $k \gets k + 1$, $\bN^\Sigma(\gamma_{k}) \gets \bN_{s+1}$, $\eta(\gamma_{k}) \gets 1$
                    \State New cluster: $\underline{\bN}(\gamma_{k})\gets \bN_{s+1}$, $\overline{\bN}(\gamma_{k})\gets \bN_{s+1}$
                    \State New cluster: $\calX(\gamma_k) \gets \{(\bN_s, \bx_s)\}$, $c^\Sigma(\gamma_{l^*}) \gets \int_{\tau_s}^{\tau_{s+1}}g_{t}(\bM(t))dt+\mathbf{1}\{s=S\}(h(\bM(T)))$
                \EndIf 
            \EndFor
            \For{$l \gets 1$ to $k$}
                \State $\calN_{s+1} \gets \calN_{s+1} \cup \left\{\bN^\Sigma(\gamma_{l}) / \eta(\gamma_{l})\right\}$\Comment{Cluster centroid}
                \For{$(\bN_s, \bx_s) \in \calX(\gamma_l) $}
                \State Define cost function $c_s(\bN_s,\bx_s) \gets c^\Sigma(\gamma_{l}) / \eta(\gamma_{l}) +\Gamma_{s}(\bx_{s})$
                \EndFor
            \EndFor
        \EndFor
        \EndProcedure
    \end{algorithmic}
\end{algorithm}

\subsection{Proof of Proposition~\ref{prop:cluster}}

Let $\gamma_l^{(q)}$ denote the $l^{\text{th}}$ cluster after $q$ states have been assigned to it. Let us denote by $\Phi_l$ the size of the $l^{\text{th}}$ cluster upon termination of the algorithm. Thus, the cluster at termination is denoted by $\gamma_l^{{\Phi_l}}$.

We proceed by contradiction. Assume that there exists a cluster $l$ such that:
 \[\|\overline{\bN}(\gamma_l^{(\Phi_l)})-\underline{\bN}(\gamma_l^{(\Phi_l)})\|_\infty>\varepsilon\]
Note that, by construction (Step 2), each cluster gets expanded at each iteration. Therefore, the quantity $\|\overline{\bN}(\gamma^{(i)}_l)-\underline{\bN}(\gamma^{(i)}_l)\|_\infty$ is non decreasing in $i$. Thus, there must exist $i\in\{1,\cdots,\Phi_l$\} such that:
 \begin{align*}
&\|\overline{\bN}(\gamma^{(i)}_l)-\underline{\bN}(\gamma^{(i)}_l)\|_\infty\leq \varepsilon \\
&\|\overline{\bN}(\gamma_l^{(i+1)})-\underline{\bN}(\gamma_l^{(i+1)})\|_\infty> \varepsilon 
 \end{align*}
By definition of the $\ell_\infty$ norm, we know that there exists a coordinate $j=1,\cdots,r$ such that:
 \begin{align*}
&|(\overline{\bN}(\gamma^{(i)}_l))_j-(\underline{\bN}(\gamma^{(i)}_l))_j|\leq \varepsilon\\
&|(\overline{\bN}(\gamma_l^{(i+1)}))_j-(\underline{\bN}(\gamma_l^{(i+1)}))_j|> \varepsilon. 
 \end{align*}
Denote the $i+1^{\text{th}}$ sample that got clustered into $\gamma_l$ as $\bN_{s+1}^*$. By construction (Step 2), at most one of $(\overline{\bN}(\gamma^{(i)}_l))_j$ and $(\underline{\bN}(\gamma^{(i)}_l))_j$ can change after $\bN_{s+1}^*$ gets added to the cluster. Without loss of generality, assume that $\bN_{s+1}^*$ changed the minimum state, that is:
\[(\bN_{s+1}^*)_j=(\underline{\bN}(\gamma_l^{(i+1)}))_j \neq (\underline{\bN}(\gamma^{(i)}_l))_j \qquad (\overline{\bN}(\gamma_l^{(i+1)}))_j = (\overline{\bN}(\gamma^{(i)}_l))_j \]
However, this implies:
\begin{align*}
\max(\|\bN_{s+1}^*-\underline{\bN}(\gamma^{(i)}_l)\|_\infty, \|\bN_{s+1}^*-\overline{\bN}(\gamma^{(i)}_l)\|_\infty)
&\geq \|\bN_{s+1}^*-\overline{\bN}(\gamma^{(i)}_l)\|_\infty\\
& \geq |(\bN_{s+1}^*)_j - (\overline{\bN}(\gamma^{(i)}_l))_j|\\
& = |(\underline{\bN}(\gamma_l^{(i+1)}))_j - (\overline{\bN}(\gamma_l^{(i+1)}))_j|\\
&>\varepsilon.
\end{align*}
This implies that $\bN_{s+1}^*$ would not be accepted in the cluster (Step 3), a contradiction. Therefore, we have:
 \[\|\overline{\bN}(\gamma_l^{(\Phi_l)})-\underline{\bN}(\gamma_l^{(\Phi_l)})\|_\infty\leq \varepsilon \;\; \forall l\]
 This completes the proof.\hfill\Halmos

\subsection{Proof of Proposition~\ref{prop:approx_error}}

By definition of Algorithm~\ref{alg:dp_clust}, the clustered and true states coincide at decision epoch $s=1$. By Proposition~\ref{prop:cluster}, the error at decision epoch $s=2$ is less than $\varepsilon$. We then proceed by induction, assuming that the result is true for some decision epoch $s$ and proving it for epoch $s+1$.

Define a true state $\bN^*_{s+1} \in \calN^*_{s+1}$. By definition of $\calN^*_{s+1}$, there exist a true state $\bN^*_{s} \in \calN^*_{s}$ and a decision $\bx_{s}\in\calF_s$ at the preceding epoch that leads to $\bN^*_{s+1}$, that is:
$$\bM^*(\tau_s)\gets\bN^*_s\ \text{and}\ \frac{d\bM^*(t)}{dt}=f(\bM^*(t),\bx_{s}),\ \forall t\in\left[\tau_s,\tau_{s+1}\right]\implies\bM^*(\tau_{s+1})=\bN^*_{s+1}.$$

Per the induction hypothesis, there exists a clustered state $\bN_{s}\in\calN_{s}$ such that:
\begin{equation}\label{eq:proof1}
\|\bN^*_s - \bN_s\|_\infty\leq\varepsilon\sum_{\sigma=3}^{s+1} \exp\left(\sum_{\nu=\sigma}^{s} L_{\nu-1}(\tau_{\nu} - \tau_{\nu-1})\right)
\end{equation}

In order to capture the propagation of approximation error between decision epochs $s$ and $s+1$, we make use of Gronwall's inequality. Let us specify it in Lemma~\ref{lem:gronwall}.
\begin{lemma}\label{lem:gronwall}
Let $L$ be a constant and $\bu:[a,b]\mapsto\md{R}^n$ be a real-valued $n$-dimensional continuous function defined over an interval $[a,b]$. Assume that $u$ is differentiable on the interior of $(a,b)$ and satisfies $\|\bu'(t)\|_\infty \leq L \|\bu(t)\|_\infty$ for all $t\in(a,b)$. Then, for all $t\leq b$, $\|\bu(t)\|_\infty \leq \|\bu(a)\|_\infty e^{L(t-a)}$.
\end{lemma}

Let us consider the state $\widehat\bN_{s+1}$ obtained at decision epoch $s+1$ starting from the clustered state $\bN_{s}\in\calN_{s}$ at epoch $s$ with decision $\bx_s\in\calF$:
$$\bM(\tau_s)\gets\bN_s\ \text{and}\ \frac{d\bM(t)}{dt}=f(\bM(t),\bx_{s}),\ \forall t\in\left[\tau_s,\tau_{s+1}\right]\ \text{and}\ \widehat\bN_{s+1}\gets\bM(\tau_{s+1}).$$

By assumption, the transition function is $L_s$-Lipschitz over $[\tau_s,\tau_{s+1}]$. Therefore, we have:
\begin{equation*}
\|f(\bM^*(t),\bx_s) - f(\bM(t), \bx_s)\|_\infty\leq L_s\|\bM^*(t)- \bM(t)\|_\infty,\quad  \forall t \in[\tau_s, \tau_{s+1}]
\end{equation*}

Per Lemma~\ref{lem:gronwall}, we obtain:
\begin{align}
    \|\bN^*_{s+1}-\widehat\bN_{s+1}\|_\infty&=\|\bM^*(\tau_{s+1})-\bM(\tau_{s+1})\|_\infty\nonumber\\
    &\leq e^{L_s(\tau_{s+1}-\tau_s)}\|\bM^*(\tau_{s})-\bM(\tau_{s})\|_\infty\nonumber\\
    &=e^{L_s(\tau_{s+1}-\tau_s)}\|\bN^*_s-\bN_s\|_\infty\label{eq:proof2}
\end{align}

Combining Equations~\eqref{eq:proof1} and~\eqref{eq:proof2}, we obtain:
\begin{align}
    \|\bN^*_{s+1}-\widehat\bN_{s+1}\|_\infty&\leq
    \varepsilon\sum_{\sigma=3}^{s+1} \exp\left(\sum_{\nu=\sigma}^{s} L_{\nu-1}(\tau_{\nu} - \tau_{\nu-1})\right)\cdot e^{L_s(\tau_{s+1}-\tau_s)}\nonumber\\
    &=\varepsilon\sum_{\sigma=3}^{s+1} \exp\left(\sum_{\nu=\sigma}^{s+1} L_{\nu-1}(\tau_{\nu} - \tau_{\nu-1})\right)\label{eq:proof3}
\end{align}

Now, by construction of the state-clustering algorithm (Algorithm~\ref{alg:dp_clust}), we assign state $\widehat\bN_{s+1}$ to a cluster with centroid $\bN_{s+1}\in\calN_{s+1}$ at decision epoch $s+1$. Per Proposition~\ref{prop:cluster}, we therefore know that:
\begin{equation}\label{eq:proof4}
    \|\bN_{s+1}-\widehat\bN_{s+1}\|_\infty\leq\varepsilon
\end{equation}

Combining Equations~\eqref{eq:proof3} and~\eqref{eq:proof4}, we obtain, per the triangular inequality:
\begin{align*}
    \|\bN^*_{s+1}-\bN_{s+1}\|_\infty&\leq\|\bN^*_{s+1}-\widehat\bN_{s+1}\|_\infty+\|\widehat\bN_{s+1}-\bN_{s+1}\|_\infty\\
    &\leq\varepsilon\left(1+\sum_{\sigma=3}^{s+1} \exp\left(\sum_{\nu=\sigma}^{s+1} L_{\nu-1}(\tau_{\nu} - \tau_{\nu-1})\right)\right)\\
    &=\sum_{\sigma=3}^{s+2} \exp\left(\sum_{\nu=\sigma}^{s+1} L_{\nu-1}(\tau_{\nu} - \tau_{\nu-1})\right),
\end{align*}
where we use the convention that an empty sum is equal to zero.

This completes the induction.\hfill\Halmos

\subsection{Proof of Theorem~\ref{thm:exact}}

Consider a node of the branching tree, and let $\widehat{\bz}$ denote the solution of the linear relaxation of the set partitioning formulation along with the corresponding branching constraints. We distinguish three cases:
\begin{enumerate}
    \item One of the variables $y_l$ is fractional for $l\in\{1,\cdots,q\}$. Then. the incumbent solution is eliminated by the bi-partite branching disjunction on $y_l$ (Step 2 of Algorithm~\ref{alg:exact}).
    \item The variable $\by$ satisfies all integrality constraints but one of the resource allocation constraints $\bx_{is}(\widehat\bz)\notin\calF_{is}$ for some segment $i\in\{1,\cdots,n\}$ and some decision epoch $s\in\{1,\cdots,S\}$. Then there must exist a component $k$ such that $\left\lfloor\bx_{is}(\widehat{\bz})\right\rfloor^k_{is}<\bx_{is}(\widehat{\bz})<\left\lceil\bx_{is}(\widehat{\bz})\right\rceil^k_{is}$. The incumbent solution is therefore eliminated by the bi-partite branching disjunction on $x^k_{is}(\bz)$ (Step 3 of Algorithm~\ref{alg:exact}).
    \item The variable $\by$ satisfies all integrality constraints, $\bx_{is}(\widehat\bz)\in\calF_{is}$ for all segments $i\in\{1,\cdots,n\}$ and all decision epochs $s\in\{1,\cdots,S\}$, but there exists a fractional plan-based variable $z^p_i\in(0,1)$. As noted in Section~\ref{subsec:BP}, there must exist a segment $i$, an epoch $s$, a component $k$ and a plan $p_0$ such that:
    \vspace{-6pt}
    \begin{equation*}
    \bx_{is}(\widehat{\bz})=\sum_{p\in\calP_i}\balpha^p_{is}\widehat{z}^p_i\in\calF_{is}\ \text{and}\ \alpha^{k,p_0}_{is}\neq x^k_{is}(\widehat{\bz})\ \text{and}\ \widehat{z}^{p_0}_i>0.
    \end{equation*}
    Then, the tri-partite branching disjunction (Equation~\eqref{eq:tripartite}, Step 4 of Algorithm~\ref{alg:exact}) makes the incumbent solution infeasible. To see this, let $\widehat\calP_i=\{p\in\calP_i:\widehat z^p_i>0\}$ denote the set of plans included in the incumbent solution for segment $i$, partitioned into $\widehat\calP^<_i=\{p\in\widehat\calP_i:\alpha^{pk}_{is}<x^k_{is}(\widehat\bz)\}$, $\widehat\calP^=_i=\{p\in\widehat\calP_i:\alpha^{pk}_{is}=x^k_{is}(\widehat\bz)\}$, and $\widehat\calP^>_i=\{p\in\widehat\calP_i:\alpha^{pk}_{is}> x^k_{is}(\widehat\bz)\}$. Since $\alpha^{k,p_0}_{is}\neq x^k_{is}(\widehat{\bz})$, we know that $\widehat\calP^<_i\neq\emptyset$ \textit{or} $\widehat\calP^>_i\neq\emptyset$. But since, in addition, $\sum_{p\in\calP_i}\alpha^{pk}_{is}\widehat{z}^p_i=x^k_{is}$, we actually know that $\widehat\calP^<_i\neq\emptyset$ \textit{and} $\widehat\calP^>_i\neq\emptyset$. This implies that the incumbent solution is infeasible in all three subsequent branches. Indeed:
    \begin{itemize}
        \item[--] In the left branch, we eliminate the plans in $\widehat\calP^>_i$ and in $\widehat\calP^=_i$, so $\alpha^{kp}_{is}<x^k_{is}(\widehat\bz)$ for all $p\in\calP$ such that $z_i^p>0$. Therefore, there exists no linear combination of remaining plan-based variables $\bz$ such that $x^k_{is}(\bz)=x^k_{is}(\widehat\bz)$.
        \item[--] In the right branch, the solution is infeasible for similar reasons.
        \item[--] In the middle branch, we eliminate the plans in $\widehat\calP^<_i$ and in $\widehat\calP^>_i$, so the branching constraints force $z_i^{p_0}=0$. Therefore, the incumbent solution is infeasible in the $\bz$ space.
    \end{itemize}
\end{enumerate}

This concludes that the branching scheme eliminates any infeasible incumbent solution. Note that the tri-partite branching proceeds differently in the left and right branches, on the one hand, and in the middle branch, on the other hand. Indeed, the left branch and the right branch force new resource allocation decisions. In contrast, the middle branch may result in the same resource allocation decisions but nonetheless restores integrality of $z_i^{p_0}$ variable and updates the cost estimates accordingly---thus contributing to tightening the relaxation bound and to progressing toward a feasible solution.

Finally, the branching scheme ensures that the solution is feasible upon termination. Recall that the algorithm terminates when $\sum_{p \in \calP_i}\widehat{z}^p_i| \alpha^{pk}_{is}-x^k_{is}(\widehat{\bz})|=0$ for all segments $i\in\{1,\cdots,n\}$, all decision epochs $s\in\{1,\cdots,S\}$, and all components $k$. Then, we know that $\alpha^{pk}_{is}=x^k_{is}(\widehat{\bz})$ whenever $\widehat{z}_i^p>0$. This ensures that a single plan $p\in\calP_i$ is such that $\widehat{z}_i^p>0$ for all segments $i\in\{1,\cdots,n\}$, and the solution is therefore integral. Finiteness of the algorithm follows from the discreteness of the solution space.
\newpage

\section{Additional computational results}

\subsection{Details on the experimental setup for the vaccination centers problem}\label{app:setup}

Recall that, for vaccine allocation problem, we apply our methodology to the full United States, comprising 51 regions (50 states plus Washington, DC), over a three-month period. The vaccination centers problem, however, is too complex to be solved at that scale. Instead, we solve it for smaller groups of states, defined by the CDC. Table~\ref{tab:groups} reports the composition of each of the 10 groups, along with their share of the US population, the budget of vaccines in each group, and the resulting number of discretized decisions (assuming pallets of 25,000 vaccines). Note that we consider here an overall budget of 2.5 million vaccines per week, which we distribution across the ten groups proportionally to the US population.

\begin{table}[h!]
\centering
\footnotesize
\caption{Group composition for the vaccination centers problem.}
\label{tab:groups}
\begin{tabular}{ccccc}
\toprule
Group & States & \% US Population & Vaccines & Decisions \\\hline 
& Connecticut, Maine, Massachusetts, & & &\\
A & New Hampshire, Rhode Island, Vermont, & 11.8\% & 300,000 & 12\\
& New York & & & \\\hline
 & Delaware, District of Columbia, Maryland, & & & \\
B & Pennsylvania, Virginia, West Virginia, & 9.8\% & 250,000 & 10\\
& New Jersey & & &\\\hline
C & North Carolina, South Carolina, & 12.7\% & 325,000   & 13\\
& Georgia, Florida  & & &\\\hline
D & Kentucky, Tennessee, & 4.9\% & 125,000 & 5\\
& Alabama, Mississippi & & &\\\hline
E & Illinois, Indiana, Michigan, & 15.6\% & 400,000 & 16\\
& Minnesota, Ohio, Wisconsin & & &\\\hline
F & Arkansas, Louisiana, New Mexico, & 12.7\% & 325,000 & 13\\
& Oklahoma, Texas & & & \\\hline
G & Iowa, Kansas, & 3.9\% & 100,000 & 4 \\
& Missouri, Nebraska & & & \\\hline
H & Colorado, Montana, North Dakota, & 3.9\% & 100,000 & 4\\
& South Dakota, Utah, Wyoming & & &\\\hline
I & Arizona, California,  & 19.6\% & 500,000 & 20\\
& Hawaii, Nevada & & & \\\hline
J & Alaska, Idaho, & 4.9\% & 125,000 & 5 \\
& Oregon, Washington & & & \\
\bottomrule
\end{tabular}
    \vspace{-12pt}
\end{table}

\subsection{Additional results on the branch-and-price algorithm}\label{app:results}

Tables~\ref{tab:BP_vaccine}--\ref{tab:BP_traffic} show extended results for the four problems, complementing the results shown in Table~\ref{tab:BP} of the paper. The observations are consistent with those from Section~\ref{subsec:Benefits}, thus showing the robustness of the findings across problem instances.

\begin{table}[h!]
\centering
\footnotesize\renewcommand{\arraystretch}{1.0}
\caption{Detailed performance assessment for the vaccine allocation problem.}
\label{tab:BP_vaccine}
\begin{tabular}{ccclcccccccc}
\toprule
& & & & & \multicolumn{4}{c}{CPU times (s)}  & \multicolumn{3}{c}{Solution quality} \\ \cmidrule(lr){6-9} \cmidrule(lr){10-12}
$n$ & $S$ & $D$ & \multicolumn{1}{c}{Method} & Nodes & Init. & RMP & PP & Total & Upper bound & Solution & Gap \\\hline
51 &   6 &   6 &  Column Generation &          1 &      15 &    0.02 &    0.14 &    0.2 &   45,995 &   45,563 &   0.94\% \\
 & & & Bi-partite B\&P &       14 &      15 &    0.11 &     0.69 &     5.25 &   45,821 &   45,780 &     0.09\% \\
 & & &  Tri-partite B\&P &       14 &      15 &    0.12 &     0.59 &     5.61 &   45,821 &   45,780 &     0.09\% \\\hline
  51 &   6 &  11 &  Column Generation &          1 &     27 &    0.02 &    0.55 &   0.61 &   48,240 &   48,108 &   0.27\% \\
 & & &  Bi-partite B\&P &       14 &     27 &    0.14 &     1.24 &     5.94 &   48,208 &   48,208 &     0\% \\
& & & Tri-partite B\&P &       14 &     27 &    0.17 &     1.5 &     7.01 &   48,208 &   48,208 &     0\% \\\hline
  51 &   6 &  21 &  Column Generation &          1 &     58 &    0.01 &    0.75 &    0.8 &   50,649 &   50,143 &   1\% \\
& & &  Bi-partite B\&P &     9,990 &     57 &   52 &    710 &  1,249 &   50,542 &   50,437 &     0.21\% \\
 & & &  Tri-partite B\&P &     9,990 &     57 &   51 &   709 &  1,250 &   50,542 &   50,437 &     0.21\% \\\hline
 51 &   8 &   6 &  Column Generation &          1 &     21 &    0.02 &     0.22 &   0.29 &   93,694 &   92,261 &   1.53\% \\
 & & & Bi-partite B\&P &       28 &     21 &     0.27 &     1.43 &     7.11 &   93,378 &   93,263 &     0.12\% \\
 & & &  Tri-partite B\&P &       46 &     21 &     0.42 &     2.22 &     9.62 &   93,353 &   93,263 &     0.09\% \\\hline
 51 &   8 &  11 &  Column Generation &          1 &     34 &     0.02 &    0.58 &   0.68 &   99,296 &   98,801 &   0.5\% \\
 & & & Bi-partite B\&P &      732 &     34 &    7.8 &    52 &   137 &   99,081 &   98,969 &     0.11\% \\
 & & &  Tri-partite B\&P &      744 &     34 &    7.4 &    53 &   139 &   99,068 &   98,969 &     0.1\% \\\hline
  51 &   8 &  21 &  Column Generation &          1 &     87 &    0.02 &    2.04 &   2.13 &  103,992 &  103,668 &   0.31\% \\
& & &  Bi-partite B\&P &     1,868 &     87 &   20 &   346 &   557 &  103,771 &  103,668 &     0.1\%\\
 & & &  Tri-partite B\&P &     1,868 &     87 &   20 &   346 &   555 &  103,771 &  103,668 &     0.1\% \\\hline
  51 &  10 &   6 &  Column Generation &          1 &     23 &    0.02 &    0.31 &   0.38 &  146,281 &  145,698 &   0.4\% \\
 & & &  Bi-partite B\&P &       16 &     23 &    0.17 &     1.39 &     6.22 &  145,947 &  145,810 &     0.09\% \\
 & & &  Tri-partite B\&P &       16 &     23 &    0.18 &     1.4 &     6.32 &  145,947 &  145,810 &     0.09\% \\\hline
 51 &  10 &  11 &  Column Generation &          1 &     43 &    0.02 &    0.81 &   0.92 &  156,464 &  155,624 &   0.54\% \\
 & & &  Bi-partite B\&P &     4,778 &     43 &   53 &   426 &  1,064 &  156,183 &  156,026 &     0.1\% \\
 & & &  Tri-partite B\&P &     4,778 &     43 &    55 &   420 &  1058 &  156,183 &  156,026 &     0.1\% \\\hline
 51 &  10 &  21 &  Column Generation &          1 &    118 &    0.02 &    2.55 &   2.78 &  164,927 &  163,111 &   1.1\% \\
 & & &  Bi-partite B\&P &     9,990 &    118 &   62 &   1,262 &  2,495 &  164,698 &  164,065 &     0.38\% \\
 & & &  Tri-partite B\&P &     9,990 &    118 &   62 &  1,251 &  2,486 &  164,698 &  164,065 &     0.38\% \\\hline
  51 &  12 &   6 &  Column Generation &          1 &     27 &     0.03 &    0.5 &    0.6 &  198,566 &  196,855 &   0.86\% \\
 & & &  Bi-partite B\&P &     1,224 &     27 &   17 &    63 &   238 &  197,468 &  197,272 &     0.01\% \\
 & & & Tri-partite B\&P &     1,302 &     27 &   19 &     68 &   257 &  197,468 &  197,272 &     0.01\% \\\hline
51 &  12 &  11 &  Column Generation &          1 &     51 &    0.03 &    0.95 &   1.04 &  212,653 &  211,412 &   0.58\% \\
 & & & Bi-partite B\&P &      122 &     51 &    1.37 &    18 &    35 &  212,027 &  211,978 &     0.02\% \\
 & & & Tri-partite B\&P &      122 &     51 &    1.33 &    18 &    35 &  212,027 &  211,978 &     0.02\% \\\hline
 51 &  12 &  21 &  Column Generation &          1 &    144 &     0.03 &    4.1 &   4.32 &  225,323 &  222,592 &   1.21\% \\
 & & &  Bi-partite B\&P &     9,990 &    144 &  101 &   1,765 &  3,475 &  225,048 &  224,495 &     0.25\% \\
 & & &  Tri-partite B\&P &     9,990 &    144 &  103 &  1,808 &  3,535 &  225,048 &  224,495 &     0.25\% \\
\bottomrule
\end{tabular}
\end{table}

\begin{table}[h!]
\centering
\footnotesize\renewcommand{\arraystretch}{1.0}
\caption{Detailed performance assessment for vaccination centers problem.}
\label{tab:BP_facility}
\begin{tabular}{cccccclcccccccc}
\toprule
 & & & & & & & & \multicolumn{4}{c}{CPU times (s)}  & \multicolumn{3}{c}{Solution quality} \\ \cmidrule(lr){9-12} \cmidrule(lr){13-15}
Group & $n$ & $S$ & $D$ & $K$ & Flex.? & \multicolumn{1}{c}{Method} & Nodes & Init. & RMP & PP & Total & Upper bound & Solution & Gap \\\hline
A &  7 &  6 &  13 &  2 &         No &  Column Generation &        1 &       58 &       0.12 &       0.33 &    0.96 &   7,224 &    6,109 &     18.66\% \\
& & & & & & Tri-partite B\&P &    36,017 &  58 &   794 &  440 &  10,813 &   7,094 &   6,665 &        6.05\% \\\hline
A &  7 &  6 &  13 &  2 &         Yes &  Column Generation &        1 &       45 &       0.19 &       0.43 &    1.11 &   8200 &    6660 &     29.91\% \\
& & & & & & Tri-partite B\&P &    11,899 &   45 &     224 &     178 &   10,808 &  8,132 &  7,917 &        2.64\% \\\hline
A &  7 &  6 &  13 &  3 &         No &  Column Generation &        1 &       46 &       0.18 &       0.27 &    0.86 &   8,200 &    6,703 &     39.32\% \\
& & & & & & Tri-partite B\&P &    16,586 &  45 &   383 &  241 &  10,839&   7,867 &   7,511 &        4.53\% \\\hline
A &  7 &  6 &  13 &  3 &         Yes &  Column Generation &        1 &       47 &       0.27 &       0.31 &    1.11 &   8,200 &    6,822 &     19.07\% \\
& & & & & &  Tri-partite B\&P &    16,358 &  47 &     363 &     209 &   10,955 &  8,048 &  7,716 &        4.13\% \\\hline
A &  7 &  8 &  13 &  2 &         No &  Column Generation &        1 &       79 &       0.30 &       1.34 &    2.17 &  16,562 &   12,471 &     48.54\% \\
& & & & & & Tri-partite B\&P &    24,453 &  78 &  1,138 &  701 &  10,827 &  16,381 &  15,842 &        3.29\% \\\hline
A &  7 &  8 &  13 &  2 &         Yes &  Column Generation &        1 &       54 &       0.33 &       1.01 &    2.02 &  18,662 &   16,417 &     17.69\% \\
& & & & & & Tri-partite B\&P &    33,034 &  53 &   979 &  754 &  10,803 &  18,628 &  18,443 &        0.99\% \\\hline
A &  7 &  8 &  13 &  3 &         No &  Column Generation &        1 &       54 &       0.41 &       1.02 &    1.87 &  18,652 &   14,619 &     27.66\% \\
& & & & & & Tri-partite B\&P &     8,638 &  53 &   289 &  305 &  10,891 &  17,740 &  16,929 &        4.57\% \\\hline
A &  7 &  8 &  13 &  3 &         Yes &  Column Generation &        1 &       53 &       0.45 &       1.05 &    1.93 &  18,662 &   14,918 &     31.59\% \\
& & & & & & Tri-partite B\&P &    12,287 &  53 &   343 &  375 &  10,939 &  18,274 &  17,187 &        5.94\% \\\hline
H &  6 &  6 &   5 &  2 &         No &  Column Generation &        1 &       49 &       0.03 &       0.04 &    0.55 &    786 &     652 &     16.97\% \\
& & & & & & Tri-partite B\&P &     4,678 &  49 &   147 &   24 &   1,026 &    764 &    763 &        0.09\% \\\hline
H &  6 &  6 &   5 &  2 &         Yes &  Column Generation &        1 &       36 &       0.05 &       0.03 &    0.46 &    786 &     652 &     29.15\% \\
& & & & & & Tri-partite B\&P &     2,919 &   36 &      75 &      11 &     604 &   764 &   763 &        0.09\% \\\hline
H &  6 &  6 &   5 &  3 &         No &  Column Generation &        1 &       42 &       0.03 &       0.02 &    0.48 &    678 &     652 &      3.75\% \\
& & & & & & Tri-partite B\&P &    45,144 &  42 &  1,278 &  113 &  10,802 &    675 &    674 &        0.23\% \\\hline
H &  6 &  6 &   5 &  3 &         Yes &  Column Generation &        1 &       43 &       0.02 &       0.03 &    0.57 &    678 &     652 &      3.75\% \\
& & & & & & Tri-partite B\&P &    58,496 &  43 &    1,382 &     137 &   10,801 &   675 &   674 &        0.23\% \\\hline
H &  6 &  8 &   5 &  2 &         No &  Column Generation &        1 &       45 &       0.04 &       0.08 &    0.68 &   1,978 &    1,671 &     15.55\% \\
& & & & & & Tri-partite B\&P &       62 &  44 &     1.43 &    0.89 &      6.99 &   1,978 &   1,978 &        0.00\% \\\hline
H &  6 &  8 &   5 &  2 &         Yes &  Column Generation &        1 &       42 &       0.08 &       0.07 &    0.58 &   1,978 &    1,753 &     49.08\% \\
& & & & & & Tri-partite B\&P &    20,280 &  41 &  2,002 &   76 &  10,818 &   1,886 &   1,865 &        1.11\% \\\hline
H &  6 &  8 &   5 &  3 &         No &  Column Generation &        1 &       43 &       0.05 &       0.05 &    0.53 &   1791 &    1780 &      0.60\% \\
& & & & & & Tri-partite B\&P &    14,149 &  42 &   298 &   71 &   1,835 &   1,787 &   1,787 &        0.05\% \\\hline
H &  6 &  8 &   5 &  3 &         Yes &  Column Generation &        1 &       45 &       0.03 &       0.05 &    0.54 &   1,791 &    1,780 &      0.60\% \\
& & & & & & Tri-partite B\&P &    11,927 &  45 &   231 &   63 &   1,506 &   1,787 &   1,787 &        0.05\%\\
\bottomrule
\end{tabular}
\begin{tablenotes}
    \vspace{-6pt}
    \item ``Flex.'' refers to flexibility in vaccine allocations across selected facilities. In the no-flexibility case, the capacity of each facility is equal to the vaccine budget (shown in Table~\ref{tab:groups}) divided by the number of facilities $K$. In the flexibility case, an extra pallet is available in each facility (Equation~\eqref{eq:SPcapacity}), with the same overall budget (Equation~\eqref{eq:SPbudget}).
\end{tablenotes}
\end{table}

\begin{table}[h!]
\centering
\footnotesize\renewcommand{\arraystretch}{1.0}
\caption{Detailed performance assessment for the content promotion  problem.}
\label{tab:BP_online}
\begin{tabular}{cccclcccccccc}
\toprule
  &  &  & & & & \multicolumn{4}{c}{CPU times (s)}  & \multicolumn{3}{c}{Solution quality} \\ \cmidrule(lr){7-10} \cmidrule(lr){11-13}
$n$ & $S$ & $D$ & $K$ & \multicolumn{1}{c}{Method} & Nodes & Init. & RMP & PP & Total & Upper bound & Solution & Gap \\\hline
20 &   6 &  11 &  4 &  Column Generation &         1 &       49 &       0.03 &       2.17 &    2.35 &  1.015 &   1.004 &   1.06\% \\
& & & & Bi-Partite B\&P &        8 &           49 &    0.08 &     4.04 &      5.9 &  1.015 &  1.015 &     0\% \\
& & & & Tri-partite B\&P &        8 &           49 &    0.09 &     4.11 &      6.1 &  1.015 &  1.015 &     0\% \\\hline
20 &   6 &  11 &  6 &  Column Generation &         1 &       56 &       0.02 &       1.76 &    2.01 &  1.017 &   1.007 &   0.96\% \\
& & & & Bi-Partite B\&P &    25,958 &           56 &  417 &  2,430 &   7,200 &  1.017 &  1.016 &     0.1\% \\
& & & & Tri-partite B\&P &    26,194 &           56 &  426 &  2,407 &   7,201 &  1.017 &  1.016 &     0.1\% \\\hline
20 &   6 &  21 &  2 &  Column Generation &         1 &       79 &       0.04 &       4.62 &    4.84 &  1.017 &   0.978 &   3.9\% \\
& & & & Bi-Partite B\&P &      696 &           79 &    9.35 &   198 &    327 &  1.015 &  1.009 &     0.53\% \\
& & & & Tri-partite B\&P &      142 &           79 &    1.67 &    52 &     74 &  1.013 &  1.013 &     0.09\% \\\hline
20 &   6 &  21 &  4 &  Column Generation &         1 &       72 &       0.03 &       3.77 &    4.06 &  1.052 &   1.04 &   1.13\% \\
& & & & Bi-Partite B\&P &    11,400 &           72 &  179 &  1,861 &   7,201 &  1.051 &  1.047 &     0.35\% \\
& & & & Tri-partite B\&P &     5,908 &           72 &  106 &  1,158 &   3,913 &  1.049 &  1.049 &     0.06\% \\\hline
20 &   6 &  21 &  6 &  Column Generation &         1 &       61 &       0.03 &       4.85 &    5.08 &  1.069 &   1.056 &   1.29\% \\
& & & & Bi-Partite B\&P &    19,660 &           61 &  223 &  2,761 &   7,202 &  1.067 &  1.065 &     0.17\% \\
& & & & Tri-partite B\&P &    19,640 &           61 &  222 &  2,766 &   7,201 &  1.067 &  1.065 &     0.17\% \\\hline
20 &   8 &  11 &  2 &  Column Generation &         1 &       58 &       0.04 &       5.01 &    5.37 &  1.527 &   1.44 &   5.7\% \\
& & & & Bi-Partite B\&P &    17,960 &           58 &  294 &  1,949 &   7,201 &  1.523 &  1.511 &     0.82\% \\
& & & & Tri-partite B\&P &    17,819 &           58 &  338 &  2,057 &   7,201 &  1.52 &  1.515 &     0.31\% \\\hline
20 &   8 &  11 &  4 &  Column Generation &         1 &       57 &       0.03 &       3.80 &    4.12 &  1.545 &   1.5 &   2.93\% \\
& & & & Bi-Partite B\&P &    18,876 &           57 &  158 &  1,845 &  10,801 &  1.545 &  1.533 &     0.74\% \\
& & & & Tri-partite B\&P &    19,016 &           57 &  156 &  1,783 &  10,801 &  1.545 &  1.53 &     0.74\% \\\hline
20 &   8 &  11 &  6 &  Column Generation &         1 &       79 &       0.05 &       5.64 &    6.05 &  1.551 &   1.529 &   1.41\% \\
& & & & Bi-Partite B\&P &    21,122 &           79 &  187 &  2,686 &  10,802 &  1.55 &  1.542 &     0.51\% \\
& & & & Tri-partite B\&P &    23,670 &           79 &  196 &  2,588 &  10,802 &  1.55 &  1.542 &     0.51\% \\\hline
20 &   8 &  21 &  4 &  Column Generation &         1 &      103 &       0.08 &      12 &   13 &  1.575 &   1.528 &   3.01\% \\
& & & & Bi-Partite B\&P &     6,546 &          103 &   77 &  1,642 &  10,802 &  1.572 &  1.564 &     0.45\% \\
& & & & Tri-partite B\&P &     6,536 &          103 &   74 &  1,685 &  10,805 &  1.572 &  1.564 &     0.45\% \\\hline
20 &   8 &  21 &  6 &  Column Generation &         1 &       79 &       0.05 &       8.61 &   10 &  1.602 &   1.587 &   0.91\% \\
& & & & Bi-Partite B\&P &     5,462 &           79 &   58 &  1,156 &  10,805 &  1.598 &  1.595 &     0.23\% \\
& & & & Tri-partite B\&P &     5,252 &           79 &   59 &  1,423 &  10,806 &  1.598 &  1.595 &     0.23\% \\\hline
20 &  10 &  11 &  2 &  Column Generation &         1 &       78 &       0.08 &       8.64 &    9.08 &  2.078 &   1.977 &   4.82\% \\
& & & & Bi-Partite B\&P &    16,446 &           78 &  369 &  3,337 &  10,803 &  2.07 &  2.045 &     1.22\% \\
& & & & Tri-partite B\&P &    15,688 &           78 &  347 &  3,627 &  10,803 &  2.07 &  2.045 &     1.22\% \\\hline
20 &  10 &  11 &  4 &  Column Generation &         1 &       82 &       0.06 &       6.85 &    7.24 &  2.098 &   2.070 &   1.31\% \\
& & & & Bi-Partite B\&P &    11,450 &           82 &  131 &  1,521 &  10,802 &  2.098 &  2.09 &     0.35\% \\
& & & & Tri-partite B\&P &    11,304 &           82 &  130 &  1,570 &  10,802 &  2.097 &  2.09 &     0.35\% \\\hline
20 &  10 &  11 &  6 &  Column Generation &         1 &       90 &       0.08 &       7.87 &    8.65 &  2.105 &   2.088 &   0.81\% \\
& & & & Bi-Partite B\&P &    10,394 &           90 &  122 &  1,435 &  10,802 &  2.105 &  2.097 &     0.41\% \\
& & & & Tri-partite B\&P &    10,536 &           90 &  120 &  1,384 &  10,802 &  2.105 &  2.097 &     0.41\% \\\hline
20 &  10 &  21 &  2 &  Column Generation &         1 &      134 &       0.08 &      15 &   16 &  2.017 &   1.874 &   7.08\% \\
& & & & Bi-Partite B\&P &    12,606 &          134 &  229 &  3,221 &  10,801 &  2.009 &  1.988 &     1.04\% \\
& & & & Tri-partite B\&P &    11,864 &          134 &  229 &  3,279 &  10,801 &  2.009 &  1.988 &     1.04\% \\\hline
20 &  10 &  21 &  4 &  Column Generation &         1 &      111 &       0.09 &      19 &   32 &  2.115 &   2.058 &   2.7\% \\
& & & & Bi-Partite B\&P &     1,050 &          111 &   14 &   454 &  10,823 &  2.112 &  2.099 &     0.61\% \\
& & & & Tri-partite B\&P &     1,042 &          111 &   14 &   527 &  10,827 &  2.112 &  2.099 &     0.61\% \\\hline
20 &  10 &  21 &  6 &  Column Generation &         1 &      124 &       0.08 &      15 &   18 &  2.164 &   2.147 &   0.69\% \\
& & & & Bi-Partite B\&P &     1,016 &          124 &   13 &   331 &  10,818 &  2.162 &  2.151 &     0.5\% \\
& & & & Tri-partite B\&P &     1,016 &          124 &   14 &   347 &  10,821 &  2.162 &  2.151 &     0.5\%\\ 
\bottomrule
\end{tabular}
\end{table}

\begin{table}[h!]
\centering
\footnotesize\renewcommand{\arraystretch}{1.0}
\caption{Detailed performance assessment for the congestion mitigation problem.}
\label{tab:BP_traffic}
\begin{tabular}{cccclcccccccc}
\toprule
  &  &  & & & & \multicolumn{4}{c}{CPU times (s)}  & \multicolumn{3}{c}{Solution quality} \\ \cmidrule(lr){7-10} \cmidrule(lr){11-13}
$n$ & $S$ & $B^1$ & $B^2$ & \multicolumn{1}{c}{Method} & Nodes & Init. & RMP & PP & Total & Upper bound & Solution & Gap \\\hline
5 &  2 &  4 &  2 &  Column Generation &          1 &     7.54 &     0.05 &    0.02 &    0.1 &   0.50 &     0.50 &   0\% \\
& & & & Bi-partite B\&P &        1 &     7.54 &    0.01 &    0.02 &    0.41 &  0.50 &  0.50 &     0\% \\
& & & & Tri-partite B\&P &        1 &     7.54 &    0.01 &    0.02 &      0.41 &  0.50 &  0.50 &     0\% \\\hline
5 &  2 &  4 &  4 &  Column Generation &          1 &     8.37 &    0.01 &    0.04 &   0.05 &   0.63 &     0.63 &   0\% \\
& & & & Bi-partite B\&P &        1 &     8.37 &   0.03 &    0.04 &    0.48 &  0.63 &  0.63 &     0\% \\
& & & & Tri-partite B\&P &        1 &     8.37 &    0.01 &    0.04 &     0.45 &  0.63 &  0.63 &     0\% \\\hline
5 &  2 &  6 &  2 &  Column Generation &          1 &     8.04 &     0.01 &    0.03 &   0.06 &   0.68 &     0.67 &   0.6\% \\
& & & & Bi-partite B\&P &        2 &     8.04 &   0.04 &    0.05 &    0.54 &  0.68 &  0.68 &     0.2\% \\
& & & & Tri-partite B\&P &        5 &     8.04 &   0.02 &     0.07 &      0.8 &  0.68 &  0.68 &     0.02\% \\\hline
5 &  2 &  6 &  4 &  Column Generation &          1 &    10 &    0.01 &    0.09 &   0.11 &   0.82 &     0.81 &   1.23\% \\
& & & & Bi-partite B\&P &        2 &    10 &   0.02 &    0.13 &    0.56 &  0.82 &  0.81 &     0.16\% \\
& & & & Tri-partite B\&P &        8 &    10s &   0.05 &     0.23 &      0.8s &  0.82 &  0.82 &     0\% \\\hline
5 &  2 &  8 &  2 &  Column Generation &          1 &     8.85 &    0.01 &    0.26 &   0.28 &   0.87 &     0.87 &   0\% \\
& & & & Bi-partite B\&P &        1 &     8.85 &   0.01 &    0.06 &    0.48 &  0.87 &  0.87 &     0\% \\
& & & & Tri-partite B\&P &        1 &     8.85 &    0.01 &    0.06 &     0.51 &  0.87 &  0.87 &     0\% \\\hline
5 &  2 &  8 &  4 &  Column Generation &          1 &    12 &    0.01 &    0.17 &   0.19 &   1.01 &     1.00 &   0.72\% \\
& & & & Bi-partite B\&P &        8 &    12 &   0.06 &    0.36 &    0.95 &  1.01 &  1.01 &     0.27\% \\
& & & & Tri-partite B\&P &       14 &    12 &   0.08 &     0.74 &     1.49 &  1.01 &  1.01 &     0.05\% \\\bottomrule
5 &  4 &  4 &  2 &  Column Generation &          1 &     14 &     0.02 &    0.41 &   0.44 &   1.34 &     1.34 &   0.31\% \\
& & & & Bi-partite B\&P &        1 &     14 &   0.02 &    0.4 &    0.83 &  1.34 &  1.34 &     0.31\% \\
& & & & Tri-partite B\&P &        6 &     14 &   0.05 &    0.92 &     1.74 &  1.34 &  1.34 &    0\% \\\hline
5 &  4 &  4 &  4 &  Column Generation &          1 &    49 &    0.03 &    2.3 &   2.35 &   1.61 &     1.55 &   4.16\% \\
& & & & Bi-partite B\&P &       28 &    49 &   0.2 &    8.7 &    10.8 &  1.61 &  1.59 &     1.06\% \\
& & & & Tri-partite B\&P &       73 &    49 &   0.79 &   19.4 &    25 &  1.60 &  1.60 &     0.02\% \\\hline
5 &  4 &  6 &  2 &  Column Generation &          1 &    33 &    0.02 &     1.35 &   1.39 &   1.90 &     1.90 &   0.19\% \\
& & & & Bi-partite B\&P &        4 &    33 &   0.04 &    2.12 &    2.73 &  1.90 &  1.90 &     0.11\% \\
& & & & Tri-partite B\&P &        7 &    33 &   0.08 &    3.38 &     4.07 &  1.90 &  1.90 &     0.04\% \\\hline
5 &  4 &  6 &  4 &  Column Generation &          1 &   185 &    0.02 &    6.21 &   6.27 &   2.20 &     2.05 &   6.61\% \\
& & & & Bi-partite B\&P &       48 &   185 &   0.34 &   32 &   34 &  2.19 &  2.19 &     0.21\% \\
& & & & Tri-partite B\&P &       57 &   185 &   0.43 &   39 &    42 &  2.19 &  2.19 &     0.06\% \\\hline
5 &  4 &  8 &  2 &  Column Generation &          1 &    83 &    0.03 &    4.02 &   4.1 &   2.47 &     2.44 &   1.11\% \\
& & & & Bi-partite B\&P &        2 &    83 &   0.05 &    5.4 &    5.92 &  2.47 &  2.44 &     1.09\% \\
& & & & Tri-partite B\&P &       63 &    83 &   0.63 &   29 &    33 &  2.46 &  2.46 &     0.07\% \\\hline
5 &  4 &  8 &  4 &  Column Generation &          1 &   550 &    0.03 &   12 &  12 &   2.82 &     2.80 &   0.95\% \\
& & & & Bi-partite B\&P &        4 &   550 &   0.049 &   18 &   18 &  2.82 &  2.80 &     0.89\% \\
& & & & Tri-partite B\&P &       70 &   550 &   0.63 &   96 &    100 &  2.81 &  2.81 &     0.07\% \\\bottomrule
5 &  6 &  4 &  2 &  Column Generation &          1 &    24 &    0.03 &    2.02 &   2.09 &   2.05 &     1.94 &   5.32\% \\
& & & & Bi-partite B\&P &      136 &    24 &    1.76 &   24 &   38 &  2.04 &  2.03 &     0.46\% \\
& & & & Tri-partite B\&P &      203 &    24 &   2.9 &     32 &    56 &  2.03 &  2.03 &     0.002\% \\\hline
5 &  6 &  4 &  4 &  Column Generation &          1 &   111 &    0.04 &    6.9 &   6.94 &   2.36 &     2.30 &   2.8\% \\
& & & & Bi-partite B\&P &      724 &   111 &  20 &  328 &  656 &  2.36 &  2.33 &     1.02\% \\
& & & & Tri-partite B\&P &      964 &   111 &  42 &  442 &   760 &  2.34 &  2.34 &     0.09\% \\\hline
5 &  6 &  6 &  2 &  Column Generation &          1 &    74 &    0.04 &    5.9 &   6.07 &   2.94 &     2.80 &   4.72\% \\
& & & & Bi-partite B\&P &       16 &    74 &    0.21 &   17 &    18 &  2.94 &  2.94 &     0.005\% \\
& & & & Tri-partite B\&P &       16 &    74 &   1.06 &   18 &    20 &  2.94 &  2.94 &     0.005\% \\\hline
5 &  6 &  6 &  4 &  Column Generation &          1 &   459 &    0.04 &   20 &  20 &   3.32 &     2.93 &  11.86\% \\
& & & & Bi-partite B\&P &      206 &   459 &   4.1 &  350 &  389 &  3.31 &  3.31 &     0.19\% \\
& & & & Tri-partite B\&P &      243 &   459 &   4.94 &  396 &   443 &  3.31 &  3.31 &     0.1\% \\\hline
5 &  6 &  8 &  2 &  Column Generation &          1 &   209 &    0.04 &   12 &  12 &   3.93 &     3.80 &   3.39\% \\
& & & & Bi-partite B\&P &      556 &   209 &  18 &  512 &  770 &  3.92 &  3.90 &     0.36\% \\
& & & & Tri-partite B\&P &      757 &   209 &   31 &  676 &  1,088 &  3.91 &  3.90 &     0.1\% \\\hline
5 &  6 &  8 &  4 &  Column Generation &          1 &  1,369 &    0.04 &   47 &  47 &   4.41 &     4.30 &   2.53\% \\
& & & & Bi-partite B\&P &       56 &  1,369 &   0.8 &  273 &  278 &  4.40 &  4.40 &     0.1\% \\
& & & & Tri-partite B\&P &       56 &  1,369 &   0.7 &  272 &   277 &  4.40 &  4.40 &     0.1\%\\
\bottomrule
\end{tabular}
\end{table}

\subsection{Robustness to parameter estimation errors}\label{app:robust}

Recall that we consider in this paper a \textit{deterministic} prescriptive contagion analytics model. We conduct in this appendix a sensitivity analysis to establish the robustness of our results under parameter misspecification. Specifically, we assess the performance of resource allocation decisions under a perturbed contagion model in Table~\ref{tab:vaccine_robust}, for the vaccine allocation problem. In these experiments, we consider all solutions obtained with each decision-making procedure under the original DELPHI-V model (Figure~\ref{subfig:vaccines}); we then perturb all DELPHI-V parameters according to uniform distributions by up to 20\% from their original values; and we evaluate the performance of each solution on the perturbed model. We replicate this process 20 times and report average performance, using the same nomenclature as in Table~\ref{tab:vaccine0}.

Results confirm the strong performance of our solution. Note that the number of lives saved increases slightly as compared to those estimated in Table~\ref{tab:vaccine0} with the nominal parameters, due to the non-linearities in the system dynamics. However, the relative performance of the optimized vaccine allocation solution is very similar to the one obtained under the nominal parameters. Specifically, the optimized vaccine allocation solution yields improvements of 42\% (resp. 12\%) with a supply of 2.5M vaccines (resp. 7M vaccines) per week as compared to the cost-based allocation benchmark, an improvements of 1.6\% (resp. 5.3\%) with a supply of 2.5M vaccines (resp. 7M vaccines) per week as compared to the coordinate descent benchmark. These results suggest that our solutions and our findings are highly robust against model misspecification.

\begin{table}[h!]
\footnotesize\renewcommand{\arraystretch}{1.0}
\caption{Average death toll comparison under $20\%$ perturbation of contagion parameters (vaccine allocation problem, full country, $D=21$, average over 20 samples).}
\label{tab:vaccine_robust}
\begin{center}
\begin{tabular}{clcccccc}
\toprule
&&\multicolumn{2}{c}{\textbf{$S=6$}}& \multicolumn{2}{c}{\textbf{$S=8$}} & \multicolumn{2}{c}{\textbf{$S=12$}} \\ \cmidrule(lr){3-4}\cmidrule(lr){5-6}\cmidrule(lr){7-8}
Budget & \multicolumn{1}{c}{Method} & Time (sec.) & Deaths & Time (sec.) & Deaths & Time (sec.) & Deaths \\\toprule 
| & Do nothing & | & 576.57K& | & 607.48K & | & 657.02K \\\toprule
2.5M& Uniform allocation & \multicolumn{1}{c}{|} & -6.34K& \multicolumn{1}{c}{|} & -11.96K & \multicolumn{1}{c}{|} & -23.90K \\
& Cost-based allocation & \multicolumn{1}{c}{|} & -7.28K& \multicolumn{1}{c}{|} &  -14.13K & \multicolumn{1}{c}{|} & -31.15K \\\cmidrule{2-8}
& MIQO implementation & n/a & n/a& n/a & n/a & n/a & n/a \\
& Discretization ($\delta=0.002$) & 1,000$^*$ & -5.47K& n/a & n/a & n/a & n/a\\
& Discretization ($\delta=0.001$) & 1,000$^*$ & -5.81K& n/a & n/a & n/a & n/a\\
& Coordinate Descent & 5.43 & \textbf{-12.90K}& 4.73 & -21.75K & 6.03 & -43.49K \\\cmidrule{2-8}
&Branch-and-price ($\varepsilon=0.002$) & 426.2 & \textbf{-12.98K}& 1234.2 & \textbf{-22.57K} & 3671.4 & \textbf{-44.17K}  \\\toprule
7M&Uniform allocation& \multicolumn{1}{c}{|} & -13.56K& \multicolumn{1}{c}{|} & -25.29K & \multicolumn{1}{c}{|} & -47.61K\\
&Cost-based allocation& \multicolumn{1}{c}{|} & -18.04K& \multicolumn{1}{c}{|} & -33.78K  & \multicolumn{1}{c}{|}& -63.60K \\\cmidrule{2-8}
&MIQO implementation& n/a & n/a& n/a & n/a & n/a & n/a \\ 
&Discretization ($\delta=0.002$) & 1,000$^*$ & -6.07K& n/a & n/a & n/a & n/a \\ 
&Discretization ($\delta=0.001$) & 1,000$^*$ & -3.18K& n/a & n/a & n/a & n/a \\ 
&Coordinate Descent  & 5.46 & -21.25K& 12.7 & -38.46K & 118.5 & -67.43K \\\cmidrule{2-8}
&Branch-and-price ($\varepsilon=0.002$) & 67.5 & \textbf{-22.03K}& 156.8 & \textbf{-39.73K} & 708.7 & \textbf{-70.99K} \\\bottomrule
\end{tabular}
\end{center}
\begin{tablenotes}
    \vspace{-6pt}
    \item $*$ and ``n/a'': no optimal and feasible solution, respectively. Bold font:  solutions within $1\%$ of the best-found solution.
\end{tablenotes}
\end{table}

\section{Benefits of \texorpdfstring{$\ell_\infty$}{}-based state clustering}
\label{app:clustering}

In this appendix, we provide a detailed comparison of our state clustering approach (referred to as $\ell_\infty$ clustering, for disambiguation) to a $k$-means benchmark, detailed in Algorithm~\ref{alg:dp_clust_kmeans}. Recall that the purpose of the state clustering algorithm is to group $Q$ states $\bN_{s+1}^1,\cdots,\bN_{s+1}^Q$ into $k$ clusters at each epoch $s=1,\cdots,S$. From a computational standpoint, our $\ell_\infty$ algorithm scales linearly to retain tractability given the large number of states (Algorithm~\ref{alg:dp_clust}), and so does $k$-means clustering. The key difference is that our $\ell_\infty$ clustering algorithm assigns states based on their $\ell_\infty$ distance to the element-wise minimum and maximum states in each cluster; in turn, it provides guarantees on the diameter of each cluster and on the global approximation error (Section~\ref{subsec:dynamic_program}). In comparison, the $k$-means benchmark minimizes the average distance to cluster centroids. Note, also, that the $k$-means benchmark requires a number of clusters as an input, whereas our $\ell_\infty$ clustering algorithm assigns states to clusters dynamically.

\begin{algorithm}[h!]
    \caption{$k$-means state-clustering}
    \label{alg:dp_clust_kmeans}
    \renewcommand{\arraystretch}{1.0}\footnotesize
    \begin{algorithmic}[1]
        \Procedure{kMeansStateClustering}{$\calN_1$, $k$}
            \For{$s \gets 1$ to $S$} \Comment{Forward pass}
            \State Initialization: $\calN_{s+1} \gets \emptyset$, $\widetilde{\calN}_{s+1} \gets \emptyset$, $k \gets 0$ 
            \For{$\bN_s\in\calN_s,\ \bx_s\in\calF_s$}
                \State Next state: $\bM(\tau_s)\gets\bN_s$; $\frac{d\bM(t)}{dt}=f(\bM(t),\bx_{s}),\ \forall t\in\left[\tau_s,\tau_{s+1}\right]$; $\bN_{s+1}\gets\bM(\tau_{s+1})$
                \State $\widetilde{\calN}_{s+1} \gets \widetilde{\calN}_{s+1}  \cup \{\bN_{s+1}\}$
                \State Compute cost: $\overline{C}_s^q=\int_{\tau_s}^{\tau_{s+1}}g_{t}(\bM^q(t))dt+\mathbf{1}\{s=S\}(h(\bM^q(T)))$
            \EndFor
            \State $\gamma_1,\cdots,\gamma_k \gets \text{kMeans}(\widetilde{\calN}_{s+1}, k)$ \Comment{$\gamma_1,\cdots, \gamma_k$ denote the $k$-means clusters}
            \State Store centroids $\bmu(\gamma_\ell)$, and binary parameters $\xi^q(\gamma_\ell)$ indicating if $\bN_{s+1}^q$ is assigned to cluster $\gamma_\ell$, $\ell=1,\cdots,k$
            \For{$l \gets 1$ to $k$}
                \State $\calN_{s+1} \gets \calN_{s+1} \cup \{\bmu(C_\ell)\}$
                \State Compute average cost $c_s^\mu(C_\ell) \gets \sum_{q=1}^Q\xi^q(\gamma_l)\overline{C}^q_s$
            \EndFor
            \For{$\bN_s\in\calN_s,\ \bx_s\in\calF_s$}
                \State Define cost function $c_s(\bN_s,\bx_s)  \gets \sum_{l=1}^k \xi^q(\gamma_l)c_s^\mu(C_\ell) +C_{s}(\bx_{s})$
            \EndFor
        \EndFor
        \EndProcedure
    \end{algorithmic}
\end{algorithm}

First, we compare the $\ell_\infty$ clustering algorithm to the $k$-means clustering in terms of approximation accuracy, for the vaccine allocation problem. Table~\ref{tab:DP_kmeans} reports the Median Absolute Error (MAE) and the Maximum Absolute Error between the clustered states and the true states obtained with exhaustive enumeration (Algorithm~\ref{alg:dp_full}). Note that we consider here small-scale examples for which exhaustive enumeration remains tractable. For a fair comparison, we select values of $k$ to closely match the number of states from $\ell_\infty$ clustering. As expected, $k$-means achieves a smaller median absolute error but the $\ell_\infty$ algorithm achieves a smaller maximum absolute error. For instance, in the six-period instance with the largest number of clusters, both algorithms achieve comparable MAE but the $k$-means benchmark results in a five-fold increase in the maximum absolute error; in the four-period instance, the $k$-means benchmark achieves a very small MAE but still increases the maximum absolute error by a factor 2 to 3. These results reflect that our $\ell_\infty$ algorithm is precisely designed to bound the diameter of its clusters (Proposition~\ref{prop:cluster}).

\begin{table}[h!]
\centering
\footnotesize\renewcommand{\arraystretch}{1.0}
\caption{$\ell_\infty$ clustering versus and $k$-means clustering (vaccine allocation problem).}
\label{tab:DP_kmeans}
\resizebox{\textwidth}{!}{
\begin{tabular}{cccllcccccccccccc}
\toprule
 & & & & & & & \multicolumn{5}{c}{\textbf{Maximum Absolute Error ($\times 10^{-2}$)}} & \multicolumn{5}{c}{\textbf{Median Absolute Error ($\times 10^{-4}$)}} \\ \cmidrule(lr){8-12}\cmidrule(lr){13-17}
$n$ & $S$ & $D$ & \multicolumn{1}{c}{Method} & \multicolumn{1}{c}{Hyperparameter} & $|S|$ & Time (sec.) & $s=2$ & $s=3$ & $s=4$ & $s=5$ & $s=6$ & $s=2$ & $s=3$ & $s=4$ & $s=5$ & $s=6$ \\\hline 
51 & 4 & 6 & $\ell_\infty$ Clust. & $\varepsilon=0.002$ &2,350 & 3.33 &0.12&0.19&0.38 & |  & |&0.00&0.29&0.27& |& |  \\
&& & $k$-Means & $k=12$ &2,193 & 3.41 &0.00&0.51&1.49 & |  & | & 0.00&0.03&0.02& |  & |  \\
&& & $\ell_\infty$ Clust. & $\varepsilon=0.005$ &1,458 & 2.20  &0.28&0.57&0.68& |  & | &0.00&1.04&0.95 & |  & |  \\
&& & $k$-Means & $k=8$ &1,581 & 2.81  &0.00&5.08&5.14 & |  & |&0.00&0.01&0.32 & |  & |  \\
&& & $\ell_\infty$ Clust. & $\varepsilon=0.01$ & 982 & 1.57  &0.68&0.87&0.86 & |  & |&0.01&3.72&4.02 & |  & | \\
&& & $k$-Means & $k=5$ &1,071 & 1.95  &4.96&9.22&14.15 & |  & |&0.00&0.95&0.97& |  & |  \\\hline
51 & 4 & 11 & $\ell_\infty$ Clust. & $\varepsilon=0.002$ &4,971 & 13.3 &0.12&0.27&0.49 & |  & |&0.08&0.37&0.36& |& |  \\
&& & $k$-Means & $k=25$ &4,437 & 12.4 &0.02&0.65&1.57 & |  & | & 0.0&0.01&0.12& |  & |  \\
&& & $\ell_\infty$ Clust. & $\varepsilon=0.005$ &2,533 & 6.74  &0.27&0.68&0.83& |  & | &0.29&1.59&1.42 & |  & |  \\
&& & $k$-Means & $k=13$ &2,601 & 9.34  &0.01&4.47&4.38 & |  & |&0.03&0.08&0.24 & |  & |  \\
&& & $\ell_\infty$ Clust. & $\varepsilon=0.01$ & 1,521 & 4.29 &0.72&0.86&1.37 & |  & |&0.29&3.58&4.23 & |  & | \\
&& & $k$-Means & $k=8$ &1,683 & 6.02 &3.75&6.23&9.79 & |  & |&0.01&0.48&0.83& |  & |  \\\hline
51 & 4 & 21 & $\ell_\infty$ Clust. & $\varepsilon=0.002$ &10,615 & 74.6 &0.08&0.27&0.67 & |  & |&0.08&0.52&0.41& |& |  \\
&& & $k$-Means & $k=53$ &9,231 & 59.3 &0.00&1.14&1.52 & |  & | & 0.00&0.00&0.01& |  & |  \\
&& & $\ell_\infty$ Clust. & $\varepsilon=0.005$ &4,477 & 26.9  &0.37&0.62&0.79& |  & | &0.01&1.86&1.57 & |  & |  \\
&& & $k$-Means & $k=22$ &4,488 & 27.4  &0.01&3.37&6.72 & |  & |&0.02&0.14&0.18 & |  & |  \\
&& & $\ell_\infty$ Clust. & $\varepsilon=0.01$ & 2,473 & 14.8  &0.82 &0.89&1.47 & |  & |&0.28&3.78&4.91 & |  & | \\
&& & $k$-Means & $k=13$ &2,703 & 18.5  &3.82&6.76&10.47 & |  & |&0.00&0.51&0.62& |  & |  \\\hline
51 & 6 & 6 & $\ell_\infty$ Clust. & $\varepsilon=0.002$ &5,555 & 10.9 &0.09&0.16&0.38 & 0.26  & 0.29&0.00&0.31&0.27& 0.07& 0.08  \\
&& & $k$-Means & $k=19$ &5,202 & 9.11 &0.00&0.32&1.58 & 1.66  & 1.79 & 0.00&0.00&0.02& 0.07  & 0.08 \\
&& & $\ell_\infty$ Clust. & $\varepsilon=0.005$ &2,903 & 5.51  &0.32&0.63&0.71& 0.84  & 0.87 &0.00&1.02&0.88 & 0.49  & 0.16 \\
&& & $k$-Means & $k=10$ &2,907 & 5.24  &0.00&0.87&1.86 & 3.04  & 3.42&0.00&0.01&0.22& 0.13  & 0.12  \\
&& & $\ell_\infty$ Clust. & $\varepsilon=0.01$ & 1,736 & 3.20  &0.73&0.92&0.94 & 1.47  & 2.07&0.00&3.65&4.13 & 2.16  & 1.57 \\
&& & $k$-Means & $k=6$ &1,887 & 3.94  &0.00&9.22&9.05 & 17.73  & 21.47&0.00&0.24&0.52& 0.47  & 0.36  \\
\bottomrule
\end{tabular}
}
\end{table}

Table \ref{tab:vaccine_kmeans} then reports the performance of the branch-and-price algorithm, with each state clustering algorithm in the pricing problem. Again, to perform an apples-to-apples comparison, the death toll is computed using the full continuous-state contagion models, as opposed to either state-clustering approximation. Note, at the outset, that both algorithms lead to similar computational times in branch-and-price, which is expected since both scale linearly with the number of states to cluster. However, the $\ell_\infty$ algorithm consistently leads to stronger solutions than the $k$-means algorithm, resulting in hundreds of extra lives saved over a six-week horizon. This result underscores the benefits of controlling the worst-case deviation between the clustered states and the true states via the $\ell_\infty$ algorithm---as noted theoretically in Proposition \ref{prop:cluster} and experimentally in Table \ref{tab:DP_kmeans}. The smaller cluster diameters mitigate error propagation over dynamical systems---as established in Proposition~\ref{prop:approx_error}---and can ultimately result in superior resource allocation decisions from the optimization methodology developed in this paper.

\begin{table}[h!]
\footnotesize\renewcommand{\arraystretch}{1.0}
\caption{Comparison of vaccine allocation solutions obtained with $\ell_\infty$ clustering vs. $k$-means clustering.}
\label{tab:vaccine_kmeans}
\begin{center}
\resizebox{\textwidth}{!}{
\begin{tabular}{llccccccccc}
\toprule
   &&\multicolumn{3}{c}{\textbf{$n=51$, $S=6$, $D=6$}} & \multicolumn{3}{c}{\textbf{$n=51$, $S=6$, $D=11$}} & \multicolumn{3}{c}{\textbf{$n=51$, $S=6$, $D=21$}} \\ \cmidrule(lr){3-5}\cmidrule(lr){6-8}\cmidrule(lr){9-11}
Hyperparameter & \multicolumn{1}{c}{Method} & $|S|$ & Time (sec.) & Deaths & $|S|$ & Time (sec.) & Deaths & $|S|$ & Time (sec.) & Deaths \\\hline 
|& Do nothing & | & | & 573.37K & | & | & 573.37K & | & | & 573.37K \\\hline
$\varepsilon=0.01$ & Branch-and-price ($k$-means)& 1,887 & 4.99 & -8.31K & 3,111& 13.4 & -10.57K & 4,641 & 32.6 & -10.95K  \\
& Branch-and-price ($\ell_\infty$)& 1,736 & 4.91 & \textbf{-10.37K} & 2,656 & 12.5 & \textbf{-10.94K} & 4,328 & 26.1 & \textbf{-11.10K} \\
& Extra lives saved & |& |&$+24.8\%$  & |& |&$+3.5\%$ & |& |&$+1.4\%$\\\hline
$\varepsilon=0.002$ & Branch-and-price ($k$-means)& 5,202 & 13.11 & -10.31K & 12,342 & 62.4 & -10.70K & 29,172 & 439.1 & -11.02K \\
 & Branch-and-price ($\ell_\infty$)& 5,555  & 12.61 & \textbf{-10.66K} & 12,834 & 65.4 & \textbf{-11.07K} & 29,020 & 426.2 & \textbf{-11.28K}  \\
& Extra lives saved & |& |&$+9.7\%$  & |& |&$+3.4\%$ & |& |&$+1.8\%$\\
\bottomrule
\end{tabular}
}
\end{center}
\begin{tablenotes}
    \vspace{-6pt}
    \item $*$ and ``n/a'': no optimal and feasible solution, respectively. Bold font:  solutions within $1\%$ of the best-found solution.
\end{tablenotes}
\end{table}
\end{APPENDICES}

\end{document}